\documentclass[12pt,reqno]{amsart}

\usepackage{amsfonts}
\usepackage{amsmath}
\usepackage{amssymb}
\usepackage{amsthm}
\usepackage{bm}
\usepackage{enumerate} 
\usepackage{hyperref}
\usepackage{cite}
\usepackage[margin=1.2in]{geometry}
\usepackage{mathrsfs} 
\usepackage{graphicx} 
\usepackage{tikz-cd}

\usepackage{xcolor}

\theoremstyle{plain}
\newtheorem{theorem}{Theorem}[section]
\newtheorem{lemma}[theorem]{Lemma}
\newtheorem{proposition}[theorem]{Proposition}
\newtheorem{corollary}[theorem]{Corollary}

\theoremstyle{definition}
\newtheorem{definition}[theorem]{Definition}

\newtheorem{remark}[theorem]{Remark}

\newtheorem{example}[theorem]{Example}

\numberwithin{equation}{section}

\newcommand\N{\mathbb{N}}
\newcommand\Z{\mathbb{Z}}
\newcommand\R{\mathbb{R}}

\newcommand\C{\mathbb{C}}

\newcommand\D{\mathcal{D}}

\renewcommand\S{\mathcal{S}}

\newcommand\ev[2]{\langle#1,#2\rangle}

\DeclareMathOperator\supp{supp}
\DeclareMathOperator\im{Im}

\DeclareMathOperator{\spn}{span}
\DeclareMathOperator\re{Re}

\DeclareMathOperator\QA{QA}

\DeclareMathOperator{\sinc}{sinc}

\begin{document}

\title[Uncertainty principles and Gelfand-Shilov spaces]{Optimal Beurling--H\"ormander-type uncertainty principles and nontriviality of Gelfand-Shilov spaces}

\author[L. Neyt]{Lenny Neyt}
\address{University of Vienna\\ Faculty of Mathematics\\ Oskar-Morgenstern-Platz 1 \\ 1090 Wien\\ Austria}
\thanks{The research of L. Neyt was funded in whole by the Austrian Science Fund (FWF) 10.55776/ESP8128624. For open access purposes, the author has applied a CC BY public copyright license to any author-accepted manuscript version arising from this submission.}
\email{lenny.neyt@univie.ac.at}


\begin{abstract}
We study how rapidly a nonzero function $f \in L^2(\R)$ and its Fourier transform $\widehat{f}$ can simultaneously decay. 
For each suitable weight function $\omega$, we construct a convex weight function $\varpi$ and establish the existence of a threshold $\Phi_{(\omega,\varpi^*)} \in  (0, \infty)$, where $\varpi^*$ denotes the convex conjugate of $\varpi$. 
For $a,b > 0$, the conditions
\[
    \int_{-\infty}^\infty |f(x)| e^{\omega(a|x|)} \, \mathrm{d}x < \infty
    \qquad \text{and} \qquad
    \int_{-\infty}^\infty |\widehat{f}(\xi)| e^{\varpi^*(b|\xi|)} \, \mathrm{d}\xi
    < \infty
\]
force $f = 0$ when $ab > \Phi_{(\omega,\varpi^*)}$, whereas the functions satisfying both conditions form a dense subspace of $L^2(\R)$ when $ab < \Phi_{(\omega,\varpi^*)}$. 
Similar results are also established in higher dimensions. 
We characterize the weights $\omega$ for which $\varpi$ is equivalent to $\omega$ and use this to determine when the Beurling--H\"ormander uncertainty principle is optimal. 
A discrete form of the same construction yields a characterization of the nontriviality of Gelfand-Shilov spaces $\S^{M}_{A}(\R^d)$ and $\Sigma^{M}_{A}(\R^d)$ in terms of their defining weight sequences $M$ and $A$, under suitable regularity assumptions.
\end{abstract}

\maketitle

\section{Introduction}
\label{sec:Intro}

\subsection{Uncertainty principles of Beurling--H\"ormander type}
\label{sec:BHIntro}

The uncertainty principle in harmonic analysis encompasses a family of results giving an exact mathematical interpretation of the following phenomenon:
A function $f \in L^2(\R^d)$ and its Fourier transform
	\[ \widehat{f}(\xi) = \int_{\R^d} f(x) e^{- i \ev{x}{\xi}}  \mathrm{d}x \]
cannot simultaneously be too well-localized.
Wiener first observed this type of behavior, and it has led to a fruitful and active branch of research; see, for example, the surveys \cite{F-S-Uncertainty, H-J-Uncertainty, T-Uncertainty}.
One of the first and most influential results in this direction was given by Hardy \cite{H-ThmFourier}.
He observed that a function and its Fourier transform cannot both have arbitrarily strong Gaussian decay.
In fact, the threshold is given by the Gaussian functions themselves. 
Shortly after, Morgan \cite{M-NoteFourier} refined this result to the nonsymmetric case.
Let $p \in (1, \infty)$ and $p^* = p / (p - 1)$ be its conjugate exponent.
Morgan showed the existence of a constant $C_p > 0$ such that if $\lambda^2 > C_p$, then, for any $f \in L^2(\R)$,
	\begin{equation}
		\label{eq:MorganUncertainty} 
		\int_{-\infty}^\infty |f(x)| e^{\lambda^{p} |x|^p} \, \mathrm{d}x < \infty \quad \text{and} \quad \int_{-\infty}^\infty |\widehat{f}(\xi)| e^{\lambda^{p^*} |\xi|^{p^*}} \, \mathrm{d}\xi < \infty \quad \Longrightarrow \quad f = 0 . 
	\end{equation}
When $\lambda^2 < C_p$, however, he showed the existence of nontrivial $f$ satisfying the specific time-frequency decay.

More recently, in \cite{S-ClassEntFuncRapDecRealAxis}, Sedletski\v{\i} extended the work of Hardy and Morgan to regularly varying functions \cite{B-G-T-RegVar, S-RegVar}.
We recall that an eventually positive function $\omega : [0, \infty) \to [0, \infty)$ is called \emph{regularly varying} when $\omega(t) = t^p L(t)$ for some $p \in \R$, called the \emph{index}, and with $L : [0, \infty) \to [0, \infty)$ a \emph{slowly varying function}, that is, $L(\lambda t) / L(t) \to 1$ for every $\lambda > 0$ when $t \to \infty$.
For $p > 1$, the conjugate of $\omega$ is asymptotic to $\eta(t) = t^{p^*} L^*(t)$, where $L^*(t) = \frac{p-1}{p^{p^*}}(L^{p^*-1})^{\#}(t^{p^*-1})$, with $\#$ denoting the de Bruijn conjugate (see \cite[Section 1.5.7, p.~28]{B-G-T-RegVar}). 
Sedletski\v{\i} showed that for some constant $C'_p > 0$, if $\lambda^2 > C'_p$, then, for any $f \in L^2(\R)$,
	\begin{equation}
		\label{eq:SedletskiiUncertainty}
		\int_{-\infty}^\infty |f(x)| e^{\lambda^p |x|^p L(\lambda |x|)} \, \mathrm{d}x < \infty ~ \text{and} ~ \int_{-\infty}^\infty |\widehat{f}(\xi)| e^{\lambda^{p^*} |\xi|^{p^*} L^*(\lambda |\xi|)} \, \mathrm{d}\xi < \infty  \Longrightarrow f = 0 . 
	\end{equation}
If $\lambda^2 < C'_p$, there exists a nontrivial $f$ satisfying the specific time-frequency decay.
Therefore, the uncertainty principles \eqref{eq:MorganUncertainty} and \eqref{eq:SedletskiiUncertainty} are, in essence, optimal.

Beurling proved the following elegant result \cite[Theorem on p.~372]{B-CollectedWorks2}:
if $f \in L^1(\R)$, then
	\begin{equation}
		\label{eq:BeurlingUncertainty} 
		\iint_{\R^2} |f(x)| |\widehat{f}(\xi)| e^{|\ev{x}{\xi}|} \, \mathrm{d}x \mathrm{d}\xi < \infty \quad \Longrightarrow \quad f = 0 . 
	\end{equation}
We refer to \cite{B-D-J-HermiteUncertaintyPrinc} for a multidimensional version of \eqref{eq:BeurlingUncertainty}, and to \cite{G-BeurlingUncertaintyPrinciple, H-HeisenbergUncertaintyPrincipleBeurling} for extensions of it.
A straightforward consequence goes as follows.
Let $\omega : [0, \infty) \to [0, \infty)$ be convex and let
	\[ \omega^*(t) = \sup_{s \geq 0} [ st - \omega(s)] , \qquad t \geq 0 , \]
be its convex conjugate. 
For any $f \in L^2(\R)$,
	\begin{equation}
		\label{eq:BeurlingHormanderUncertainty} 
		\int_{-\infty}^\infty |f(x)| e^{\omega(|x|)} \, \mathrm{d}x < \infty \quad \text{and} \quad \int_{-\infty}^\infty |\widehat{f}(\xi)| e^{\omega^*(|\xi|)} \, \mathrm{d}\xi < \infty \quad \Longrightarrow \quad f = 0 . 
	\end{equation}
In particular, the above result recovers \eqref{eq:MorganUncertainty} and \eqref{eq:SedletskiiUncertainty}, up to constants; the fact that $\omega^*(t) \sim t^{p^*} L^*(t)$ follows from \cite[Theorem 1.8.10, p.~48]{B-G-T-RegVar}.
Even though \eqref{eq:BeurlingUncertainty} is formulated in \cite{B-CollectedWorks2}, the editors state that no proof of it was preserved.
Afterwards, H\"ormander published a proof of \eqref{eq:BeurlingUncertainty} and \eqref{eq:BeurlingHormanderUncertainty} in \cite{H-UniquenessThmBeurlingFourier}.
H\"ormander's arguments were based on notes he took when Beurling explained the results to him during their time together at the Institute for Advanced Study.
H\"ormander also states that Beurling remarked that using a result of Lindel\"of, he could construct examples showing that \eqref{eq:BeurlingHormanderUncertainty} is rather precise.
The exact details on this, however, appear to have been lost, and the problem of optimality remained open.
One of the primary aims of this paper is to give a definitive answer to this question.

The apparent optimality of \eqref{eq:BeurlingHormanderUncertainty}, up to constants, is supported by the results of Hardy, Morgan, and Sedletski\v{\i}.
But, when $\omega(t)$ is dominated by any $t^p$ for $p \in (1, \infty)$ (or equivalently, when $\omega^*(t)$ dominates any $t^p$ for $p \in (1, \infty)$), this conclusion seems to break down.
In the limit case $\omega(t) = t$, we may interpret \eqref{eq:BeurlingHormanderUncertainty} as
	\begin{equation}
		\label{eq:FourierCompactSupport}
		\int_{-\infty}^\infty |f(x)| e^{\omega(|x|)} \, \mathrm{d}x < \infty \quad \text{and} \quad \widehat{f} \, \text{ has compact support} \quad \Longrightarrow \quad f =  0 . 
	\end{equation}
However, when $\omega(2t) = O(\omega(t))$, another classical result due to Beurling \cite{B-First} says that \eqref{eq:FourierCompactSupport} already holds whenever
	\begin{equation}
		\label{eq:QuasiAnalytic} 
		\int_1^\infty \frac{\omega(t)}{t^2} \, \mathrm{d}t = \infty . 
	\end{equation}
See also the work of Hirschman in \cite{H-BehFourierTransQuasiAnalytic}.
Using the Phragm\'en-Lindel\"of principle, one shows that, for any $\eta > \frac{\pi}{2}$,
	\begin{equation}
		\label{eq:MaxGrowthBound} 
		\int_{-\infty}^\infty |f(x)| e^{|x|} \, \mathrm{d}x < \infty \quad \text{and} \quad  |\widehat{f}(\xi)| \lesssim e^{-e^{\eta |\xi|}} \qquad \Longrightarrow \qquad f = 0 . 
	\end{equation}
If $\eta < \frac{\pi}{2}$, there exist nontrivial $f$ satisfying the time-frequency decay as in \eqref{eq:MaxGrowthBound}.
See Theorem \ref{t:MainUncertaintyPrincipleCaset} for a more precise statement of this fact.
Observing that $\omega^*(t) \asymp e^{t-1}$ for $\omega(t) = t \log_+ t$, with $\log_+ = \max\{0, \log\}$, we then see the appearance of a discrepancy in the optimality of \eqref{eq:BeurlingHormanderUncertainty}.

Our primary contribution is establishing a broad optimal uncertainty principle that in particular encompasses all of the above cases.
The principal idea is to use, on one side, the convex conjugate of a weight that measures the rate at which \eqref{eq:QuasiAnalytic} blows up.
Specifically, to $\omega$ we will associate the convex function $\varpi : [0, \infty) \to [0, \infty)$ defined as $\varpi(t) = 0$ for $t \in [0, 1]$ and
	\begin{equation}
		\label{eq:varpi} 
		\varpi(t) = t \int_1^t \frac{\omega(s)}{s^2} \, \mathrm{d}s , \qquad t > 1 .
	\end{equation}
When $\omega(t) = t$, we find exactly $\varpi(t) = t \log_+ t$, while for regularly varying $\omega(t) = t^p L(t)$, $p \in (1, \infty)$, we have $\varpi^* \asymp \omega^*$.
In fact, the functions $\omega$ and $\varpi$ are equivalent exactly when the lower Matuszewska index $\beta(\omega)$ \cite{B-G-T-RegVar} of $\omega$ is strictly larger than $1$; see Proposition \ref{p:SuffNeccCondEquivVarpi}.
Therefore, using $\varpi^*$ instead of $\omega^*$ in \eqref{eq:BeurlingHormanderUncertainty} recovers, up to constants, the optimal results \eqref{eq:MorganUncertainty}, \eqref{eq:SedletskiiUncertainty}, and \eqref{eq:MaxGrowthBound}.

Let us now specify in exact terms how we will express our results.
We will interpret an (optimal) uncertainty principle as in the following definition; see also Definition \ref{def:UncertaintyPairMultiDim} for the multidimensional version.

\begin{definition}
\label{def:UncertaintyPair}
Let $\omega, \eta : [0, \infty) \to [0, \infty)$ be nondecreasing.
We call $(\omega, \eta)$ an \emph{uncertainty pair} if there exists a constant $\Phi_{(\omega, \eta)} \in [0, \infty)$ such that for any $a, b > 0$ the following hold:
	\begin{itemize}
		\item[(I)] If $a b > \Phi_{(\omega, \eta)}$, should $f \in L^2(\R)$ satisfy
			\begin{equation}
				\label{eq:UncertaintyDecayIntro}
		 		\int_{-\infty}^\infty |f(x)| e^{\omega(a |x|)} \, \mathrm{d}x < \infty \quad \text{and} \quad \int_{-\infty}^\infty |\widehat{f}(\xi)| e^{\eta(b |\xi|)} \, \mathrm{d}\xi < \infty 
			\end{equation}
		then $f = 0$.
		\item[(II)] If $a b < \Phi_{(\omega, \eta)}$, there exists a dense subspace of $L^2(\R)$ consisting of functions $f$ satisfying \eqref{eq:UncertaintyDecayIntro}.
	\end{itemize}
We call $\Phi_{(\omega, \eta)}$ the \emph{uncertainty threshold} of $(\omega, \eta)$.
The uncertainty pair $(\omega, \eta)$ is called \emph{optimal} if $\Phi_{(\omega, \eta)} > 0$.
\end{definition}

Clearly, $(\omega, \eta)$ is an (optimal) uncertainty pair if and only if $(\eta, \omega)$ is an (optimal) uncertainty pair, with the same uncertainty threshold.

One of our main results is the following optimal uncertainty principle.
It will be shown in Section \ref{sec:UncertaintyPrinciples}.
We note that it also holds in multiple dimensions, as in Definition \ref{def:UncertaintyPairMultiDim}.
Here and below, we will call an eventually nonnegative function $\omega$ on $[0, \infty)$ \emph{moderate} if
	\[ \limsup_{t \to \infty} \frac{\omega(2t)}{\omega(t)} < \infty  \quad \text{and} \quad \liminf_{t \to \infty} \frac{\omega(2t)}{\omega(t)} > 1 . \]

\begin{theorem}
\label{t:MainUncertaintyPrincipleIntro}
Let $\omega : [0, \infty) \to [0, \infty)$ be an unbounded nondecreasing function that is moderate, such that $t \mapsto \omega(e^t)$ is convex, and for which \eqref{eq:QuasiAnalytic} holds.
Let $\varpi$ be defined as in \eqref{eq:varpi}.
Then $(\omega, \varpi^*)$ is an optimal uncertainty pair.
\end{theorem}

By Lemma \ref{l:RegVar>1Convex}, the previous result covers the optimality of \eqref{eq:SedletskiiUncertainty}.
In general, the exact value of the uncertainty threshold $\Phi_{(\omega, \varpi^*)}$ does not follow from our arguments.
But in Theorem \ref{t:UncertaintyThresholdBound} we provide explicit upper and lower bounds in terms of the weight function.
When $\omega(t) = t^p L(t)$ is convex and regularly varying of index $p \geq 1$, and we put $\rho = \min\{p, p^*\}$ when $p > 1$, the exact value of the threshold is:
	\begin{equation}
		\label{eq:UncertaintyThresholdRegVar} 
		\Phi_{(\omega, \varpi^*)} = \begin{cases} \frac{\pi}{2} , & p = 1 , \\ \frac{\left[\sin\left(\frac{\pi}{2}(\rho - 1)\right)\right]^{1 /\rho}}{(p-1)^{1/p}} , & p > 1 , \end{cases} 
	\end{equation}
see Theorem \ref{t:LowBoundUncertaintyThresholdRegVar}.
Note that for $p > 1$ this follows from the results of Sedletski\v{\i}.
The exact behavior at the uncertainty threshold $\Phi_{(\omega, \varpi^*)}$ was considered by Morgan for $\omega(t) = t^p$ in \cite{M-NoteFourier}, while in the more general case it is unclear.
When $\omega(t) = t^2$, the Gaussian functions display the optimal symmetric time-frequency decay.
We also refer to the recent work of the author with Toft and Vindas in \cite{N-T-V-HermiteExpansions, N-T-V-QuantitativeUncertaintyPrinciples} concerning quantitative uncertainty principles for time-frequency Gaussian decay.

An answer to the optimality question of \eqref{eq:BeurlingHormanderUncertainty} can now be summarized in the specific case of regularly varying functions as follows.
Its proof will also be given in Section \ref{sec:UncertaintyPrinciples}; see also Theorem \ref{t:SuffNeccCondBeurlingHormanderOptimal} for a more general version.
Additionally, the case of $\omega(t) = t$ is covered in Theorem \ref{t:MainUncertaintyPrincipleCaset} (see also \eqref{eq:LowerBoundWeightFunct}).

\begin{theorem}
\label{t:MainUncertaintyPrincipleRegVar}
Let $\omega : [0, \infty) \to [0, \infty)$ be an unbounded nondecreasing convex function that is regularly varying and $t = o(\omega(t))$.
Then $(\omega, \omega^*)$ is an uncertainty pair.
Moreover
	\begin{equation*}
		\label{eq:CharSharpness}
		(\omega, \omega^*) \text{ is optimal} \quad \iff \quad \text{the index of $\omega$ is strictly larger than $1$} . 
	\end{equation*}
\end{theorem}

\medskip

\subsection{The nontriviality of Gelfand-Shilov spaces}
\label{sec:GSIntro}

In his seminal work \cite{S-ThDist}, Schwartz introduced the Fr\'echet space $\S(\R^d)$ consisting of all smooth functions $f$ satisfying
	\begin{equation}
		\label{eq:Schwartz} 
		\sup_{x \in \R^d} |x^\alpha f^{(\beta)}(x)| < \infty , \qquad \forall \alpha, \beta \in \N^d . 
	\end{equation}
As the Fourier transform restricts to an isomorphism on $\S(\R^d)$, this allows one to uniquely extend it by transposition to the strong dual $\S'(\R^d)$, the space of tempered distributions.
This procedure allowed Schwartz to rigorously define the Fourier transform on a large class of generalized functions.
However, in their studies of parabolic initial-value problems \cite{G-S-First, G-S-GenFunc3, G-S-GenFunc2}, Gel'fand and Shilov introduced broader spaces than $\S'(\R^d)$ on which the Fourier transform can be uniquely extended.
These spaces then turned out to be the right framework for the analysis of decay and regularity properties of global solutions to large classes of linear and semi-linear PDE on $\R^d$; see, for example, \cite{C-R-TF, G-GS, N-R-GlobalPseudoDiffCalcEuclidean}.

The idea of Gel'fand and Shilov is to quantify the bounds in \eqref{eq:Schwartz} in terms of $\alpha$ and $\beta$, thereby making the test function space smaller than $\S(\R^d)$ and broadening its dual to include elements beyond tempered distributions. 
In this paper, a \emph{weight sequence} $M = (M_p)_{p \in \N}$ stands for a sequence of positive real numbers that is logarithmically convex and is eventually greater than $h^p$ for any $h > 0$; see Definition \ref{def:WeightSeq}.
Then, given two weight sequences $M$ and $A$, we consider those smooth functions $f$ satisfying, for some $h, k > 0$,
	\begin{equation}
		\label{eq:GelfandShilov}
		\sup_{x \in \R^d} |x^\alpha f^{(\beta)}(x)| \lesssim k^\alpha A_{|\alpha|} \cdot h^\beta M_{|\beta|} , \qquad \alpha, \beta \in \N^d . 
	\end{equation}	
Then, we may define the spaces
	\begin{align}
		\label{eq:GSRoumieuDef}
		\S^M_A(\R^d) &= \{ f \in C^\infty(\R^d) \mid f \text{ satisfies \eqref{eq:GelfandShilov} for some $h, k > 0$} \} \\
	\intertext{and}
		\label{eq:GSBeurlingDef}
		\Sigma^M_A(\R^d) &= \{f \in C^\infty(\R^d) \mid f \text{ satisfies \eqref{eq:GelfandShilov} for all $h, k > 0$} \} .
	\end{align}
The spaces $\S^M_A(\R^d)$ and $\Sigma^M_A(\R^d)$ are called the \emph{Gelfand-Shilov spaces (of Roumieu and Beurling type)}.
Under mild conditions on $M$ and $A$, the Fourier transform defines a bijection between $\S^M_A(\R^d)$ and $\S^A_M(\R^d)$ and between $\Sigma^M_A(\R^d)$ and $\Sigma^A_M(\R^d)$; see Lemma \ref{l:GSFourierInvariant}.

The structural and topological properties of the Gelfand-Shilov spaces are well-understood; see, for instance, \cite{C-C-K-CharGSSpFourier, D-N-V-NuclGSSpaces, D-N-V-InclGSSp, P-P-V-Parametrix}.
But an immediate question is now: When are $\S^M_A(\R^d)$ and $\Sigma^M_A(\R^d)$ nontrivial, that is, when do they contain nonzero elements?
This is a longstanding open question that has only received partial answers.
By the Denjoy-Carleman theorem \cite[Theorem 4.2, p.~56]{K-UltradistI}, both $\S^M_A(\R^d)$ and $\Sigma^M_A(\R^d)$ can contain compactly supported nontrivial elements if and only if	
	\begin{equation}
		\label{eq:DenjoyCarleman}
		\sum_{p = 1}^\infty \frac{M_{p-1}}{M_p} < \infty .
	\end{equation}
In this case, no additional conditions on $A$ are needed, while \eqref{eq:DenjoyCarleman} implies that $p! \lesssim \varepsilon^p M_p$ for any $\varepsilon > 0$ \cite[Lemma 4.1, p.~55]{K-UltradistI}.
Hirschman gave a technical characterization in terms of $A$ for nontriviality in \cite{H-BehFourierTransQuasiAnalytic} when \eqref{eq:DenjoyCarleman} fails but still $p! \lesssim h^p M_p$ for some $h > 0$.
In the special case when $M = (p!)_{p \in \N}$---henceforth denoted as $p!$---Debrouwere and Vindas proved in \cite{D-V-NonTrivSpAnalytic} that $\S^{p!}_A(\R^d)$, respectively $\Sigma^{p!}_A(\R^d)$, is nontrivial if and only if for some, respectively for all, $k > 0$,
	\begin{equation}
		\label{eq:Mp=p!}
		\sup_{p \geq 2} \frac{(\log p)^p}{k^p A_p}  < \infty . 
	\end{equation}
Therefore, the situation for $M$ above and up to $p!$ is well-understood.

When $M_p \lesssim \varepsilon^p p!$ for any $\varepsilon > 0$, any $f$ satisfying \eqref{eq:GelfandShilov} extends to an entire function on $\C^d$ and the situation changes drastically.
Here, the question of nontriviality of the Gelfand-Shilov spaces is most interesting, yet primarily unresolved.
An elegant characterization is available for what is often referred to as the \emph{classical} Gelfand-Shilov spaces $\S^r_s(\R^d)$ and $\Sigma^r_s(\R^d)$.
In this case, the weight sequences are chosen as $M = (p!^{r})_{p \in \N}$ and $A = (p!^{s})_{p \in \N}$ for $r,s > 0$.
Then, Gel'fand and Shilov showed \cite{G-S-First, G-S-GenFunc2}:
	\begin{equation}
		\label{eq:GSClassicalNonTriv}
		\S^r_s(\R^d) \neq \{0\} \quad \Longleftrightarrow \quad r + s \geq 1 , \qquad \Sigma^r_s(\R^d) \neq \{0\} \quad \Longleftrightarrow \quad r + s > 1. 
	\end{equation}
Katznelson and Mandelbrojt later gave sufficient conditions in \cite{K-M-GS1, K-M-GS2} for nontriviality that in particular recover \eqref{eq:GSClassicalNonTriv}. 
Additionally, they state that ``from a certain point of view'' their results are optimal, though no decisive argument for this is given.
In fact, compared to the results of this paper, this claim is not true.

The following necessary condition for the nontriviality of $\S^M_A(\R^d)$, respectively $\Sigma^M_A(\R^d)$, is well known:
	\begin{equation}
		\label{eq:GSNonTrivNecc}
		p! \lesssim h^p M_p A_p , \qquad p \in \N ,
	\end{equation}
for some, respectively all, $h > 0$.
See Proposition \ref{p:NonTrivGSNeccCond} for a short proof of this fact in our cases.
For the classical Gelfand-Shilov spaces, this condition is also sufficient, as follows from \eqref{eq:GSClassicalNonTriv}.
However, considering the weight sequence $M_p = p! / (\log_+ p)^p$, $p\geq2$, indicates that \eqref{eq:GSNonTrivNecc} cannot serve as a general characterization for nontriviality.
Indeed, here the condition translates to \eqref{eq:Mp=p!}, and thus the nontriviality of $\S^{p!}_A(\R^d)$ and $\Sigma^{p!}_A(\R^d)$.

The idea behind our characterization of nontriviality is to replace $M$ or $A$ in \eqref{eq:GSNonTrivNecc} with altered weight sequences $\overline{M}$ and $\overline{A}$.
To do so, we write $m_p = M_p / M_{p-1}$, $p \geq 1$, and we let $m(t)$ denote the counting function of this sequence.
Then $\overline{M}$ is the weight sequence for which
	\begin{equation}
		\label{eq:overlineMIntro} 
		\sum_{k = 1}^\infty m\left(\frac{t}{k}\right) 
	\end{equation}
is the counting function of the sequence $\overline{m}_p = \overline{M}_p / \overline{M}_{p-1}$.
See Definition \ref{def:overlineM} for exact details.
When $M = p!$, we find $\overline{M}$ is equivalent to $(p! / (\log_+ p)^p)_{p \in \N}$, while for $M = (p!^r)_{p \in \N}$, $0<r<1$, the weight sequence $\overline{M}$ is equivalent to $M$.
In fact, we have $M \asymp \overline{M}$ exactly when the upper Matuszewska index $\alpha(M)$ \cite{JG-S-S-ORegVar} of $M$ is strictly smaller than $1$; see Proposition \ref{p:CondEquivSeq}.

We will prove the following characterization of nontriviality of the Gelfand-Shilov spaces.
See also the version presented in Theorem \ref{t:GeneralCharNonTrivGS}, where we relax the conditions.
Below, the classical condition $(M.1)^*$ for a weight sequence $N$ means that $(p! / N_p)_{p \in \N}$ is also a weight sequence;
$(M.2)$ says that $N_{p+q} \leq C H^{p+q} N_p N_q$ for certain $C, H > 0$ and all $p, q \in \N$;
while $(M.2)^*$ means that eventually $2 n_p \leq n_{Qp}$ for some $Q \in \Z_+$, where $n_p = N_p / N_{p-1}$, $p \in \Z_+$.
For two weight sequences $N, K$ we write $N \preceq K$, respectively $N \triangleleft K$, if $N_p \lesssim h^p K_p$ for some, respectively all, $h > 0$.

\begin{theorem}
\label{t:CharNonTrivGSIntro}
Let $M$ and $A$ be weight sequences with one satisfying $(M.1)^*$ and the other $(M.2)$ and $(M.2)^*$.
Then 
	\begin{equation}
		\label{eq:NiceCharNonTrivRoumieu} 
		\S^{M}_{A}(\R^d) \neq \{0\} \quad \Longleftrightarrow \quad p! \preceq \overline{M} \cdot A ~ \text{ and } ~ p! \preceq M \cdot \overline{A} 
	\end{equation}
and
	\begin{equation}
		\label{eq:NiceCharNonTrivBeurling} 
		\Sigma^{M}_{A}(\R^d) \neq \{0\} \quad \Longleftrightarrow \quad p! \triangleleft \overline{M} \cdot A ~ \text{ and } ~ p! \triangleleft M \cdot \overline{A} .
	\end{equation}
\end{theorem}

When both sequences satisfy $(M.2)$ and $(M.2)^*$, then \eqref{eq:GSNonTrivNecc} for some, respectively all, $h > 0$ is indeed the characterizing condition for nontriviality of $\S^M_A(\R^d)$, respectively $\Sigma^M_A(\R^d)$.
This is stated in the following result; see also Theorem \ref{t:GSNonTrivCharClassicalCond} for a more general result.
Otherwise, the condition may become insufficient; see Theorem \ref{t:UncertaintyGSOriginalM}.

\begin{theorem}
\label{t:CharNonTrivGSIntro2}
Let $M$ and $A$ be weight sequences satisfying $(M.2)$ and $(M.2)^*$.
Then 
	\[
		\S^{M}_{A}(\R^d) \neq \{0\} \quad \Longleftrightarrow \quad p! \preceq M \cdot A
	\]
and
	\[
		\Sigma^{M}_{A}(\R^d) \neq \{0\} \quad \Longleftrightarrow \quad p! \triangleleft M \cdot A . 
	\]
\end{theorem}

\medskip

\subsection{Proof strategy and structure of this paper}

We now give an outline of how we will prove our main results.

There is an inherent connection between the two problems above, which was already observed in the work of Gel'fand and Shilov \cite{G-S-GenFunc3, G-S-GenFunc2}.
See van Eijndhoven \cite{E-FuncAnalCharGS} and Chung, Chung, and Kim \cite{C-C-K-CharGSSpFourier} for full characterizations.
Moreover, both types of behaviour can be equivalently expressed through bounds on entire functions.
If, for a function $f$, its Fourier transform $\widehat{f}$ decays faster than any exponential, by the Fourier inversion formula, the original function $f$ extends to an entire function.
Then, the specific decay of $\widehat{f}$ controls the growth of this extension in the imaginary direction.
When $f$ also has specified decay on the real line, say $f$ satisfies \eqref{eq:UncertaintyDecayIntro}, we are led to entire functions $F$ satisfying Paley--Wiener-type bounds
	\begin{equation}
		\label{eq:EntireFuncBoundsIntro} 
		|F(x + iy)| \lesssim \exp\left[ \eta^*\left(\frac{|y|}{b}\right) - \omega(a |x|) \right] , \qquad z = x+iy \in \C . 
	\end{equation}
This will be made more precise in Section \ref{sec:UncertaintyPrinciples}; see, in particular, Lemma \ref{l:TFDecayAndEntireFunc}.
A similar observation can be made for elements in the Gelfand-Shilov spaces and the entire extension of their Fourier transform.
Here, we may associate a weight function to a weight sequence $M$ as
	\[ \omega_M(t) = \sup_{p \in \N} \log \frac{t^p M_0}{M_p} . \]
See Section \ref{sec:OptimalDecayRegularity}, and in particular Lemma \ref{l:GSToEntire}, for the exact details.

Therefore, both problems are equivalent to the question of existence of entire functions satisfying \eqref{eq:EntireFuncBoundsIntro}. 
This is a classical question that goes back to the work of Gel'fand and Shilov \cite{G-S-First, G-S-GenFunc3, G-S-GenFunc2}, and was further considered by H\"ormander \cite{H-PaleyWiener}, Babenko \cite{B-NewProblemQuasiAnal, B-ClassesSpInfDiff}, Katznelson and Mandelbrojt \cite{K-M-GS1, K-M-GS2}, and Evgrafov \cite{E-AsympEstEntFunc}, among others.
In \cite{E-AsympEstEntFunc}, it was shown that for \eqref{eq:EntireFuncBoundsIntro} to hold for a nontrivial entire $F$, we must necessarily have
	\[ \limsup_{t \to \infty} \frac{\varpi(t)}{\eta^*(t)} < \infty . \]
In Proposition \ref{p:NonTrivEntireLowerBound} we establish a more quantitative version of this.
One of the main contributions of this paper is to show that it is then, in fact, also sufficient for the general class of weight functions that we consider; see Theorem \ref{t:CharNonTrivEntireWeightFunc}.
To show this, in Propositions \ref{p:ConstructionFunction} and \ref{p:ConstructWindow}, we construct an infinite product of dilates of a single entire function.
We establish \eqref{eq:EntireFuncBoundsIntro}, with $\eta^* = \varpi$, by verifying decay on the real line and a growth bound uniform on horizontal lines.
Then, a uniform bound can be found through interpolation estimates essentially due to Hirschman \cite{H-ConvexTh}; see Lemma \ref{l:InterpolaEst}.
This method also allows us to give a general lower bound for the threshold when these bounds can hold; see Theorem \ref{t:PaleyWienerThresholdLowerBound}.
In the special case where $\omega$ is a convex regularly varying function of index at least $1$, we establish the exact values in Theorem \ref{t:LowBoundPaleyWienerThresholdRegVar}.

Let us now review the structure of this paper.
It is split into two parts: In Sections \ref{sec:Weights}, \ref{sec:EntireFunc}, and \ref{sec:UncertaintyPrinciples} we work with weight functions; in Sections \ref{sec:WeightSeq}, \ref{sec:OptimalDecayRegularity}, and \ref{sec:NonTrivGS} we deal with weight sequences.
First, in Section \ref{sec:Weights}, we formally introduce the weight functions $\omega$ we work with and the associated weight function $\varpi$ as in \eqref{eq:varpi}. 
Our main result here is the characterization of when $\omega$ and $\varpi$ are equivalent (Proposition \ref{p:SuffNeccCondEquivVarpi}), or when $\varpi$ dominates $\omega$ (Proposition \ref{p:CharCondDomVarpi}).
Then, in Section \ref{sec:EntireFunc}, we discuss the optimality of the Paley--Wiener-type bounds as explained above.
Section \ref{sec:UncertaintyPrinciples} then passes from entire functions to Fourier pairs.
Here, we establish the results as discussed in Section \ref{sec:BHIntro}.

Sections \ref{sec:WeightSeq}--\ref{sec:NonTrivGS} transfer these results to weight sequences. 
Lemma \ref{l:overlineM} shows that the associated function $\omega_{\overline{M}}$ of $\overline{M}$ is equivalent to the weight obtained by applying \eqref{eq:varpi} to $\omega_M$.
Further in Section \ref{sec:WeightSeq}, we characterize when $M$ and $\overline{M}$ are equivalent or when $\overline{M}$ is dominated by $M$; see Proposition \ref{p:CondEquivSeq}.
Section \ref{sec:OptimalDecayRegularity} considers the optimal bounds for \eqref{eq:GelfandShilov}.
In particular, we see that this is obtained when $A_p = \overline{M}_p^* = p! / \overline{M}_p$; see Theorem \ref{t:UncertaintyGS}.
See also Theorem \ref{t:UncertaintyGSOriginalM} for when $\overline{M}_p^*$ may be replaced by $M_p^*$.
Finally, in Section \ref{sec:NonTrivGS}, we establish the results as discussed in Section \ref{sec:GSIntro}.

\medskip

\section{Weight functions}
\label{sec:Weights}

We fix, in exact terms, how we will express our results through weight functions.
Moreover, we consider relevant transforms and review several examples.

\subsection{Weight functions}
\label{sec:WeightFunctions}

For two nondecreasing functions $\omega, \eta : [0, \infty) \to [0, \infty)$ we write $\omega \preceq \eta$ if, for some $C, H > 0$,
	\[ \omega(t) \leq \eta(H t) + C , \qquad t \geq 0 . \]
We define the constant
	\[ \mathfrak{c}(\omega, \eta) = \inf\{ H \in [0, \infty] \mid \exists C > 0, ~ \forall t > 0 : \omega(t) \leq \eta(H t) + C \} . \]
We call $\omega$ and $\eta$ \emph{equivalent} if $\omega \preceq \eta$ and $\eta \preceq \omega$.
When $\mathfrak{c}(\omega, \eta) = 0$, we say $\eta$ \emph{dominates} $\omega$ and we write $\omega \triangleleft \eta$.

\begin{definition}
\label{def:WeightFunc}
A \emph{weight function} stands for a nondecreasing unbounded continuous function $\omega : [0, \infty) \to [0, \infty)$ for which:
	\begin{itemize}
		\item[$(i)$] $\omega(0) = 0$;
		\item[$(ii)$] $\log t = o(\omega(t))$;
		\item[$(iii)$] $t \mapsto \omega(e^t)$ is convex on $\R$.
	\end{itemize}
We call $\omega$ \emph{quasianalytic} if \eqref{eq:QuasiAnalytic} holds, that is, if
	\[ \int_1^\infty \frac{\omega(t)}{t^2} \, \mathrm{d}t = \infty . \]
\end{definition}

Let $\omega$ be a weight function.
As the function $\varphi(t) = \omega(e^t)$ is nondecreasing and convex on $\R$, it has a nondecreasing right derivative $\varphi_+'(t)$ everywhere on $\R$.
Let $\omega_+'$ denote the right derivative of $\omega$ on $(0, \infty)$.
It follows that
	\[ \varphi_+'(\log t) = t \omega_+'(t) , \qquad t > 0 . \]
We may therefore introduce the following function.
	
\begin{definition}
Let $\omega$ be a weight function.
We define $\mu = \mu_\omega$ as the function $[0, \infty) \to [0, \infty)$ with $\mu(0) = 0$ and 
	\[ \mu(t) = t \omega_+'(t) , \qquad t > 0 , \]
where $\omega_+'$ denotes the right derivative of $\omega$.
\end{definition}

For a weight function $\omega$, the function $\mu$ is nondecreasing and
	\begin{equation}
		\label{eq:omegaAsIntegral} 
		\omega(t) = \int_0^t \frac{\mu(s)}{s} \, \mathrm{d}s , \qquad t \geq 0 . 
	\end{equation}
Then, for any $\lambda > 1$, we find the bound
	\begin{equation}
		\label{eq:muBoundomega} 
		\mu(t) \leq \frac{\omega(\lambda t) - \omega(t)}{\log \lambda} \leq \mu(\lambda t) .
	\end{equation}
When $\omega$ is convex, then $\omega'_+(t) \geq \omega(t) / t$, so automatically
	\begin{equation}
		\label{eq:OmegaBoundMuConvex}
		\omega(t) \leq \mu(t) . 
	\end{equation}

We will often assume a weight function satisfies one of the following two conditions.

\begin{definition}
Let $\omega$ be a weight function.
	\begin{itemize}
		\item[$(i)$] $\omega$ is called \emph{lower moderate} if $2 \omega(t) \leq \omega(H t) + C$ for certain $C, H > 0$.
		\item[$(ii)$] $\omega$ is called \emph{upper moderate} if $\omega(2t) \leq H \omega(t) + C$ for certain $C, H > 0$.
	\end{itemize}
When $\omega$ is both lower and upper moderate, we call $\omega$ \emph{moderate}.
\end{definition}

Each of these conditions is invariant under equivalence.
A lower moderate weight function $\omega$ necessarily dominates some polynomial, that is, $t^\sigma \triangleleft \omega$ for some $\sigma > 0$.
On the other hand, an upper moderate weight function $\omega$ is necessarily polynomially dominated, that is, $\omega \triangleleft t^\sigma$ for some $\sigma > 0$.

We may now associate indices corresponding to the respective behavior.
We use the convention that $\omega(\infty) = \infty$.

\begin{definition}
Let $\omega$ be a weight function.
For any $\lambda > 0$, we set
	\[ \underline{\mathfrak{m}}_\lambda(\omega) = \inf \{ H \in (0, \infty] \mid \exists C > 0 ~ \forall t \geq 0 : \lambda \omega(t) \leq \omega(H t) + C \} \]
and
	\[ \overline{\mathfrak{m}}_\lambda(\omega) = \inf \{ H \in (0, \infty] \mid \exists C > 0 ~ \forall t \geq 0 : \omega(\lambda t) \leq H \omega(t) + C \} . \]
\end{definition}

For a weight function $\omega$, we have both $\underline{\mathfrak{m}}_\lambda(\omega) \leq 1$ and $\overline{\mathfrak{m}}_\lambda(\omega) \leq 1$ for any $\lambda \leq 1$, and $\underline{\mathfrak{m}}_1(\omega) = \overline{\mathfrak{m}}_1(\omega) = 1$.
Now $\omega$ is lower, respectively upper, moderate if and only if $\underline{\mathfrak{m}}_\lambda(\omega) < \infty$, respectively $\overline{\mathfrak{m}}_\lambda(\omega) < \infty$, for some/all $\lambda > 1$.
In this case, the respective coefficients depend continuously on $\lambda$.

\begin{lemma}
\label{l:Continuitymathfrakm}
Let $\omega$ be a weight function.
If $\omega$ is lower, respectively upper, moderate then the function $\lambda \mapsto \underline{\mathfrak{m}}_\lambda(\omega)$, respectively $\lambda \mapsto \overline{\mathfrak{m}}_\lambda(\omega)$, is nondecreasing, submultiplicative, and continuous on $(0, \infty)$.
\end{lemma}

\begin{proof}
A straightforward argument shows that the functions are nondecreasing and submultiplicative.
Therefore, to demonstrate the continuity of $\lambda \mapsto \underline{\mathfrak{m}}_\lambda(\omega)$, respectively $\lambda \mapsto \overline{\mathfrak{m}}_\lambda(\omega)$, it suffices to verify $\lim_{\lambda \to 1} \underline{\mathfrak{m}}_\lambda(\omega) = 1$, respectively $\lim_{\lambda \to 1} \overline{\mathfrak{m}}_\lambda(\omega) = 1$.

We start by verifying $\lim_{\lambda \to 1} \underline{\mathfrak{m}}_\lambda(\omega) = 1$ when $\omega$ is lower moderate.
Suppose $C_0, H_0 > 1$ are such that $2 \omega(t) \leq \omega(H_0 t) + C_0$.
As $\mu$ is nondecreasing
	\[ \mu(t) \geq \frac{1}{\log H_0} \int_{t / H_0}^t \frac{\mu(s)}{s} \, \mathrm{d}s = \frac{1}{\log H_0}\left(\omega(t) - \omega(t / H_0)\right) \geq \frac{1}{2 \log H_0} \omega(t) - \frac{C_0}{2 \log H_0} . \]
Set $c_0 = 1 / (2 \log H_0)$ and $c_1 = C_0 / (2 \log H_0)$.
For any $\varepsilon > 0$, by \eqref{eq:muBoundomega},
	\[ \omega(e^\varepsilon t) - \omega(t) \geq \varepsilon \mu(t) \geq \varepsilon c_0 \omega(t) - \varepsilon c_1 , \]
that is
	\[ (1 + \varepsilon c_0) \omega(t) \leq \omega(e^\varepsilon t) + \varepsilon c_1 . \]
Therefore $\underline{\mathfrak{m}}_{1 + \varepsilon c_0}(\omega) \leq e^\varepsilon$.
This shows that $\lim_{\lambda \to 1^+} \underline{\mathfrak{m}}_\lambda(\omega) = 1$.
On the other hand, since $(\underline{\mathfrak{m}}_{1/\lambda}(\omega))^{-1} \leq \underline{\mathfrak{m}}_\lambda(\omega)$ for $0 < \lambda < 1$, the result here follows.

We now show $\lim_{\lambda \to 1} \overline{\mathfrak{m}}_\lambda(\omega) = 1$ when $\omega$ is upper moderate.
We remark that this was also shown in \cite[Lemma 2.7]{N-T-V-QuantitativeUncertaintyPrinciples}, but we give an alternative argument.
Assume $C_1, H_1 > 1$ are such that $\omega(et) \leq H_1 \omega(t) + C_1$.
Using that $\mu$ is nondecreasing, we find
	\[ \mu(t) \leq \int_t^{e t} \frac{\mu(s)}{s} \, \mathrm{d}s = \omega(e t) - \omega(t) \leq (H_1 - 1) \omega(t) + C_1 . \]
For any $\varepsilon > 0$, by \eqref{eq:muBoundomega},
	\[ \omega(e^\varepsilon t) \leq \varepsilon \mu(e^\varepsilon t) + \omega(t) \leq \varepsilon (H_1 - 1) \omega(e^\varepsilon t) + \varepsilon C_1 + \omega(t) . \]
Suppose $\varepsilon$ is small enough so that $\varepsilon (H_1 -  1) < 1$.
Then
	\[ \omega(e^\varepsilon t) \leq \frac{1}{1 -  \varepsilon(H_1 - 1)} \omega(t) + \frac{\varepsilon C_1}{1 - \varepsilon(H_1 - 1)} ,  \]
so
	\[ \overline{\mathfrak{m}}_{e^\varepsilon}(\omega) \leq \frac{1}{1 - \varepsilon(H_1 - 1)} . \]
Therefore $\lim_{\lambda \to 1^+} \overline{\mathfrak{m}}_\lambda(\omega) = 1$.
As $(\overline{\mathfrak{m}}_{1/\lambda}(\omega))^{-1} \leq \overline{\mathfrak{m}}_\lambda(\omega)$ for $0 < \lambda < 1$, the proof is complete.
\end{proof}

We will often assume a weight function $\omega$ to be convex.
Note that then $\underline{\mathfrak{m}}_\lambda(\omega) \leq \lambda$ and $\overline{\mathfrak{m}}_{1/\lambda}(\omega) \leq 1 / \lambda$ for $\lambda \geq 1$.
In particular, $\omega$ is lower moderate.
Also, $\omega$ satisfies
	\begin{equation}
		\label{eq:LowerBoundWeightFunct} 
		\omega(t) \asymp t \quad \text{or} \quad t = o(\omega(t)) . 
	\end{equation}
For any weight function $\omega$, when $t = o(\omega(t))$, we define the \emph{convex conjugate} of $\omega$ as
	\[ \omega^*(t) = \sup_{s \geq 0} [st - \omega(s)] , \qquad t \geq 0 . \]
Then $\omega^*$ is a convex weight function.
It holds that $(\omega^*)^* = \omega$ exactly when $\omega$ is convex.
For two weight functions $\omega$ and $\eta$ dominating $t$: if $\omega \asymp \eta$ then also $\omega^* \asymp \eta^*$, while if $\omega \triangleleft \eta$ then $\eta^* \triangleleft \omega^*$. 

We now consider a typical family of examples of moderate weight functions.

\begin{definition}
\label{def:RegVar}
An eventually positive, locally bounded, measurable function $\eta : [0, \infty) \to [0, \infty)$ is called \emph{regularly varying} \cite{B-G-T-RegVar} if
	\[ \lim_{t \to \infty} \frac{\eta(\lambda t)}{\eta(t)} \in (0, \infty) , \qquad \forall \lambda > 0 . \]
One shows \cite[Theorem 1.4.1, p.~17]{B-G-T-RegVar} that $\eta(t) = t^p L(t)$ for some $p \in \R$, called the \emph{index} of $\eta$, and $L$ a slowly varying function, that is,
	\[ \lim_{t \to \infty} \frac{L(\lambda t)}{L(t)} = 1 , \qquad \forall \lambda > 0 . \]
\end{definition}
	
A regularly varying weight function $\omega$ of index $p > 0$ is always moderate, and we have
	\[ \underline{\mathfrak{m}}_\lambda(\omega) = \lambda^{\frac{1}{p}} \qquad \text{and} \qquad \overline{\mathfrak{m}}_\lambda(\omega) = \lambda^p , \qquad \lambda > 0 . \]
The weight function can be convex only if $p \geq 1$.
When a regularly varying weight function has index $p > 1$, it is always equivalent to a convex weight function.

\begin{lemma}
\label{l:RegVar>1Convex}
Let $\omega : [0, \infty) \to [0, \infty)$ be regularly varying of index $p > 1$.
Then there exists a regularly varying convex weight function $\eta$ such that $\mathfrak{c}(\omega, \eta) = \mathfrak{c}(\eta, \omega) = 1$.
\end{lemma}

\begin{proof}
By the Smooth Variation Theorem \cite[Theorem 1.8.2, p.~45]{B-G-T-RegVar} there exists a regularly varying nondecreasing smooth function $\widetilde{\eta} : [0, \infty) \to [0, \infty)$ such that $\omega(t) \sim \widetilde{\eta}(t)$ and for which
	\[ \lim_{t \to \infty} \frac{t^2 \widetilde{\eta}^{(2)}(t)}{\widetilde{\eta}(t)} = p(p - 1) > 0 . \]
Therefore, $\widetilde{\eta}$ is eventually convex.
Altering $\widetilde{\eta}$ on a finite interval, we may find our $\eta$.
\end{proof}

We now consider some well-known examples of regularly varying weight functions.

\begin{example}
\label{ex:ConvexConjugates}
$(i)$ The function asymptotically defined by
	\[ \omega_{p, \ell}(t) = t^p \log^\ell_+ t , \] 
is a moderate weight function if $p > 1$ and $\ell \in \R$ or if $ p = 1$ and $\ell \geq 0$.
When $p > 1$, we find
	\[ \omega^*_{p, \ell}(t) \asymp \omega_{p^*, - \frac{\ell}{p - 1}}(t) = t^{p^*} (\log_+ t)^{-\frac{\ell}{p - 1}} . \]
If $p = 1$ and $\ell > 0$, then
	\[ \omega^*_{1, \ell}(t) \asymp t^{\frac{\ell - 1}{\ell}} \exp\left[ t^{1 / \ell} \right] . \]
	
$(ii)$ Let $\log_{+, n}$, $n \in \Z_+$, denote the $n$th iteration of $\log_+$.
The functions
	\[ \omega_n(t) = t \log_{+, n} t , \qquad n \in \Z_+ , \]
are weight functions.
Let $\exp_n$, $n \in \Z_+$, denote the $n$th iteration of $\exp$.
Then
	\[ \omega^*_n(t) \asymp \frac{\exp_n t}{\prod_{j = 1}^{n-1} \exp_j t} , \qquad n \in \Z_+ . \]
\end{example}

We formally introduce $\varpi$ as in \eqref{eq:varpi} from the Introduction \ref{sec:Intro}.

\begin{definition}
Let $\omega : [0, \infty) \to [0, \infty)$ be a nondecreasing function. Then we define the function $\varpi = \varpi_\omega$ as $\varpi(t) = 0$ for $t \in [0, 1)$ and
	\[ \varpi(t) = t \int_1^t \frac{\omega(s)}{s^2} \, \mathrm{d}s  , \qquad t \geq 1 . \]
\end{definition}

Several properties of $\varpi$ are gathered in the following lemma.

\begin{lemma}
\label{l:PropertiesVarpi}
Let $\omega, \eta$ be weight functions.
	\begin{itemize}
		\item[$(i)$] $\varpi_\omega$ is a convex weight function.
		\item[$(ii)$] $t = o(\varpi_\omega(t))$ if and only if $\omega$ is quasianalytic.
		\item[$(iii)$] If $\omega \preceq \eta$, then also $\varpi_\omega \preceq \varpi_\eta$.
		In particular, if $\omega \asymp \eta$ then $\varpi_\omega \asymp \varpi_\eta$.
		\item[$(iv)$] We have $\omega \preceq \varpi$ and $\mathfrak{c}(\omega, \varpi) \leq 2$.
		If $\omega$ is lower moderate, then
			\begin{equation}
				\label{eq:CoeffomegatovarpiUpperBound}
				\mathfrak{c}(\omega, \varpi) \leq \inf_{\lambda > 0} \left[ \frac{\lambda + 1}{\lambda} \underline{\mathfrak{m}}_{\lambda}(\omega) \right] . 
			\end{equation}
		\item[$(v)$] If $\omega$ is convex, then
			\begin{equation} 
				\label{eq:UpperLowerBoundsVarpi2}
				\varpi(t) \leq \omega(t) \log t , \qquad t \geq 1 . 
			\end{equation}
		\item[$(vi)$] If $\omega$ is upper moderate, then $\varpi$ is moderate.
	\end{itemize}
\end{lemma}

\begin{proof}
$(i)$ and $(ii)$ are clear from the definition.

$(iii)$  Suppose $\omega(t) \leq \eta(H t) + C$ for some $C, H > 0$.
Then, for $t \geq 1$,
	\begin{equation}
		\label{eq:EstOfComparisonCoeffvarpi} 
		\varpi_\omega(t)  \leq t \int_{1}^t \frac{\eta(Hs) + C}{s^2} \, \mathrm{d}s \leq \varpi_\eta(H t) + \left(C - H \int_1^H \frac{\eta(s)}{s^2} \, \mathrm{d}s\right)t . 
	\end{equation}
Using $t = O(\varpi_\eta(t))$ and the convexity of $\varpi_\eta$, we conclude that $\varpi_\omega \preceq \varpi_\eta$.

$(iv)$
For any $\kappa > 1$ and $t \geq \kappa$,
	\[ \varpi(t) \geq t \int_{t/\kappa}^t \frac{\omega(u)}{u^2} \, \mathrm{d}u \geq (\kappa-1) \omega\left(\frac{t}{\kappa}\right) . \]
Putting $\kappa = 2$ gives $\omega(t) \leq \varpi(2t)$ for $t \geq 2$, so that $\omega \preceq \varpi$ and $\mathfrak{c}(\omega, \varpi) \leq 2$.
Suppose now that $\omega$ is lower moderate.
Fix $\kappa > 1$ and suppose $(\kappa-1)^{-1} \omega(t) \leq \omega(Ht) + C$ for certain $C, H > 0$.
Then, for $t \geq \kappa/H$,
	\[ \omega(t) \leq (\kappa - 1) \omega(H t) + (\kappa - 1) C \leq \varpi(\kappa H t) + (\kappa - 1) C . \]
This shows \eqref{eq:CoeffomegatovarpiUpperBound}.
	
$(v)$
Suppose $\omega$ is convex.
Then, for $t \geq 1$,
	\[ \varpi(t) = t \int_1^t \frac{\omega\left(\frac{s}{t} t\right)}{s^2} \, \mathrm{d}s \leq t \int_1^t \frac{\frac{s}{t} \omega(t)}{s^2} \, \mathrm{d}s = \omega(t) \log t . \]

$(vi)$ 
Fix some $\lambda > 0$.
Suppose $\omega(\lambda t) \leq H \omega(t) + C$ for some $C, H > 0$. 
Then, for $t \geq 1/\lambda$,
	\begin{align}
		\varpi(\lambda t) 
		&= \lambda t \int_1^{\lambda t} \frac{\omega(s)}{s^2} \, \mathrm{d}s 
		= t \int_{1/\lambda}^t \frac{\omega(\lambda s)}{s^2} \, \mathrm{d}s 
		\leq t \int_{1/\lambda}^t \frac{H \omega(s) + C}{s^2} \, \mathrm{d}s 
		\nonumber \\
		&\leq H \varpi(t) + \left( C \lambda - H \int_1^{1/\lambda} \frac{\omega(s)}{s^2} \, \mathrm{d}s \right) t . 
		\label{eq:EstOfvarpiAtLambda}
	\end{align}
Using $t = O(\varpi(t))$, we conclude that $\varpi$ is upper moderate.
As $\varpi$ is convex, it is also lower moderate, and thus moderate.
\end{proof}

\begin{lemma}
\label{l:IndicesForVarpi}
Let $\omega, \eta$ be quasianalytic weight functions.
Then, for any $\lambda > 0$,
	\[ \underline{\mathfrak{m}}_\lambda(\varpi) \leq \underline{\mathfrak{m}}_\lambda(\omega) \quad \text{and} \quad \overline{\mathfrak{m}}_\lambda(\varpi) \leq \overline{\mathfrak{m}}_\lambda(\omega) , \]
and if $\omega \preceq \eta$, then
	\[ \mathfrak{c}(\varpi_\omega, \varpi_\eta) \leq \mathfrak{c}(\omega, \eta) . \]
\end{lemma}

\begin{proof}
Fix $\lambda > 0$.
Suppose first $C, H > 0$ are such that $\lambda \omega(t) \leq \omega(H t) + C$.
For any $\varepsilon > 0$, for $t$ large enough,
	\[ \lambda \varpi(t) \leq t \int_1^t \frac{\omega(H s) + C}{s^2} \, \mathrm{d}s \leq (1 + \varepsilon) \varpi( H t) . \]
Letting $\varepsilon \to 0^+$ and using Lemma \ref{l:Continuitymathfrakm}, we see that $\underline{\mathfrak{m}}_\lambda(\varpi) \leq \underline{\mathfrak{m}}_\lambda(\omega)$.
Similar arguments, using \eqref{eq:EstOfComparisonCoeffvarpi} and \eqref{eq:EstOfvarpiAtLambda}, yield $\mathfrak{c}(\varpi_\omega, \varpi_\eta) \leq \mathfrak{c}(\omega, \eta)$ and $\mathfrak{m}_\lambda(\varpi) \leq \mathfrak{m}_\lambda(\omega)$.
\end{proof}

We now characterize when $\omega \asymp \varpi$ or when $\varpi$ dominates $\omega$.
First, we introduce the Matuszewska indices associated with a weight function.

\begin{definition}
\label{def:Matuszewska}
Let $\omega$ be a weight function. The Matuszewska indices \cite[Section 2.1.2, p.~68]{B-G-T-RegVar} of $\omega$ are defined as:
	\begin{itemize}
		\item[$(i)$] the \emph{upper Matuszewska index $\alpha(\omega)$} is given by
			\[ \alpha(\omega) = \inf \{ \alpha \in \R \mid \exists C, t_0 > 0 \text{ s.t. } \omega(\lambda t) \leq C \lambda^\alpha \omega(t) \text{ for all } \lambda \geq 1 \text{ and } t \geq t_0 \} . \]
	
		\item[$(ii)$] the \emph{lower Matuszewska index $\beta(\omega)$} is given by
			\[ \beta(\omega) = \sup \{ \beta \in \R \mid \exists C, t_0 > 0 \text{ s.t. } \omega(\lambda t) \geq C \lambda^\beta \omega(t) \text{ for all } \lambda \geq 1 \text{ and } t \geq t_0 \} . \]
	\end{itemize}
\end{definition}

For a weight function $\omega$, clearly always $\beta(\omega) \leq \alpha(\omega)$.
When $\omega$ is regularly varying of index $\rho$, then $\alpha(\omega) = \beta(\omega) = \rho$.
If $\eta$ is another weight function equivalent to $\omega$, then the Matuszewska indices of $\omega$ and $\eta$ are identical.

Whether a weight function is upper, respectively lower, moderate is determined by the upper, respectively lower, Matuszewska index.

\begin{lemma}
Let $\omega$ be a weight function.
The following statements are true:
	\begin{itemize}
		\item[$(i)$] $\omega$ is upper moderate if and only if $\alpha(\omega) < \infty$;
		\item[$(ii)$] $\omega$ is lower moderate if and only if $\beta(\omega) > 0$.
	\end{itemize}
\end{lemma}

\begin{proof}
See \cite[Theorem 2.1.7, p.~71]{B-G-T-RegVar}.
\end{proof}

We now see that $\omega \asymp \varpi$ is equivalent to $\beta(\omega) > 1$.

\begin{proposition}
\label{p:SuffNeccCondEquivVarpi}
Let $\omega$ be a weight function.
Then $\omega \asymp \varpi$ if and only if $\beta(\omega) > 1$.
In this case, for any $\varepsilon$,
	\begin{equation}
		\label{eq:UpperBoundInclusionCoeffvarpitoomega} 
		\mathfrak{c}(\varpi, \omega) \leq \underline{\mathfrak{m}}_{B + \varepsilon}(\omega) ,  \qquad  B = \inf_{1 < \beta < \beta(\omega)} \left[ \limsup_{t \to \infty} \frac{1}{\beta - 1} \left( \sup_{\lambda \geq 1} \frac{\lambda^\beta \omega(t)}{\omega(\lambda t)} \right) \right] .
	\end{equation}
When $\omega$ is convex,  we also have
	\begin{equation}
		\label{eq:LowerBoundInclusionCoeffvarpitoomega}
		\mathfrak{c}(\varpi, \omega) \geq \sup_{\lambda > 1} \frac{\underline{m}_{\lambda \log \lambda}(\omega)}{\lambda} . 
	\end{equation}
\end{proposition}

\begin{proof}
Suppose first that $\beta(\omega) > 1$.
Then $t = o(\omega(t))$ and $\omega$ is lower moderate.
Take $1 < \beta < \beta(\omega)$ and suppose $B' > 0$ is such that $\lambda^\beta \omega(t) \leq B' \omega(\lambda t)$ for all $\lambda \geq 1$ and $t \geq t_0 \geq 1$.
Fix $\varepsilon > 0$.
Let $t_1 \geq t_0$ be such that $\omega(t_0) t \leq \varepsilon \omega(t)$ for any $t \geq t_1$.
Finally, let $C_0, H_0 > 0$ be such that $(\varepsilon + B' / (\beta - 1)) \omega(t) \leq \omega(H_0 t) + C_0$.
Then, for $t \geq t_1$,
	\begin{align*} 
		\varpi(t) 
		&= t \int_1^t \frac{\omega(s)}{s^2} \, \mathrm{d}s 
		\leq t \int_1^{t_0} \frac{\omega(s)}{s^2} \, \mathrm{d}s + t \int_{t_0}^t \frac{\omega(s)}{s^2} \, \mathrm{d}s \\
		&\leq \omega(t_0) t + B' t \int_{t_0}^t \left(\frac{s}{t}\right)^{\beta} \frac{\omega(t)}{s^2} \, \mathrm{d}s
		\leq \left[\varepsilon + \frac{B'}{\beta - 1}\right]\omega(t)
		\leq \omega\left(H_0 t \right) + C_0 .
	\end{align*}
This shows that $\varpi \preceq \omega$, and using Lemma \ref{l:Continuitymathfrakm} we also see the bound in \eqref{eq:UpperBoundInclusionCoeffvarpitoomega}.
By Lemma \ref{l:PropertiesVarpi}$(iv)$, we conclude that $\omega$ and $\varpi$ are equivalent.

Conversely, suppose $\varpi \asymp \omega$.
Since the lower Matuszewska index is invariant under equivalence, and by Lemma \ref{l:PropertiesVarpi}$(iii)$, we may assume $\omega$ is convex.
Suppose $\varpi(t) \leq \omega(H_1 t) + C_1$ for some $C_1, H_1 > 0$.
As $\omega$ is convex, the function $t \mapsto \omega(t) / t$ is nondecreasing.
Then, for any $\lambda > 1$ and $t \geq 1$,
	\[
		\frac{\omega(t)}{t} 
		\leq \frac{1}{\log \lambda} \int_{t}^{\lambda t} \frac{\omega(s)}{s^2} \, \mathrm{d}s 
		\leq \frac{1}{\lambda \log \lambda} \frac{\varpi(\lambda t)}{t} 
		\leq \frac{1}{\lambda \log \lambda} \frac{\omega(\lambda H_1 t)}{t} + \frac{C_1}{\lambda \log \lambda} . 
	\]
As $\omega \asymp \varpi$, necessarily $\omega$ is quasianalytic, and then $t = o(\varpi(t))$ by Lemma \ref{l:PropertiesVarpi}$(ii)$.
Using the convexity of $\omega$, we get that for every $\lambda > 1$ and $\varepsilon > 0$, if $t$ is large enough,
	\begin{equation}
		\label{eq:alogatimesomega}
		(\lambda \log \lambda) \omega(t) \leq \omega((1 + \varepsilon) \lambda H_1 t) .
	\end{equation}
Using Lemma \ref{l:Continuitymathfrakm}, this already shows \eqref{eq:LowerBoundInclusionCoeffvarpitoomega}.
Put now $\kappa = \frac32 H_1 e^{2 H_1}$, which we may now assume to be strictly larger than $1$.
Set $\delta = \frac12 \log_\kappa (4/3)$.
Choosing $\lambda = e^{2 H_1}$ in \eqref{eq:alogatimesomega}, we get, for $t$ large enough,
	\[
		\kappa^{1 + \delta} \omega(t) < \omega(\kappa t) . 
	\]
By \cite[Theorem 2.16]{JG-S-S-ORegVar} we have that $\beta(\omega) > 1 + \delta$.
\end{proof}

To determine when $\varpi$ dominates $\omega$, we introduce the following behavior.

\begin{definition}
\label{def:GeomPers}
Let $\omega$ be a weight function.
Then $\omega$ is called \emph{geometrically persistent} if, for some $\lambda > 1$,
	\begin{equation}
		\label{eq:GeometricPersistent} 
		\lim_{t \to \infty} \sum_{k = 1}^{\lfloor \log_\lambda t \rfloor} \frac{\lambda^k \omega\left(\frac{t}{\lambda^k}\right)}{\omega(t)} = \infty . 
	\end{equation}
\end{definition}

Let $\omega$ and $\eta$ be equivalent weight functions such that $t = o(\omega(t))$ and $t = o(\eta(t))$.
Then $\omega$ is geometrically persistent if and only if $\eta$ is so.

\begin{proposition}
\label{p:CharCondDomVarpi}
Let $\omega$ be a weight function.
The following are equivalent: 
	\begin{itemize}
		\item[(a)] $\omega \triangleleft \varpi$;
		\item[(b)] $\omega(t) = o(\varpi(t))$;
		\item[(c)] $\varpi$ is regularly varying of index 1;
		\item[(d)] $\omega$ is geometrically persistent.
	\end{itemize}	
\end{proposition}

\begin{proof}
(a) $\Longrightarrow$ (d). 
Using the convexity of $\varpi$, for any $R > 1$ one finds a $t_R > 0$ such that $\varpi(t) / \omega(t) \geq R$ for $t \geq t_R$.
Fix any $\lambda > 1$.
We calculate, for $t \geq R$,
	\begin{align*}
		\varpi(t)
		&\leq t \sum_{k = 1}^{\lceil \log_\lambda t \rceil} \int_{t / \lambda^k}^{t / \lambda^{k-1}} \frac{\omega(s)}{s^2} \, \mathrm{d}s
		\leq t \sum_{k = 1}^{\lceil \log_\lambda t \rceil} \omega\left(\frac{t}{\lambda^{k-1}}\right) \left(\frac{\lambda^k}{t} - \frac{\lambda^{k-1}}{t}\right) \\
		&\leq (\lambda - 1) \left( \omega(t) + \sum_{k=1}^{\lfloor \log_\lambda t \rfloor} \lambda^k \omega\left(\frac{t}{\lambda^k}\right) \right) .
	\end{align*}
Consequently, for $t \geq \max\{t_R, R\}$,
	\[ \frac{R}{\lambda - 1} - 1 \leq \sum_{k = 1}^{\lfloor \log_\lambda t \rfloor} \frac{\lambda^k \omega\left(\frac{t}{\lambda^k}\right)}{\omega(t)}  . \]
Letting $R \to \infty$ shows that $\omega$ is geometrically persistent.

(d) $\Longrightarrow$ (b).
Let $\lambda > 1$ be such that \eqref{eq:GeometricPersistent} holds.
Then
	\begin{align*}
		\varpi(t)
		&\geq t \sum_{k = 1}^{\lfloor \log_\lambda t \rfloor} \int_{t / \lambda^k}^{t / \lambda^{k-1}} \frac{\omega(s)}{s^2} \, \mathrm{d}s
		\geq t \sum_{k = 1}^{\lfloor \log_\lambda t \rfloor}  \omega\left(\frac{t}{\lambda^k}\right) \left(\frac{\lambda^k}{t} - \frac{\lambda^{k-1}}{t}\right) \\
		&= \left(1 - \frac{1}{\lambda}\right) \sum_{k = 1}^{\lfloor \log_\lambda t \rfloor} \lambda^k \omega\left(\frac{t}{\lambda^k}\right) .
	\end{align*}
Therefore $\omega(t) = o(\varpi(t))$.

(b) $\Longrightarrow$ (c).
Put $L(t) = \varpi(t) / t$, then
	\[ \frac{t L'(t)}{L(t)} = \frac{\omega(t)}{\varpi(t)} = o(1) . \]
For any $\lambda > 1$, we now get
	\begin{equation*}
		\label{eq:varpiRegVarIndex1} 
		\log \frac{L(\lambda t)}{L(t)} = \int_t^{\lambda t} \frac{s L'(s)}{L(s)} \frac{\mathrm{d}s}{s} \to 0 , \qquad t \to \infty . 
	\end{equation*}
Thus $L$ is slowly varying, and since $\varpi(t) = t L(t)$, it follows that $\varpi$ is regularly varying of index 1.

(c) $\Longrightarrow$ (a).
For any $\lambda, t > 1$,
	\[ \frac{\varpi(\lambda t)}{\lambda t} - \frac{\varpi(t)}{t} = \int_{t}^{\lambda t} \frac{\omega(s)}{s} \frac{\mathrm{d}s}{s} \geq \left( 1 - \frac{1}{\lambda} \right) \frac{\omega(t)}{t} . \]
Therefore
	\[ \frac{\omega(t)}{\varpi(t)} \leq \frac{\lambda}{\lambda - 1} \left( \frac{\varpi(\lambda t)}{\lambda \varpi(t)} - 1 \right) = o(1) , \qquad t \to \infty , \]
that is, $\omega(t) = o(\varpi(t))$.
Now, for any $R \geq 1$, when $t$ is large enough,
	\[ \varpi(t) \geq \frac{1}{2R} \varpi(R t) \geq \omega(R t) . \]
We conclude $\omega \triangleleft \varpi$.
\end{proof}

\begin{corollary}
If a weight function $\omega$ is geometrically persistent, then so is $\varpi$.
\end{corollary}

\begin{proof}
By Proposition \ref{p:CharCondDomVarpi}, we have that $\varpi$ is a regularly varying function of index 1.
A straightforward argument shows that any such function is geometrically persistent.
\end{proof}

\begin{corollary}
\label{c:PolynomialBounds}
Let $\omega$ be a weight function.
	\begin{itemize}
		\item[$(i)$] If $\omega \asymp \varpi$ then $t^\beta = O(\omega(t))$ for some $\beta > 1$.
		\item[$(ii)$] If $\omega \triangleleft \varpi$ then $\omega(t) = o(t^{1 + \varepsilon})$ for every $\varepsilon > 0$.
	\end{itemize}
\end{corollary}

\begin{proof}
$(i)$ By Proposition \ref{p:SuffNeccCondEquivVarpi} we may take some $1 < \beta < \beta(\omega)$.
Then there exists some $C > 0$ such that $\omega(\lambda s) \geq C \lambda^\beta \omega(s)$ for any $\lambda \geq 1$ and $s \geq t_0$.
Therefore
	\[ \omega(t) \geq \frac{C \omega(t_0)}{t_0^{\beta}} t^\beta , \qquad t \geq t_0 . \]
	
$(ii)$ Proposition \ref{p:CharCondDomVarpi} shows $\varpi$ is regularly varying of index 1.
Combining $\omega \triangleleft \varpi$ with Potter's estimate \cite[Theorem 1.5.6, p.~25]{B-G-T-RegVar} yields the result.
\end{proof}

For regularly varying weight functions, one finds the following result.

\begin{lemma}
\label{l:RegVarEquivOverlineOmega}
Let $\omega$ be a regularly varying weight function of index $p \geq 1$.
The following statements are true.
	\begin{itemize}
		\item[1.] $\beta(\omega) > 1$ if and only if $p > 1$. In this case, $\mathfrak{c}(\omega, \varpi) = [p - 1]^{1/p}$ and $\mathfrak{c}(\varpi, \omega) = [p - 1]^{-1/p}$.
		\item[2.] $\omega$ is geometrically persistent if and only if $p = 1$.
	\end{itemize}
\end{lemma}

\begin{proof}
1. It is clear that $\beta(\omega) = p$. By \cite[Proposition 1.5.8, p.~26]{B-G-T-RegVar} we have $\varpi(t) \sim \omega(t) / [p - 1]$, from which the second statement follows.

2. If $p = 1$, then $\omega$ is clearly geometrically persistent. The converse follows from Proposition \ref{p:CharCondDomVarpi} and Corollary \ref{c:PolynomialBounds}.
\end{proof}

Let us now consider $\varpi$ for some examples.

\begin{example}
$(i)$ Let $\omega_{p, \ell}(t) = t^p \log_+^\ell t$ be as in Example \ref{ex:ConvexConjugates}$(i)$.
We write $\varpi_{p, \ell} = \varpi_{\omega_{p, \ell}}$.
If $p > 1$ and $\ell \in \R$, then by Proposition \ref{p:SuffNeccCondEquivVarpi} and Lemma \ref{l:RegVarEquivOverlineOmega} we have that
	\[ \varpi_{p, \ell} \asymp \omega_{p, \ell} \qquad \text{and} \qquad \varpi^*_{p, \ell} \asymp \omega_{\frac{p}{p - 1}, - \frac{\ell}{p - 1}} . \]
For $p = 1$ and $\ell \geq 0$, one verifies
	\[ \varpi_{1, \ell} \asymp \omega_{1, \ell + 1} \qquad \text{and} \qquad \varpi^*_{1, \ell} \asymp t^{\frac{\ell}{\ell + 1}} \exp\left[ t^{\frac{1}{\ell + 1}} \right]  . \]
	
$(ii)$ Let $\omega_{n}(t) = t \log_{+, n}(t)$ be as in Example \ref{ex:ConvexConjugates}$(ii)$.
We write $\varpi_n = \varpi_{\omega_n}$.
One calculates
	\[ \varpi_{n + 1}(t) \asymp \omega_{n+1}(t) \cdot \log_+ t, \qquad n \in \Z_+ . \]
Let $W_n(t)$ be such that $W_n(t) \exp_n W_n(t) = t$ for $n \in \Z_+$.
When $n = 1$, then $W_1$ is exactly the Lambert $W$-function.
We find
	\[ \varpi^*_{n + 1}(t) \asymp W_n(t) \exp\left[\frac{t}{W_n(t)}\right] , \qquad n \in \Z_+ .  \]
	
$(iii)$ Consider the weight function asymptotically defined by $\omega(t) = t \exp[\log^\alpha t]$ for $0 < \alpha < 1$. 
One finds
	\[ \varpi(t) \asymp \omega(t) \log^{1 - \alpha}_+ t . \]
Consequently, the asymptotic behavior of $\varpi$ can be strictly between $\omega$ and the bound given in \eqref{eq:UpperLowerBoundsVarpi2}.
\end{example}

\section{Optimal results for bounds of Paley--Wiener type}
\label{sec:EntireFunc}

We establish optimal existence results for entire functions with prescribed decay along the real direction and controlled growth along the imaginary direction.
Specifically, we investigate when, for two weight functions $\omega$ and $\eta$, an entire function $F$ can satisfy
	\begin{equation}
		\label{eq:OptimalBoundEntire} 
		|F(z)| \lesssim \exp\left[ \eta\left(\frac{|y|}{b}\right) - \omega(a |x|) \right] , \qquad z = x + iy \in \C , 
	\end{equation}
for certain $a, b > 0$.
To express our results efficiently, we introduce the following notions.
We write $\mathcal{O}(\C)$ for the space of all entire functions on $\C$ endowed with its natural Fr\'echet space topology.

\begin{definition}
Let $\omega$ and $\eta$ be weight functions.
We call $(\omega, \eta)$ a \emph{Paley-Wiener pair} if there exists a $\Psi_{(\omega, \eta)} \in (0, \infty]$ such that for any $a, b > 0$ the following hold:
	\begin{itemize}
		\item[(I)] If $a b > \Psi_{(\omega, \eta)}$, should $F \in \mathcal{O}(\C)$ satisfy \eqref{eq:OptimalBoundEntire}, then $F = 0$.	
		
		\item[(II)] If $a b < \Psi_{(\omega, \eta)}$, the space of $F \in \mathcal{O}(\C)$ satisfying \eqref{eq:OptimalBoundEntire} is dense in $\mathcal{O}(\C)$.
	\end{itemize}
We call $\Psi_{(\omega, \eta)}$ the \emph{Paley-Wiener threshold} of $(\omega, \eta)$.
The Paley-Wiener pair $(\omega, \eta)$ is called \emph{optimal} if $\Psi_{(\omega, \eta)} < \infty$.
\end{definition}

One of the main observations of this work is that $(\omega, \varpi)$ is an optimal Paley-Wiener pair.
We also give an explicit lower bound for $\Psi_{(\omega, \varpi)}$ when $\omega$ is convex; see Theorem \ref{t:PaleyWienerThresholdLowerBound}.
The main purpose of this section is to show the following result.

\begin{theorem}
\label{t:CharNonTrivEntireWeightFunc}
Let $\omega$ be a moderate quasianalytic weight function.
Then $(\omega, \varpi)$ is an optimal Paley-Wiener pair.
The Paley-Wiener threshold $\Psi_{(\omega, \varpi)}$ is at most $\pi / 2$.
\end{theorem}

Before we continue, we first remark on the motivation of our naming convention.

\begin{remark}
Suppose $\omega$ is an upper moderate weight function satisfying $\log t = o(\omega(t))$, but it is not quasianalytic, that is,
	\[ \int_1^\infty \frac{\omega(t)}{t^2} \, \mathrm{d}t < \infty . \]
By \cite[Corollary 2.6]{B-M-T-UltradiffFourier} (see also \cite[Theorem 1.3.7]{B-LinPartDiffOpGenDist}), for any $R > 0$, there exists some nontrivial $\varphi \in C^\infty(\R)$ with $\supp \varphi \subseteq [-R, R]$ and
	\[ \int_{-\infty}^\infty |\widehat{\varphi}(\xi)| e^{\omega(b |\xi|)} \, \mathrm{d}\xi < \infty , \qquad \forall b > 0 . \]
The Paley-Wiener theorem as in \cite[Proposition 3.4]{B-M-T-UltradiffFourier} (see also \cite[Theorem 1.4.1]{B-LinPartDiffOpGenDist}) shows that the latter is equivalent to $\widehat{\varphi}$ being extendable to an entire function such that, for each $n \in \N$ and $\varepsilon > 0$,
	\[ |\widehat{\varphi}(z)| \lesssim \exp\left[ (R + \varepsilon) |y| - \omega(n |x|) \right] , \qquad z = x + iy \in \C . \]	
As $\varpi(t) \asymp t$, we see that $(\omega, \varpi)$ is a Paley-Wiener pair.
However, $\Psi_{(\omega, \varpi)} = \infty$, so $(\omega, \varpi)$ is not optimal.

On the other hand, suppose $F \in \mathcal{O}(\C)$ satisfies
	\[ |F(z)| \lesssim \exp\left[ \eta\left(\frac{|y|}{b}\right) - \omega(a |x|) \right] , \qquad z = x + iy \in \C , \]
for some nondecreasing $\eta : [0, \infty) \to [0, \infty)$ such that $\eta(t) = o(t)$.
The classical Paley-Wiener theorem implies that $F$ is the Fourier-Laplace transform of a smooth function with support in $\{0\}$.
This is only possible when $F = 0$.
Thus $(\omega, \eta)$ is not a Paley-Wiener pair.
Therefore, the choice of $\varpi$ in the nonquasianalytic case is essentially optimal.
For this, we may see Theorem \ref{t:CharNonTrivEntireWeightFunc} as an optimal extension of the well-known Paley--Wiener-type results.
\end{remark}

The finiteness of the Paley-Wiener threshold $\Psi_{(\omega, \varpi)}$ follows from Carleman's theorem (see \cite[Theorem 1.2.2, p.~2]{B-EntFunc}). 
We prove this in Section \ref{subsec:PaleyWienerThreshold}, along with several other general observations on the Paley-Wiener threshold.
 In Section \ref{subsec:InterpolEst}, we establish an interpolation estimate, Lemma \ref{l:InterpolaEst}, using an estimate of Hirschman in \cite{H-ConvexTh}. 
 This estimate is then used in Section \ref{subsec:ConstructEntFunc} to develop an explicit construction of nontrivial entire functions satisfying \eqref{eq:OptimalBoundEntire} with $\eta = \varpi$. 
 The explicit nature of our construction allows us to derive lower bounds for $\Psi_{(\omega, \varpi)}$. 
 Finally, we prove our main results in Section \ref{subsec:ProofsMainThm}.

\subsection{Upper bounds for the Paley-Wiener threshold}
\label{subsec:PaleyWienerThreshold}

We start by establishing the following expression for $\Psi_{(\omega, \eta)}$.
It shows that it is equivalent to the threshold of existence of entire functions satisfying \eqref{eq:OptimalBoundEntire}.

\begin{lemma}
\label{l:PsiExplicit}
Let $\omega, \eta$ be lower moderate weight functions.
Then $(\omega, \eta)$ is a Paley-Wiener pair if and only if there exists $F \in \mathcal{O}(\C) \setminus \{0\}$ satisfying \eqref{eq:OptimalBoundEntire} for some $a_0, b_0 > 0$.
In this case 
	\begin{equation}
		\label{eq:PsiExplicit}
		\Psi_{(\omega, \eta)} = \sup \{ ab \mid \text{there exists $F \in \mathcal{O}(\C) \setminus \{0\}$ satisfying \eqref{eq:OptimalBoundEntire} for $a$ and $b$} \} .
	\end{equation} 
\end{lemma}

\begin{proof}
If $F \in \mathcal{O}(\C)$ satisfies \eqref{eq:OptimalBoundEntire} for certain $a_0, b_0 > 0$, and if $a_1, b_1 > 0$ are such that $a_0 b_0 = a_1 b_1$, then $G(z) = F\left(\frac{a_1}{a_0} z\right)$ satisfies \eqref{eq:OptimalBoundEntire} for $a_1$ and $b_1$.
Therefore, the supremum in \eqref{eq:PsiExplicit} is well-defined.

By definition, when $ab > \Psi_{(\omega, \eta)}$ and if $F \in \mathcal{O}(\C)$ satisfies \eqref{eq:OptimalBoundEntire} for $a$ and $b$, then necessarily $F = 0$.

Let now $a$ and $b$ be such that $ab$ is strictly smaller than the right-hand side of \eqref{eq:PsiExplicit}.
Take any $a_0 > a$ and $b_0 > b$ such that $a_0 b_0$ is also strictly smaller than the right-hand side of \eqref{eq:PsiExplicit}.
Then there is some $F_0 \in \mathcal{O}(\C) \setminus \{0\}$ satisfying \eqref{eq:OptimalBoundEntire} for $a_0$ and $b_0$.
The set
	\begin{equation}
		\label{eq:DenseSetEntireFunc} 
		\spn \, \{ z^m F^{(n)}_0 \mid m, n \in \N \} 
	\end{equation}
is dense in $\mathcal{O}(\C)$.
Indeed, suppose $z_0 \in \C$ is such that $F_0(z_0) \neq 0$.
Then, for any polynomial $P$ on $\C$ it follows by Taylor approximation that
	\[ \frac{P(z)}{F_0(z_0)} \sum_{k \leq K} \frac{(z_0 - z)^k}{k!} F_0^{(k)}(z) \to P(z) \quad \text{ in $\mathcal{O}(\C)$ as $K \to \infty$} . \] 
Since each function on the left-hand side belongs to \eqref{eq:DenseSetEntireFunc}, and since the set of polynomials is dense in $\mathcal{O}(\C)$, this gives density.
It therefore remains to show that each element in \eqref{eq:DenseSetEntireFunc} satisfies \eqref{eq:OptimalBoundEntire} for $a$ and $b$.
Fix some $m, n \in \N$ and $\varepsilon > 0$.
Let $H, C > 0$ be such that $(1 - \varepsilon)^{-1} \omega(t) \leq \omega(H t) + C$ and $(1 + \varepsilon) \eta(t) \leq \eta(H t) + C$.
Take any $\delta > 0$ such that $b < b_1 \leq b_0 / (1 + \delta)$ and $a < a_1 \leq (1 - \delta) a_0$ for certain $a_1$ and $b_1$.
By the Cauchy inequalities and the fact that $\log t = o(\omega(t))$, we find, for $|z| = |x + iy|$ large enough,
	\begin{align*}
		|z^m F^{(n)}_0(z)| 
		&\leq |z|^m \cdot \frac{n!}{\delta^n} \sup_{|\zeta| \leq \delta} |F_0(z - \zeta)| 
		\lesssim \sup_{|\zeta| \leq \delta} \sum_{m' \leq m} {m \choose m'} |z - \zeta|^{m'} |F_0(z - \zeta)| \\
		&\lesssim \sup_{|\zeta| = |u + iv| \leq \delta} \exp\left[ \eta\left(\frac{|y - v|}{b_0}\right) + m \log_+ |z - \zeta| - \omega(a_0 |x - u|) \right] \\
		&\lesssim \sup_{|\zeta| = |u + iv| \leq \delta} \exp\left[ \eta\left(\frac{H}{b_0} |y - v|\right) - \omega\left(\frac{a_0}{H}  |x - u|\right) \right] \\
		&\lesssim \exp\left[ \eta\left(\frac{H}{b_1} |y|\right) - \omega\left(\frac{a_1}{H} |x|\right) \right] .
	\end{align*}
Take any $c_\varepsilon > \max\{\underline{\mathfrak{m}}_{(1 - \varepsilon)^{-1}}(\omega), \underline{\mathfrak{m}}_{(1 + \varepsilon)}(\eta)\}$.
Then $z^m F^{(n)}_0$ satisfies \eqref{eq:OptimalBoundEntire} for any $a_1 / c_\varepsilon$ and $b_1 / c_{\varepsilon}$.
We may let $c_\varepsilon \to 1$ when $\varepsilon \to 0^+$ by Lemma \ref{l:Continuitymathfrakm}.
So, we see that $z^m F^{(n)}_0$ must also satisfy \eqref{eq:OptimalBoundEntire} for $a$ and $b$.
This concludes the proof.
\end{proof}

Next, we establish a necessary bound for $\omega$ by $\eta$ when a nontrivial entire function exists satisfying \eqref{eq:OptimalBoundEntire}.

\begin{proposition}
\label{p:NonTrivEntireLowerBound}
Let $\omega, \eta$ be weight functions such that $\omega$ is quasianalytic and $\eta$ is convex.
If there exists an $F \in \mathcal{O}(\C) \setminus \{0\}$ satisfying \eqref{eq:OptimalBoundEntire} for some $a, b > 0$,
then for any $0 < \theta < ab$, when $t$ is large enough,
	\[ \varpi(t) \leq \eta\left(\frac{\pi}{2 \theta} t\right) . \]
\end{proposition}

\begin{proof}
Let $n \in \N$ be the smallest number for which $F^{(n)}(0) \neq 0$.
Put
	\[ G(z) = \frac{n!}{F^{(n)}(0)} \frac{F(z)}{z^n} . \]
Then $G$ is an entire function with $G(0) = 1$ and, for some $C_0 > 0$,
	\[ |G(z)| \leq C_0 \exp\left[ \eta\left(\frac{|y|}{b}\right) - \omega(a |x|) - n \log_+ |z| \right] , \qquad z = x + iy \in \C . \]
Let $C_1 > 0$ be such that
	\[ |G(x) G(-x)| \leq 1 + C_1 x^2 , \qquad |x| \leq 1 . \]
Using Carleman's theorem \cite[Theorem 1.2.2, p.~2]{B-EntFunc} and the convexity of $\eta$, it follows that, for any $t \geq 1$,
	\begin{align*} 
		&-\frac{\pi}{2} \, \im G'(0) 
		\leq \frac{1}{t} \int_0^\pi \log |G(t e^{i \theta})| \sin \theta \, \mathrm{d}\theta + \frac12 \int_0^t \left(\frac{1}{x^2} - \frac{1}{t^2}\right) \log |G(x) G(-x)| \, \mathrm{d}x \\
		&\quad \leq \frac{1}{t} \int_0^\pi \left[ \log C_0 + \eta\left(\frac{t}{b} \sin \theta\right) - \omega(a t |\cos \theta|) - n \log t \right] \sin \theta \, \mathrm{d}\theta \\
		&\qquad + \int_0^1 \left(\frac{1}{x^2} - \frac{1}{t^2}\right) \log[1 + C_1 x^2] \, \mathrm{d}x + \int_1^t \left(\frac{1}{x^2} - \frac{1}{t^2}\right) [\log C_0 - \omega(a x) - n \log x] \, \mathrm{d}x \\
		&\quad \leq C_2+ \frac{1}{t} \left[ \eta\left(\frac{\pi t}{2 b}\right) - \varpi(a t) - \frac{1}{t} \int_0^t \omega(a u) \, \mathrm{d}u \right] ,
	\end{align*}
for some $C_2 > 0$.
Thus, for some $C_3 > 0$,
	\[ \varpi(a t)  - C_3 t \leq \eta\left(\frac{\pi t}{2 b}\right) , \qquad t \geq 1 . \]
As $t = o(\varpi(t))$ by Lemma \ref{l:PropertiesVarpi}$(ii)$, we find by the convexity of $\varpi$ that for any $0 < \varepsilon < 1$, when $t$ is large enough,
	\[ \varpi((1 - \varepsilon) a t) \leq \eta\left(\frac{\pi t}{2 b}\right) . \qedhere \]
\end{proof}

\begin{corollary}
\label{c:PaleyWienerThresholdBelowPi/2}
Let $\omega$ be a lower moderate quasianalytic weight function.
If $(\omega, \varpi)$ is a Paley-Wiener pair, then 
	\[ \Psi_{(\omega, \varpi)} \leq \frac{\pi}{2} . \] 
\end{corollary}

\begin{proof}
Suppose $(\omega, \varpi)$ is a Paley-Wiener pair and $\Psi_{(\omega, \varpi)} > \pi/2$.
By Lemma \ref{l:PsiExplicit}, there exists some $F \in \mathcal{O}(\C) \setminus \{0\}$ satisfying \eqref{eq:OptimalBoundEntire} for certain $a$ and $b$ such that $ab > \pi/2$, say $ab = (1 + \varepsilon) \pi / 2$ with $0 < \varepsilon < 1$.
Then, by Proposition \ref{p:NonTrivEntireLowerBound}, for $t$ large enough, 
	\[ \varpi(t) \leq \varpi\left(\frac{\pi}{2 (1 - \varepsilon/2) ab} t\right) \leq \varpi\left(\frac{t}{1 + (\varepsilon - \varepsilon^2) / 2}\right) . \]
As $\varpi$ is nonconstant and nondecreasing, this gives a contradiction.
\end{proof}

We now give a general upper bound for $\Psi_{(\omega, \eta)}$ when $\eta$ is polynomially bounded.

\begin{proposition}
\label{prop:UpperBoundPaleyWienerThreshold}
Let $\omega$ and $\eta$ be weight functions such that $\eta(t) = O(t^{\alpha + \varepsilon})$ for some $\alpha > 1$ and every $\varepsilon > 0$.
Then, if $(\omega, \eta)$ is a Paley-Wiener pair,
	\[ \Psi_{(\omega, \eta)} \leq \tan\left(\frac{\alpha-1}{\alpha} \frac{\pi}{2}\right) \mathfrak{c}(\eta, \omega) .  \]
\end{proposition}

\begin{proof}
We may assume $\mathfrak{c}(\eta, \omega) < \infty$, otherwise the statement is trivial.
Fix any $\alpha' > \alpha$ and $L > \mathfrak{c}(\eta, \omega)$.
Suppose $F \in \mathcal{O}(\C) \setminus \{0\}$ satisfies \eqref{eq:OptimalBoundEntire} for certain $a, b > 0$ satisfying
	\[ ab \geq \tan\left(\frac{\alpha' - 1}{\alpha'} \frac{\pi}{2}\right) L . \]
Put $\theta_0 = \pi(\alpha' - 1) / (2\alpha') \leq \arctan(ab / L)$. 
Then, for any $r \in \R$ and $|\theta| \leq \theta_0$,
	\begin{align*}
		|F(r e^{i \theta})|
		&\lesssim \exp\left[\eta\left(\frac{|r \sin \theta|}{b}\right) - \omega(a  |r \cos \theta|)\right] \\
		&\leq \exp\left[\eta\left(\frac{a}{L}  |r \cos \theta|\right) - \omega(a |r \cos \theta|)\right]
		\lesssim 1 . 
	\end{align*}
Consequently, $F$ is bounded in $\{ z \in \C \mid |\arg z| \leq \theta_0 \text{ or } |\arg z - \pi| \leq \theta_0 \}$.
The complement of this set consists of two sectors of angle $\pi - 2\theta_0 = \pi / \alpha'$.
As $|F(z)| \lesssim e^{C |z|^{\alpha''}}$ for any $\alpha < \alpha'' < \alpha'$ and some $C > 0$, the Phragm\'en-Lindel\"of principle as in \cite[Theorem 1.4.2, p.~4]{B-EntFunc} shows that $F$ is bounded on $\C$.
By Liouville's theorem, we have that $F$ is constant.
As $F$ decays to zero on the real line, we finally conclude that $F = 0$.
Consequently, it follows that
	\[ \Psi_{(\omega, \eta)} \leq \tan\left(\frac{\alpha' - 1}{\alpha'} \frac{\pi}{2}\right) L . \]
Letting $\alpha' \to \alpha$ and $L \to \mathfrak{c}(\eta, \omega)$ shows the result.
\end{proof}

\begin{corollary}
Let $\omega$ be a weight function such that $\omega(t) = O(t^{\alpha + \varepsilon})$ for some $\alpha > 1$ and every $\varepsilon > 0$.
Then, if $(\omega, \omega)$ is a Paley-Wiener pair,
	\[ \Psi_{(\omega, \omega)} \leq \tan\left(\frac{\alpha - 1}{\alpha} \frac{\pi}{2}\right) . \]
\end{corollary}

\subsection{Interpolation estimates}
\label{subsec:InterpolEst}

We show that the bound \eqref{eq:OptimalBoundEntire} can often be deduced directly from the corresponding decay estimate on the real line and a suitable growth bound along the imaginary direction. 
More precisely, we have the following general result.
For a weight function $\omega$, we write
	\begin{equation}
		\label{eq:IsConvex}
		\kappa_\omega = \begin{cases} 1 , & \omega \text{ is convex} , \\ 0 , & \text{otherwise} . \end{cases}
	\end{equation}

\begin{lemma}
\label{l:InterpolaEst}
Let $\omega, \eta$ be weight functions.
Suppose $F \in \mathcal{O}(\C)$ is such that, for certain $a, b > 0$,
	\begin{equation}
		\label{eq:HirschmanRealLine} 
		|F(x)| \lesssim e^{- \omega(a |x|)} , \qquad x \in \R , 
	\end{equation}
and
	\begin{equation}
		\label{eq:HirschmanImaginary} 
		|F(z)| \lesssim e^{\eta\left(\frac{|y|}{b}\right)} , \qquad z = x + iy \in \C . 
	\end{equation}
Then, for any $\theta \in (0, 1)$,
	\begin{equation}
		\label{eq:HirschmanFull} 
		|F(z)| \lesssim e^{\theta \eta\left(\frac{1}{\theta b} |y|\right) - \frac{(1 + \kappa_\omega)(1 - \theta)}{2} \omega(a |x|)} , \qquad z = x + iy \in \C . 
	\end{equation}
\end{lemma}

\begin{proof}
Fix some $\theta \in (0, 1)$ and $z = x + iy \in \C$ with $y \neq 0$.
We define
	\begin{equation}
		\label{eq:GDef} 
		G(\zeta) = F\left(x +  \frac{i y \zeta}{\theta}\right) , \qquad 0 \leq \re \zeta \leq 1 . 
	\end{equation}
Then $G$ is analytic and bounded on the strip $\{ \zeta \in \C \mid 0 \leq \re \zeta \leq 1 \}$ by \eqref{eq:HirschmanImaginary}.
Note that $G(\theta) = F(z)$.
We may find some $C_0 > 0$, independent of $z$, such that, for $t \in \R$,
	\begin{equation}
		\label{eq:G(it)Bound} 
		\log |G(i t)| 
		= \log \left|F\left(x - \frac{y t}{\theta}\right)\right| \leq \log C_0 - \omega\left(a \left|x - \frac{yt}{\theta}\right|\right)
	\end{equation}
and
	\begin{equation}
		\label{eq:G(1+it)Bound}
		\log |G(1 + it)|
		= \log \left|F\left(x - \frac{y t}{\theta} + \frac{i y}{\theta}\right)\right| 
		\leq \log C_0 + \eta\left(\frac{|y|}{\theta b}\right) .
	\end{equation}
Put, for $\vartheta \in (0, 1)$,
	\begin{equation}
		\label{eq:Ptheta} 
		P_\vartheta(t) = \frac{\frac12 \tan \frac{\vartheta \pi}{2}}{\left[ \tan^2 \frac{\vartheta \pi}{2} + \tanh^2 \frac{t \pi}{2} \right] \cosh^2 \frac{t \pi}{2}} , \qquad t \in \R . 
	\end{equation}
By the interpolation estimate of Hirschman in \cite[Lemma 1]{H-ConvexTh}, we have
	\begin{equation}
		\label{eq:GInterpolEstimate} 
		\log |F(z)| = \log |G(\theta)| \leq \int_{-\infty}^\infty P_\theta(t) \log |G(it)| \, \mathrm{d}t + \int_{-\infty}^\infty P_{1-\theta}(t) \log |G(1 + it)| \, \mathrm{d}t .
	\end{equation}
We estimate the two integrals separately.
Note that $\int_{-\infty}^\infty P_\vartheta(t) \, \mathrm{d}t = 1 - \vartheta$ for $\vartheta \in (0, 1)$.
For the second integral, we have by \eqref{eq:G(1+it)Bound},
	\[ \int_{-\infty}^\infty P_{1-\theta}(t) \log |G(1 + it)| \, \mathrm{d}t \leq \theta \eta\left(\frac{|y|}{\theta b}\right) + \log C_0^\theta . \]
	
We now estimate the first integral in \eqref{eq:GInterpolEstimate}.
For any $t \in \R$ we have
	\[ \omega\left(a\left|x - \frac{yt}{\theta}\right|\right) + \omega\left(a\left|x + \frac{y t}{\theta}\right|\right) \geq (1 + \kappa_\omega) \omega(a |x|) . \]
Therefore, as $P_\theta$ is even, by \eqref{eq:G(it)Bound},
	\begin{align*}
		\int_{-\infty}^\infty P_\theta(t) \log |G(it)| \, \mathrm{d}t 
		&\leq \int_{-\infty}^\infty P_\theta(t) \left[ \log C_0 - \omega\left(a \left|x - \frac{y t}{\theta}\right|\right)\right] \, \mathrm{d}t \\
		&= \log C_0^{1-\theta} - \int_0^\infty P_\theta(t) \left[ \omega\left(a\left|x - \frac{yt}{\theta}\right|\right) + \omega\left(a\left|x + \frac{y t}{\theta}\right|\right) \right] \, \mathrm{d}t \\
		&\leq \log C_0^{1-\theta} - \frac{(1 + \kappa_\omega)(1 - \theta)}{2} \omega(a |x|) . \qedhere
	\end{align*}
\end{proof}

\subsection{Construction of entire functions satisfying optimal Paley--Wiener-type bounds}
\label{subsec:ConstructEntFunc}

We now lay out a way to construct nontrivial entire functions satisfying \eqref{eq:OptimalBoundEntire} for $\eta = \varpi$.
To do this, we need the following characterization of upper moderateness.
In particular, we get an alternative expression for the upper Matuszewska index $\alpha(\omega)$ as in Definition \ref{def:Matuszewska}.

\begin{lemma}
\label{l:ModerateIffBeta*}
Let $\omega$ be a weight function.
Then $\omega$ is upper moderate if and only if for some $\alpha > 0$ we have
	\begin{equation} 
		\label{eq:b*cond}
		\int_t^\infty \frac{\omega(s)}{s^{\alpha + 1}} \, \mathrm{d}s = O\left(\frac{\omega(t)}{t^\alpha}\right) . 
	\end{equation}
In this case, $\alpha(\omega) = \inf \{ \alpha > 0 \mid $\eqref{eq:b*cond} holds for $\alpha\}$ and $\omega(t) = o(t^\alpha)$.
\end{lemma}

\begin{proof}
See, for example, \cite[Theorem 2.11]{JG-S-S-ORegVar}.
\end{proof}

\begin{definition}
Let $\omega$ be a weight function.
For $\alpha > 0$, we write
	\[ P_\alpha(\omega) = \inf\left\{ P \in (0, \infty] \mid \exists C > 0 ~ \forall t \geq 1 : \alpha t^\alpha \int_t^\infty \frac{\omega(s)}{s^{\alpha + 1}} \, \mathrm{d}s \leq P \omega(t) + C \right\} . \]
Note that always $P_\alpha(\omega) \geq 1$.
\end{definition}

Given a suitable entire function, we may then construct a nontrivial $F$ satisfying \eqref{eq:OptimalBoundEntire} for $\eta = \varpi$ in the following manner.
Below, $[t]_+ = \max\{0, t\}$ for $t \in \R$.

\begin{proposition}
\label{p:ConstructionFunction}
Let $\omega$ be a moderate quasianalytic weight function.
Let $G$ be an entire function of finite exponential type $\upsilon \geq 1$ that is real-valued on $\R$ and satisfies:
	\begin{itemize}
		\item[$(i)$] $G(0) = 1$;
		\item[$(ii)$] $|G(x)| \leq \min\{1, |x|^{-2}\}$ for $|x| > 0$;
		\item[$(iii)$] There is some $\alpha > \max\{1, \alpha(\omega)\}$ and $C_{G} > 0$ for which $\log |G(x + iy)| \leq C_{G} |y|^{\alpha}$ when $|y| \leq 1$.
	\end{itemize}
Then
	\begin{equation*}
		\label{eq:SeqToElement} 
		F(z) = \prod_{p = 1}^\infty G\left(\frac{z}{\mu_p}\right) , \qquad \mu_p = \inf \{ t \geq 0 \mid \mu(t) \geq p \} ,
	\end{equation*}
defines an entire function with $F(0) = 1$ satisfying, for some $a, b > 0$,
	\begin{equation}
		\label{eq:OptimalBoundEntireWeightSeq} 
		|F(z)| \lesssim \exp\left[ \varpi\left(\frac{|y|}{b}\right) - \omega(a |x|) \right] , \qquad z = x + iy \in \C .
	\end{equation}
It holds that, with $\mathfrak{c}_+ = \max\{1, \mathfrak{c}(\omega, \varpi)\}$,
	\begin{equation}
		\label{eq:LowerBoundThresholdExplicitConstructGeneral}
		\sup \{ ab \mid F \text{ satisfies \eqref{eq:OptimalBoundEntireWeightSeq} for $a, b > 0$} \} \geq \frac{1}{e \mathfrak{c}_+} \sup_{0<\theta<1} \frac{\theta}{\underline{\mathfrak{m}}_{\frac{1}{1-\theta}}(\omega) \underline{\mathfrak{m}}_{\theta K_{G}}(\omega)} 
	\end{equation}
with
	\begin{equation}
		\label{eq:K}
		K_G = \frac{\upsilon}{e \mathfrak{c}_+} + \frac{\mathfrak{c}(\omega, \varpi)}{e \mathfrak{c}_+} [\upsilon + (\alpha P_\alpha(\omega) -\alpha)C_G - 2]_+ + \frac{\mathfrak{c}(\omega, \varpi)}{\mathfrak{c}_+} [\upsilon - C_G]_+ .
	\end{equation}
When $\omega$ is convex and $C_G < \upsilon$, we have
	\begin{equation}
		\label{eq:LowerBoundThresholdExplicitConstruct1}
		\sup \{ ab \mid F \text{ satisfies \eqref{eq:OptimalBoundEntireWeightSeq} for $a, b > 0$} \} \geq \frac{1}{e} \sup_{0<\theta<1} \frac{\theta}{\underline{\mathfrak{m}}_{\frac{1}{2(1 - \theta)}}(\omega) \underline{\mathfrak{m}}_{\theta K_{1, G}}(\omega)} 
	\end{equation}
with
	\begin{equation}
		\label{eq:K1}
		K_{1, G} =  \upsilon \left[1 - \frac{e-1}{\overline{m}_e(\omega)}\right]_+ + \frac{\mathfrak{c}(\omega, \varpi)}{e}[(\alpha P_\alpha(\omega) - \alpha - 1) C_G - (e - 2) \upsilon - 2]_+ .
	\end{equation}
When $\omega$ is convex and $C_G \geq \upsilon$, we have
	\begin{equation}
		\label{eq:LowerBoundThresholdExplicitConstruct2}
		\sup \{ ab \mid F \text{ satisfies \eqref{eq:OptimalBoundEntireWeightSeq} for $a, b > 0$} \} \geq \frac{1}{\mathfrak{c}_+} \sup_{0<\theta<1} \frac{\theta}{\underline{\mathfrak{m}}_{\frac{1}{2(1 - \theta)}}(\omega) \underline{\mathfrak{m}}_{\theta K_{2, G}}(\omega)} 
	\end{equation}
with
	\begin{equation}
		\label{eq:K2}
		K_{2, G} =\frac{\upsilon}{\mathfrak{c}_+} + \frac{\mathfrak{c}(\omega, \varpi)}{\mathfrak{c}_+}(2 \upsilon + (\alpha P_\alpha(\omega) - \alpha - 1) C_G - 2) .
	\end{equation}
\end{proposition}

\begin{proof}
We claim that the products
	\begin{equation}
		\label{eq:InfiniteProduct} 
		\prod_{p = 1}^q G\left(\frac{z}{\mu_p}\right) , \qquad q \geq 1 ,
	\end{equation}
converge uniformly on each compact to an entire function $F$ satisfying \eqref{eq:OptimalBoundEntire} for $\eta = \varpi$ and certain $a, b > 0$.
This would show the first part of the result.

As $G$ is of finite exponential type $\upsilon \geq 1$, and since $|x^2 G(x)| \leq 1$ for all $x \in \R$ by $(ii)$, the Phragm\'en-Lindel\"of principle as in \cite[Theorem 1.4.3, p.~4]{B-EntFunc} yields
	\begin{equation}
		\label{eq:PsiGlobalBound} 
		|G(z)| \leq \exp\left[ \upsilon |y| - 2 \log |z| \right] , \qquad z = x + iy \in \C \setminus \{0\} . 
	\end{equation}
	
Set
	\begin{equation}
		\label{eq:M(t)} 
		M(t) = \lfloor \mu (t) \rfloor = \sum_{\mu_p \leq t} 1 , \qquad t \geq 0 . 
	\end{equation}
Let $\alpha > \alpha(\omega)$ be such that $(iii)$ holds.
By \eqref{eq:muBoundomega} and Lemma \ref{l:ModerateIffBeta*}, we have that $\omega(t) = o(t^{\alpha})$, $M(t) = o(t^\alpha)$, and $\mu(t) = o(t^{\alpha})$.
Also, there are $C, P > 0$ for which $\alpha t^\alpha \int_t^\infty \omega(s) / s^{\alpha + 1} \, \mathrm{d}s \leq P \omega(t) + C$.
Integration by parts and \eqref{eq:omegaAsIntegral} show
	\begin{align}
		t^\alpha \int_{t^+}^\infty \frac{\mathrm{d} M(u)}{u^{\alpha}}
		&= \alpha t^\alpha \int_{t}^\infty \frac{M(u)}{u} \frac{\mathrm{d}u}{u^{\alpha}} - M(t) 
		\leq \alpha^2 t^\alpha \int_{t}^\infty \frac{\omega(u)}{u^{\alpha+1}} \, \mathrm{d}u - \alpha \omega(t) - \mu(t) + 1
		\nonumber \\
		&\leq \alpha (P - 1) \omega(t) - \mu(t) + \alpha C + 1 
		\label{eq:BoundSum}
	\end{align}

We verify that the products \eqref{eq:InfiniteProduct} converge in $\mathcal{O}(\C)$.
First, by $(iii)$ and \eqref{eq:BoundSum} we have, for any $q \in \Z_+$ and $z = x + iy \in \C$,
	\begin{align*}
	 	\sum_{p = 1}^q \log \left| G\left(\frac{z}{\mu_p}\right) \right| 
		&\leq \sum_{\mu_p \leq |z|} \log_+ \left| G\left(\frac{z}{\mu_p}\right) \right| + C_G |y|^\alpha \int_{|z|^+}^\infty \frac{\mathrm{d}M(u)}{u^\alpha} \\
		&\leq \sum_{\mu_p \leq |z|} \log_+ \left| G\left(\frac{z}{\mu_p}\right) \right| + \alpha C_G (P - 1) \omega(|z|) + (\alpha C + 1) C_G . 
	\end{align*}
Therefore, the sequence of partial products \eqref{eq:InfiniteProduct} is uniformly bounded on each compact set.
Now, since $G$ is real-valued on $\R$ and by $(i)$ and $(ii)$, for any $x \in \R$ we have that eventually $0 < G(x / \mu_p) \leq 1$.
Thus the infinite product converges on $\R$.
Vitali's theorem now shows that \eqref{eq:InfiniteProduct} converges to an entire function, which is precisely $F$.

We now verify that the bound \eqref{eq:OptimalBoundEntireWeightSeq} holds for $\omega$ and some $a, b > 0$.
Fix $0 < \varepsilon < 1$.
Note that $2 \log t \leq \varepsilon \omega(t)$ when $t$ is large enough.
First, we consider the real line.
By $(ii)$ and \eqref{eq:omegaAsIntegral}, for $x \in \R$ large enough,
	\begin{align*}
		\log |F(x)|
		&\leq \int_{\mu_1}^{|x|} \log \left| G\left(\frac{x}{u}\right) \right| \, \mathrm{d} M(u) + \int_{|x|^+}^\infty \log \left| G\left(\frac{x}{u}\right) \right| \, \mathrm{d} M(u) \\
		&\leq - 2 \int_{\mu_1}^{|x|} \log \left(\frac{|x|}{u}\right) \, \mathrm{d} M(u)
		= - 2 \int_{\mu_1}^{|x|} \frac{M(u)}{u} \, \mathrm{d}u \\
		&\leq 2 \int_{\mu_1}^{|x|} \frac{\mathrm{d}u}{u} - 2 \int_{\mu_1}^{|x|} \frac{\mu(u)}{u} \, \mathrm{d}u 
		\leq 2 \log |x| - 2 \log \mu_1 - 2 \omega(|x|) + 2 \omega(\mu_1) \\
		&\leq -(2 - \varepsilon) \omega(|x|) + 2 \omega(\mu_1) - 2 \log \mu_1 .
	\end{align*}
Next, we observe that $F$ is bounded on any strip.
Indeed, let $\Delta \geq \max\{1, \mu_1\}$ and suppose $z = x + iy$ is such that $|z| \geq \Delta$ and $|y| \leq \Delta$.
Then, by $(iii)$, \eqref{eq:PsiGlobalBound}, and \eqref{eq:BoundSum},
	\begin{align*}
		\log |F(z)|
		&\leq \int_{\mu_1}^\infty \log \left|G\left(\frac{x}{u} + i \frac{y}{u}\right)\right| \, \mathrm{d}M(u) 
		\leq \upsilon |y| \int_{\mu_1}^\Delta \frac{\mathrm{d}M(u)}{u} +  C_{G} |y|^\alpha \int_{\Delta^+}^\infty \frac{\mathrm{d}M(u)}{u^\alpha} \\
		&\leq \frac{\upsilon \Delta}{\mu_1}M(\Delta) + C_G [\alpha (P - 1) \omega(\Delta) + \alpha C + 1] 
		< \infty.
	\end{align*}
Finally, we consider the case outside of a fixed strip.
Suppose $z = x + iy$ with $|y| \geq \max\{1, \mu_1\}$.
Then, by $(iii)$, \eqref{eq:omegaAsIntegral}, \eqref{eq:PsiGlobalBound}, and \eqref{eq:BoundSum}, 
	\begin{align*}
		&\log |F(z)|
		\leq \int_{\mu_1}^{|y|} \log \left|G\left(\frac{x}{u} + i \frac{y}{u}\right)\right| \, \mathrm{d}M(u) + \int_{|y|^+}^\infty \log \left|G\left(\frac{x}{u} + i \frac{y}{u}\right)\right| \, \mathrm{d}M(u) \\
		&\leq \upsilon |y| \int_{\mu_1}^{|y|} \frac{\mathrm{d}M(u)}{u} - 2 \int_{\mu_1}^{|y|} \log \frac{|y|}{u} \, \mathrm{d}M(u) + C_{G} |y|^\alpha \int_{|y|^+}^\infty \frac{\mathrm{d}M(u)}{u^\alpha} \\
		&\leq \upsilon |y| \left[ \frac{\mu(|y|)}{|y|} + \frac{\omega(|y|)}{|y|} + \int_{\mu_1}^{|y|} \frac{\omega(u)}{u^2} \, \mathrm{d}u \right] 
		- 2 \omega(|y|) + 2 \log |y| - 2 \log \mu_1 \\ & \hspace{6cm}  + [\alpha (P - 1) C_{G}] \omega(|y|) - C_G \mu(|y|) + (\alpha C + 1) C_G \\
		&\leq \upsilon \varpi(|y|) + (\upsilon + \alpha(P-1)C_G - 2 + \varepsilon) \omega(|y|) + (\upsilon - C_G) \mu(|y|) \\ & \hspace{6cm} - \upsilon \left(\int_{1}^{\mu_1} \frac{\omega(u)}{u^2} \, \mathrm{d}u \right) |y| - 2 \log \mu_1 + (\alpha C + 1) C_G
	\end{align*}
Since $t = o(\varpi(t))$ by Lemma \ref{l:PropertiesVarpi}$(ii)$, if $|y|$ is large enough,
	\begin{equation}
		\label{eq:EstimateFOutsideStrip}
		\log |F(z)| \leq (\upsilon +  \varepsilon) \varpi(|y|) + (\upsilon + \alpha(P-1)C_G - 2 + \varepsilon) \omega(|y|) + (\upsilon - C_G) \mu(|y|) + O(1) . 
	\end{equation}
	
For general $\omega$, using \eqref{eq:muBoundomega} and Lemma \ref{l:PropertiesVarpi}, and the fact that $F$ is bounded on strips, we may find $c, L > 0$ such that
	\begin{align*} 
		|F(z)| &\lesssim e^{L \varpi(c |y|)} , \qquad z = x + iy \in \C . \\
\intertext{Also, as we have seen, on the real line we have}
		 |F(x)| &\lesssim e^{-(2 - \varepsilon)\omega\left(|x|\right)} , \qquad x \in \R . 
	\end{align*}
Fix $\theta \in (0, 1)$.
Lemma \ref{l:InterpolaEst} with the weight functions $(2 - \varepsilon)\omega$ and $L \varpi(c \, \cdot)$ yields
	\[ |F(z)| \lesssim e^{\theta L \varpi\left(\frac{c |y|}{\theta}\right) - (1 + \kappa_\omega)(1 - \theta) \left(1 - \frac{\varepsilon}{2}\right) \omega(|x|)} , \qquad z = x + iy \in \C, \]
with $\kappa_\omega$ as in \eqref{eq:IsConvex}.
Lemmas \ref{l:Continuitymathfrakm} and \ref{l:IndicesForVarpi} then show that the supremum of $\{ ab \mid F \text{ satisfies \eqref{eq:OptimalBoundEntireWeightSeq} for $a, b > 0$} \}$ is at least
	\[ \left[ \frac{c}{\theta} \, \underline{\mathfrak{m}}_{\frac{1}{(1 + \kappa_\omega)(1 - \theta)}}(\omega) \underline{\mathfrak{m}}_{\theta L}(\omega) \right]^{-1} . \]
Depending on whether or not $\omega$ is convex, we now give upper bounds for $c$ and $L$.
Throughout our arguments, we fix $H_0, C_0 > 0$ such that $\omega(t) \leq \varpi(H_0 t) + C_0$.
By Lemma \ref{l:PropertiesVarpi}, we may assume $H_0 < e$.
We also set $H_{0, +} = \max\{1, H_0\}$.

Suppose first that $\omega$ is not necessarily convex.
By \eqref{eq:muBoundomega}, \eqref{eq:EstimateFOutsideStrip}, and the convexity of $\varpi$, we find, for $|y|$ large enough,
	\begin{align*}
		&\log |F(z)| \\
		&~ \leq (\upsilon +  \varepsilon) \varpi(|y|) + (\upsilon + \alpha(P-1)C_G - 2 + \varepsilon) \omega(|y|) + [\upsilon - C_G]_+ \omega(e |y|) \\
		&~ \leq \left[ \frac{\upsilon + \varepsilon}{e H_{0, +}} + \frac{H_0}{e H_{0, +}} [\upsilon + \alpha(P-1)C_G - 2 + \varepsilon]_+ + \frac{H_0}{H_{0,+}} [\upsilon - C_G]_+ \right] \varpi(e H_{0, +} |y|) + O(1) .
	\end{align*}
By Lemmas \ref{l:Continuitymathfrakm} and \ref{l:IndicesForVarpi}, we may conclude \eqref{eq:LowerBoundThresholdExplicitConstructGeneral}.

We now assume $\omega$ is convex.
Then $\kappa_\omega = 1$.
Suppose first that $C_G < \upsilon$.
We have
	\[ \omega(s) \geq \omega(t) + \frac{\mu(t)}{t} (s - t) , \qquad 0 < t \leq s . \]
Now, for any $\lambda > 1$ and $t \geq 1$,
	\[ \varpi(\lambda t) = \lambda \varpi(t) + \lambda t \int_{t}^{\lambda t} \frac{\omega(s)}{s^2} \, \mathrm{d}s \geq \lambda \varpi(t) + (\lambda - 1) \omega(t) + (\lambda \log \lambda + 1 - \lambda) \mu(t) . \]
In particular, for $\lambda = e$ we have
	\[ \mu(t) + (e - 1) \omega(t) + e \varpi(t) \leq \varpi(e t) , \qquad t \geq 1 . \]
Suppose $\varpi(et) \leq H_1 \varpi(t) + C_1$ for $H_1, C_1 > 0$.
Using \eqref{eq:OmegaBoundMuConvex} and the convexity of $\varpi$, the estimate \eqref{eq:EstimateFOutsideStrip} gives, for $|y|$ large enough,
	\begin{align*} 
		&\log |F(z)| \\
		&\leq (\upsilon + \varepsilon) \varpi(|y|) + (\upsilon + (\alpha(P-1) - 1)C_G - 2 + 2\varepsilon) \omega(|y|) + \upsilon \mu(|y|) \\
		&\leq \upsilon \varpi(e |y|) - (\upsilon(e - 1) - \varepsilon) \varpi(|y|) + ((\alpha(P-1) - 1)C_G - 2 - \upsilon(e - 2) + 2\varepsilon) \omega(|y|) \\
		&\leq \left(\upsilon - \frac{e-1}{H_1} \upsilon + \frac{\varepsilon}{H_1}\right) \varpi(e |y|) + ((\alpha(P-1) - 1)C_G - 2 - \upsilon(e - 2) + 2\varepsilon) \omega(|y|) + O(1)\\
		&\leq \left(\upsilon\left(1 - \frac{e-1}{H_1}\right) + \frac{\varepsilon}{H_1}  + \frac{H_0}{e} \left[ (\alpha(P-1) - 1)C_G - 2 - \upsilon\left(e - 2\right) + 2\varepsilon \right]_+ \right) \varpi(e |y|) + O(1) .
	\end{align*}
From here, using Lemmas \ref{l:Continuitymathfrakm} and \ref{l:IndicesForVarpi}, we may conclude \eqref{eq:LowerBoundThresholdExplicitConstruct1}.

Suppose now $C_G \geq \upsilon$.
Using \eqref{eq:OmegaBoundMuConvex}, \eqref{eq:EstimateFOutsideStrip}, and the convexity of $\varpi$, we get
	\begin{align*}
		\log |F(z)|
		&\leq (\upsilon + \varepsilon) \varpi(|y|) + ( (\alpha(P-1) - 1) C_G + 2\upsilon - 2 + \varepsilon) \omega(|y|) \\
		&\leq \left( \frac{\upsilon + \varepsilon}{H_{0, +}} + \frac{H_{0}}{H_{0,+}} ( (\alpha(P-1) - 1) C_G  + 2 \upsilon - 2 + \varepsilon) \right) \varpi(H_{0, +} |y|) + O(1) .
	\end{align*}
By Lemma \ref{l:Continuitymathfrakm}, this gives \eqref{eq:LowerBoundThresholdExplicitConstruct2}.
\end{proof}

The following observation is relevant to the constants \eqref{eq:K1} and \eqref{eq:K2}.

\begin{lemma}
\label{l:(aPa - a)nonincreasing}
Let $\omega$ be an upper moderate weight function.
Then
	\[ \alpha P_\alpha(\omega) - \alpha , \qquad \alpha > \alpha(\omega) \]
is nonincreasing.
\end{lemma}

\begin{proof}
Take any $\alpha(\omega) < \alpha < \beta$.
Let $t > 0$ be such that $\omega(t) > 0$.
By a change of variables, we have
	\[ \beta t^\beta \int_{t}^\infty \frac{\omega(s)}{s^{\beta + 1}} \, \mathrm{d}s = \int_1^\infty \frac{\omega(t u^{1/\beta})}{u^2} \, \mathrm{d}u . \]
As the function $t \mapsto \omega(e^t)$ is convex, we have
	\[ \omega(t u^{1/\beta}) \leq \frac{\alpha}{\beta} \omega(t u^{1/\alpha}) + \left( 1 - \frac{\alpha}{\beta}\right) \omega(t) , \qquad u \geq 1 . \]
Plugging this into the previous identity, we find
	\[ \beta t^\beta \int_{t}^\infty \frac{\omega(s)}{s^{\beta + 1}} \, \mathrm{d}s \leq \frac{\alpha}{\beta} \left( \alpha t^\alpha \int_{t}^\infty \frac{\omega(s)}{s^{\alpha + 1}} \, \mathrm{d}s \right) + \left(1 - \frac{\alpha}{\beta}\right) \omega(t) . \]
Rearranging this inequality and dividing by $\omega(t)$ then gives
	\[ \beta \frac{\beta t^\beta \int_{t}^\infty \frac{\omega(s)}{s^{\beta + 1}} \, \mathrm{d}s}{\omega(t)} - \beta \leq \alpha \frac{\alpha t^\alpha \int_{t}^\infty \frac{\omega(s)}{s^{\alpha + 1}} \, \mathrm{d}s}{\omega(t)} - \alpha . \]
Taking the limit superior proves the claim.
\end{proof}

Our next step is to verify that an entire function $G$ as in Proposition \ref{p:ConstructionFunction} actually exists.
Moreover, we provide bounds for $\upsilon$ and $C_G$.
For a compactly supported integrable function $\psi$, we denote its Fourier-Laplace transform as
	\[ \widehat{\psi}(z) = \int_{-\infty}^\infty \psi(\xi) e^{- i z \xi} \, \mathrm{d}\xi . \]

\begin{proposition}
\label{p:ConstructWindow}
Let $\psi \in C^1(\R)$ be even and nonnegative, with support in $[-R, R]$ for some $R > 0$, and such that $\int_{-\infty}^\infty \psi(x) \mathrm{d}x = 1$.
Put $\gamma = \max\{1, \int_{-\infty}^\infty |\psi'(x)| \mathrm{d}x\}$.
Then, for any $n \in \Z_+$,
	\[ G_n(z) = 1 - \left(1 - \widehat{\psi}\left(\sqrt{n} \gamma z\right)^2\right)^n   \]
is an entire function of exponential type at most $\nu_n = 2n^{3/2}\gamma R$ that is real-valued on $\R$ and satisfies conditions $(i)$ and $(ii)$ of Proposition \ref{p:ConstructionFunction}.
Moreover, for any $1 \leq \alpha \leq 2n$,
	\begin{equation}
		\label{eq:GnAlphaEst} 
		\sup_{0 < |y| \leq 1} \frac{\log_+ |G_n(x + iy)|}{|y|^{\alpha}} \leq \begin{cases} 4^{1 - \alpha} \nu^\alpha_n, & 1 \leq \alpha \leq 2 , \\ \frac14 \nu^\alpha_n , & \alpha > 2 \end{cases} . 
	\end{equation}
\end{proposition}

To prove Proposition \ref{p:ConstructWindow}, we will need the following bound.

\begin{lemma}
\label{l:SinPowEst}
For any $n \in \Z_+$ we have
	\[ \log \left| 1 - \sin^{2n}(z) \right| \leq \frac{(2n |y|)^{2n}}{4} , \qquad z = x + iy \in \C , \, |y| \leq \frac{1}{2n} . \]
\end{lemma}

\begin{proof}
When $n = 1$ we directly find, for any $z = x + iy \in \C$,
	\[ \log | 1 - \sin^2(z) | = \log |\cos z|^2 \leq 2 \log \cosh y \leq y^2 . \]
Suppose now $n \geq 2$ and fix $z = x + iy$ with $|y| \leq 1/(2n)$.
We have
	\begin{align*} 
		|1 - \sin^{2n}(z)|^2 
		&= 1 - 2 \re \sin^{2n}(z) + |\sin^{2n}(z)|^2 \\
		&= 1 - 2 |\sin^{2n}(z)| + |\sin^{2n}(z)|^2 + 4 (\im \sin^{n}(z))^2 \\
		&= (1 - |\sin^{2n}(z)|)^2 + 4 (\im \sin^{n}(z))^2 .  
	\end{align*}
Now
	\begin{align*}
		|\im \sin^{n}(z)|
		&=
		\frac12 |\sin^{n}(z) - \overline{\sin(z)}^{n}|
		= \frac12 \left| (\sin(z) - \overline{\sin(z)}) \sum_{k=0}^{n-1} (\sin z)^{n - 1 - k} \overline{\sin z}^k \right| \\
		&\leq n |\im \sin(z)| |\sin^{n-1}(z)|
		= n |\cos(x) \sinh(y)| |\sin^{n-1}(z)| .
	\end{align*}
Also, as $n y^2 \leq 1/ (4n) \leq 1/4$ and $e^t \leq 1 + 2t$ for $0 \leq t \leq 1/4$,
	\[
		|\sin(z)|^{2n}
		\leq \cosh^{2n}(y)
		\leq e^{n y^2}
		\leq 1 + 2 n y^2
		\leq 1 + 4 n^2 \sinh^2 y .
	\]
Combining everything, we get
	\begin{align*}
		|1 - \sin^{2n}(z)|^2 
		&= 1 - 2 |\sin^{2n}(z)| + |\sin^{4n}(z)| + 4 (\im \sin^{n}(z))^2 \\
		&\leq 1 - |\sin^{2n}(z)| + 4n^2 |\sin^{2n - 2}(z)| \sinh^2(y) (|\sin z|^2 + 1 - \sin^2 x) \\
		&\leq 1 + |\sin^{2n - 2}(z)| (n^2 \sinh^2(2y) - |\sin^2(z)|) .
	\end{align*}
If $n^2 \sinh^2(2y) \leq |\sin^2(z)|$, then the right-hand side is at most $1$, and our bound trivially holds.
Suppose now that $n^2 \sinh^2(2y) > |\sin^2(z)|$.
Then, since $s(t - s) \leq t^2 / 4$ for any $0 \leq s \leq t$,
	\begin{align*}
	 	|1 - \sin^{2n}(z)|^2 
		&\leq 1 + (n^2 \sinh^2(2y))^{n - 2} |\sin^2(z)| (n^2 \sinh^2(2y) - |\sin^2(z)|) \\
		&\leq 1 + \frac14 (n^2 \sinh^2(2y))^{n} .
	\end{align*}
Then, since $\sinh(t) \leq t e^{t^2 / 2}$ for any $t \geq 0$, using $|y| \leq 1/(2n)$,
	\begin{align*} 
		\log \left| 1 - \sin^{2n}(z) \right| 
		&\leq \frac12 \log\left(1 + \frac14 (n^2 \sinh^2(2y))^{n}\right) 
		\leq \frac18 (n^2 \sinh^2(2y))^{n} \\
		&\leq \frac{e^{4 n y^2}}{8} (2 n |y|)^{2n}
		\leq \frac{e^{1/n}}{8} (2n|y|)^{2n}
		\leq \frac{(2 n |y|)^{2n}}{4} . \qedhere
	\end{align*}
\end{proof}

\begin{proof}[Proof of Proposition \ref{p:ConstructWindow}]
For any $\lambda > 0$, we put
	\[ H_\lambda(z) = 1 - \left(1 - \widehat{\psi}\left(\lambda z\right)^2\right)^n . \]
By the classical Paley-Wiener theorem, we have that $\widehat{\psi}$ is an entire function of exponential type at most $R$.
Then $H_\lambda$ is an entire function of exponential type at most $2 n \lambda R$ for any $\lambda > 0$.
Moreover, since $\psi$ is real-valued and even, we see that each $H_\lambda$ is real-valued on $\R$.
Clearly $H_\lambda(0) = 1$ for any $\lambda > 0$, so $(i)$ of Proposition \ref{p:ConstructionFunction} holds for any $H_\lambda$.

We now verify that $(ii)$ of Proposition \ref{p:ConstructionFunction} holds for $H_\lambda$ when $\lambda$ is specifically chosen.
First, we claim that
	\begin{equation}
		\label{eq:0<LaplacepsiSquared<1} 
		0 \leq \widehat{\psi}(\lambda x)^2 < 1 , \qquad x \in \R \setminus \{0\} . 
	\end{equation}
This would already show $0 \leq H_\lambda(x) < 1$ for every $x \in \R \setminus \{0\}$ and $\lambda > 0$.
Indeed, since $\psi$ is real-valued on $\R$,
	\begin{align*} 
		\widehat{\psi}(x)^2 
		&= \re \widehat{\psi}(x)^2 
		= \re \int_{-\infty}^\infty \int_{-\infty}^\infty \psi(\xi) \psi(\eta) e^{-i(\xi + \eta) x}  \, \mathrm{d}\xi \mathrm{d}\eta \\
		&= \int_{-\infty}^\infty \int_{-\infty}^\infty \psi(\xi) \psi(\eta) \cos((\xi + \eta) x) \, \mathrm{d}\xi \mathrm{d}\eta . 
	\end{align*}
Let $I$ be some nontrivial open interval such that $\psi(\xi) \geq \varepsilon > 0$ for some $\varepsilon > 0$ and any $\xi \in I$.
Then, as $\psi$ is nonnegative, for $x \in \R \setminus \{0\}$,
	\begin{align*} 
		1 - \widehat{\psi}(x)^2 
		&= \widehat{\psi}(0)^2 - \widehat{\psi}(x)^2 
		= \int_{-\infty}^\infty \int_{-\infty}^\infty \left[1 - \cos((\xi + \eta) x)\right] \psi(\xi) \psi(\eta) \, \mathrm{d}\xi \mathrm{d}\eta \\
		&\geq \varepsilon^2 \int_I \int_I \left[1 - \cos((\xi + \eta) x)\right] \, \mathrm{d}\xi \mathrm{d}\eta 
		> 0 .
	\end{align*}
This shows \eqref{eq:0<LaplacepsiSquared<1}.

Since $\psi \in C^1(\R)$ we have $|\widehat{\psi}(x)| \leq \gamma |x|^{-1}$ for all $|x| \geq 1$.
If $n = 1$, then, for $\lambda |x| \geq 1$,
	\[ |H_\lambda(x)| = |\widehat{\psi}(\lambda x)^2| \leq \frac{\gamma^2}{\lambda^2} |x|^{-2} . \]
If $n \geq 2$, we similarly find by \eqref{eq:0<LaplacepsiSquared<1}, for $\lambda |x| \geq 1$,	
\begin{align*} 
		|H_\lambda(x)| 
		&= \left| 1 - (1 - \widehat{\psi}(\lambda x)^2)^n \right|
		= \left| \widehat{\psi}(\lambda x)^2 \left(1 + (1 - \widehat{\psi}(\lambda x)^2) + \cdots + (1 - \widehat{\psi}(\lambda x)^2)^{n-1} \right) \right| \\
		&\leq n |\widehat{\psi}(\lambda x)^2|
		\leq \frac{\gamma^2 n}{\lambda^2} |x|^{-2} .
	\end{align*}
Consequently, putting $\lambda_n = \gamma \sqrt{n}$, it follows that $|H_{\lambda_n}(x)| \leq |x|^{-2}$ for every $|x| \geq 1$.
Note that $G_n = H_{\lambda_n}$.

We finish by showing that $(iii)$ of Proposition \ref{p:ConstructionFunction} holds for $1 \leq \alpha \leq 2n$ for $G_n$ and that \eqref{eq:GnAlphaEst} is true.
Define
	\[ Q_n(z) = 2 G_n(z) - 1 , \qquad n \in \Z_+ . \]
Then each $Q_n$ is an entire function of exponential type at most $\nu_n = 2 n^{3/2} \gamma R$.
Moreover, they are real-valued on $\R$ and satisfy $|Q_n(x)| \leq 1$ for $x \in \R$.
By a result of Duffin and Schaeffer \cite[Theorem 2]{D-S-SomePropFuncExponentialType}, we find the bound
	\[ |Q_n(z)| \leq \cosh\left( \nu_n |y| \right) , \qquad z = x + iy \in \C , \]
so
	\[ |G_n(z)| \leq \frac{\cosh\left( \nu_n |y| \right) + 1}{2} = \cosh^2\left( \frac{\nu_n}{2} |y| \right) , \qquad z = x + iy \in \C . \]
As $\log \cosh t \leq \min\{ t,  t^2 / 2\}$, we get
	\begin{equation}
		\label{eq:EstGnAlpha=2} 
		\log |G_n(z)| \leq \min \left\{ \nu_n |y|,  \frac{\nu^2_n |y|^2}{4}\right\} \leq 4^{1-\alpha} \nu^\alpha_n |y|^\alpha , \qquad 1 \leq \alpha \leq 2 . 
	\end{equation}
This already shows that \eqref{eq:GnAlphaEst} holds for any $1 \leq \alpha \leq 2$ and $n \in \Z_+$.

Suppose now we have some fixed $y \in \R$ such that 
	\begin{equation}
		\label{eq:ExpTypeTimes|y|}
		\nu_n |y| \leq 1.
	\end{equation}
By our assumptions on $\psi$, we have, for any $x \in \R$,
	\begin{align*}
		|1 \pm \widehat{\psi}(x + iy)| 
		&\leq \int_{-\infty}^\infty \psi(\xi) |1 \pm \cos(\xi (x + iy))| \, \mathrm{d}\xi \\
		&= \int_{-\infty}^\infty \psi(\xi) (\cosh(\xi y) \pm \cos(\xi x)) \, \mathrm{d}\xi \\
		&= \widehat{\psi}(iy) \pm \widehat{\psi}(x) .
	\end{align*}
Thus, we find the bound
	\[ |1 + \widehat{\psi}(x + iy)| + |1 - \widehat{\psi}(x + iy)| \leq 2 \widehat{\psi}(iy) = 2 \int_{-R}^R \psi(\xi) \cosh(\xi y) \, \mathrm{d}\xi \leq 2 \cosh(R |y|) . \]
Consequently, the image of $\widehat{\psi}(\cdot + iy)$ is contained in the ellipse with boundary
	\[ \{ \cos(\theta + i R |y|) \mid \theta \in [0, 2\pi] \} . \]
Therefore, for any polynomial $P$, by the maximum-modulus principle, we have
	\[ |P(\widehat{\psi}(x + iy))| \leq \sup_{\theta \in [0, 2\pi]} |P(\cos(\theta + i R |y|))| , \qquad x \in \R . \]
For the particular choices of $P(t) = 1 - (1 - t^2)^n$, $n \in \Z_+$, this yields
	\[ |G_n(x + iy)| \leq \sup_{\theta \in [0, 2\pi]} |1 - \sin^{2n}(\theta + i \sqrt{n} \gamma R |y|)| . \]
By \eqref{eq:ExpTypeTimes|y|} we have that $\sqrt{n} \gamma R |y| \leq 1/(2n)$.
Hence, Lemma \ref{l:SinPowEst} gives
	\begin{equation}
		\label{eq:EstGnAlpha>=2} 
		\log |G_n(x + iy)| \leq \frac{(2n^{\frac32} \gamma R |y|)^{2n}}{4} . 
	\end{equation}

We may now show \eqref{eq:GnAlphaEst} for any $2 \leq \alpha \leq 2n$.
First, if \eqref{eq:ExpTypeTimes|y|} is true, then \eqref{eq:EstGnAlpha>=2} gives
	\[ \log |G_n(z)| \leq \frac14 \left( \nu_n |y| \right)^{2n} \leq\frac14 \nu_n^\alpha |y|^\alpha .  \]
If \eqref{eq:ExpTypeTimes|y|} does not hold, then \eqref{eq:EstGnAlpha=2} gives
	\[ \log |G_n(z)| \leq \frac14 \nu_n^2 |y|^2 \leq \frac14 \nu^{\alpha}_n |y|^\alpha  . \]
We have verified \eqref{eq:GnAlphaEst}.
\end{proof}

\subsection{Proof of the main theorem and lower bounds for the Paley-Wiener threshold}
\label{subsec:ProofsMainThm}

By all our previous results, we may now directly establish the main theorem of this section.

\begin{proof}[Proof of Theorem \ref{t:CharNonTrivEntireWeightFunc}]
By Corollary \ref{c:PaleyWienerThresholdBelowPi/2}, we know that necessarily $\Psi_{(\omega, \varpi)} \leq \pi / 2$.
On the other hand, combining Lemma \ref{l:PsiExplicit} with Propositions \ref{p:ConstructionFunction} and \ref{p:ConstructWindow} yields $\Psi_{(\omega, \varpi)} > 0$.
\end{proof}

We now establish explicit lower bounds for the Paley-Wiener threshold.

\begin{theorem}
\label{t:PaleyWienerThresholdLowerBound}
Let $\omega$ be a moderate quasianalytic weight function.
Suppose $n \geq 1$ is such that $\alpha(\omega) < 2n$.
Then
	\begin{equation} 
		\label{eq:PsiLowerBoundGeneral}
		\Psi_{(\omega, \varpi)} \geq \frac{1}{e \mathfrak{c}_+} \sup_{0<\theta<1} \frac{\theta}{\underline{\mathfrak{m}}_{\frac{1}{1-\theta}}(\omega) \underline{\mathfrak{m}}_{\theta K_\omega}(\omega)} 
	\end{equation}
where $\mathfrak{c}_+ = \max\{1, \mathfrak{c}(\omega, \varpi)\}$, and $K_\omega$ is defined depending on $n$ as follows.
	\begin{itemize}
		\item[$(i)$] When $n = 1$, 
			\begin{equation}
				\label{eq:PsiLowerBoundAlpha<2General}
				K_\omega = \frac{2 + 2\mathfrak{c}(\omega, \varpi) (P_2(\omega) - 1)}{{e \mathfrak{c}_+}} + \frac{\mathfrak{c}(\omega, \varpi)}{\mathfrak{c}_+} .
			\end{equation}		
		Moreover, if $\omega$ is convex we have
			\begin{equation}
				\label{eq:PsiLowerBoundAlpha<2} 
				\Psi_{(\omega, \varpi)} \geq \frac{1}{e} \sup_{0 < \theta < 1} \frac{\theta}{\underline{m}_{\frac{1}{2(1 - \theta)}}(\omega) \underline{\mathfrak{m}}_{\theta K_{1, \omega}}(\omega)} 
			\end{equation}
		with
			\begin{equation}
				\label{eq:K1omega} 
				K_{1, \omega} = 2 - \frac{2(e - 1)}{\overline{\mathfrak{m}}_{e}(\omega)} + \frac{\mathfrak{c}(\omega, \varpi)}{e}[2 P_2(\omega) - 2e - 1]_+ ;
			\end{equation}
		\item[$(ii)$] When $n \geq 2$, 
			\begin{equation}
				\label{eq:PsiLowerBoundAlpha>=2General}
				K_\omega = \frac{2 n^{3/2}}{e \mathfrak{c}_+} + \frac{\mathfrak{c}(\omega, \varpi)}{e \mathfrak{c}_+} \left( 2 n^{\frac32} - 2 + \frac14 \inf_{\max\{2, \alpha(\omega)\} < \alpha \leq 2n} \left[ 2^\alpha n^{\frac{3 \alpha}{2}} (\alpha P_\alpha(\omega) - \alpha) \right] \right) .
			\end{equation}	
		Moreover, if $\omega$ is convex we have
			\begin{equation}
				\label{eq:PsiLowerBoundAlpha>=2}
				\Psi_{(\omega, \varpi)} \geq \frac{1}{\mathfrak{c}_+} \sup_{0<\theta<1} \frac{\theta}{\underline{m}_{\frac{1}{2(1 - \theta)}}(\omega) \underline{\mathfrak{m}}_{\theta K_{2, \omega}}(\omega)} , 
			\end{equation}
		with
			\begin{multline}
				\label{eq:K2omega}
				K_{2, \omega} = \frac{2 n^{3/2}}{\mathfrak{c}_+} + \frac{\mathfrak{c}(\omega, \varpi)}{\mathfrak{c}_+}\left(4 n^{\frac32} - 2  \right. \\ \left. + \frac14 \inf_{\max\{2, \alpha(\omega)\} < \alpha \leq 2n} \left[ 2^\alpha n^{\frac{3 \alpha}{2}} (\alpha P_\alpha(\omega) - \alpha - 1)\right] \right)  .
			\end{multline}
	\end{itemize}
\end{theorem}

\begin{proof}
We start with the following observation: There exists a sequence $(\psi_k)_{k \in \Z_+}$ of even functions in $C^1(\R)$ such that $\psi_k$ has support in $[-(1+1/k), 1 + 1/k]$ and $\int_{-\infty}^\infty |\psi'_k(x)| \, \mathrm{d}x = 1$.
Indeed, take some nonnegative even $\varphi \in C^1(\R)$ with support contained in $(-1, 1)$ and with $\int \varphi = 1$.
Put $\varphi_k(x) = k \varphi(k x)$ for $k \in \Z_+$.
Then we set $\psi_k = \frac12 (1_{[-1, 1]} * \varphi_k)$.
Clearly the support of $\psi_k$ is contained in $(-(1 + 1/k), 1+ 1/k)$, while $\psi'_k(x) = (\varphi_k(x+1) - \varphi_k(x-1)) / 2$.
The terms have disjoint support, so we get
	\[ \int_{-\infty}^\infty |\psi'_k(x)| \, \mathrm{d}x = \frac12 \int_{-\infty}^\infty \varphi_k(x + 1) \, \mathrm{d}x + \frac12 \int_{-\infty}^\infty \varphi_k(x - 1) \, \mathrm{d}x = 1 . \]
Therefore, when applying Proposition \ref{p:ConstructWindow}, we may assume $R \to 1^+$ and $\gamma = 1$.
In particular, we may take $\upsilon_n$ arbitrarily close to $2 n^{3/2}$.
We now consider the two cases separately.

$(i)$ We may take $\upsilon_1$ arbitrarily close to $2$.
Moreover, for any $\max\{1, \alpha(\omega)\} < \alpha \leq 2$, by \eqref{eq:GnAlphaEst}, we may construct $G_1$ as in Proposition \ref{p:ConstructWindow} such that $C_{G_1}$ as in Proposition \ref{p:ConstructionFunction} with respect to $\alpha$ is arbitrarily close to $2^{2 - \alpha}$ and thus below $\upsilon_1$.
As $\alpha P_\alpha(\omega) - \alpha$ is nonincreasing by Lemma \ref{l:(aPa - a)nonincreasing}, combining Lemmas \ref{l:Continuitymathfrakm} and \ref{l:PsiExplicit} with \eqref{eq:LowerBoundThresholdExplicitConstructGeneral} and \eqref{eq:K}, respectively \eqref{eq:LowerBoundThresholdExplicitConstruct1} and \eqref{eq:K1} when $\omega$ is convex, yields the estimates.

$(ii)$ For any $\max\{2, \alpha(\omega)\} < \alpha \leq 2n$, by \eqref{eq:GnAlphaEst}, we may construct $G_n$ as in Proposition \ref{p:ConstructWindow} such that $C_{G_n}$ as in Proposition \ref{p:ConstructionFunction} with respect to $\alpha$ is arbitrarily close to $2^{\alpha - 2} n^{\frac32 \alpha}$.
Therefore, we may have $C_{G_n} \geq \upsilon_n$.
Then the estimates follow from Lemmas \ref{l:Continuitymathfrakm} and \ref{l:PsiExplicit} combined with \eqref{eq:LowerBoundThresholdExplicitConstructGeneral} and \eqref{eq:K}, respectively \eqref{eq:LowerBoundThresholdExplicitConstruct2} and \eqref{eq:K2} when $\omega$ is convex.
\end{proof}

We finish by considering the specific case when $\omega$ is regularly varying.
In this case, we may determine the threshold exactly, instead of finding bounds through Theorem \ref{t:PaleyWienerThresholdLowerBound}.

\begin{theorem}
\label{t:LowBoundPaleyWienerThresholdRegVar}
Let $\omega$ be a convex weight function that is regularly varying of index $p \geq 1$.
	\begin{itemize}
		\item[$(i)$] When $p = 1$, then
			\[ \Psi_{(\omega, \varpi)} = \frac{\pi}{2} . \]
			
		\item[$(ii)$] When $p > 1$
			\[ \Psi_{(\omega, \varpi)} = \frac{1}{(p - 1)^{\frac{1}{p}}} \min\left\{ \tan \frac{\pi}{2 p} , \cot \frac{\pi}{2 p} \right\} . \]
	\end{itemize}
\end{theorem}

\begin{proof}
$(i)$
By Theorem \ref{t:CharNonTrivEntireWeightFunc} we have that $\Psi_{(\omega, \varpi)} \leq \pi/2$.
So, it suffices to show $\Psi_{(\omega, \varpi)} \geq \pi/2$.
We have $\omega(t) = o(\varpi(t))$ and $\varpi$ is regularly varying of index 1 by Proposition \ref{p:CharCondDomVarpi} and Lemma \ref{l:RegVarEquivOverlineOmega}.
We adapt the proof of Proposition \ref{p:ConstructionFunction} in this specific case.
Let again $\mu_p = \inf\{t \geq 0 \mid \mu(t) \geq p\}$ for $p \in \Z_+$.
Without loss of generality, we may assume $\mu_1 = 1$.
We put
	\[ F(z) = \prod_{p=1}^\infty \sinc\left(\frac{z}{\mu_p}\right) , \qquad z \in \C , \] 
with $\sinc z = \sin(z) / z = (1/2) \int_{-1}^1 e^{i z t} \, \mathrm{d}t$.
As in the proof of Proposition \ref{p:ConstructionFunction}, one shows that $F$ is an entire function.

We now give two estimates for $F$.
Let $\varepsilon > 0$ be fixed but arbitrary.
We first consider $F$ on $\R$.
By \eqref{eq:muBoundomega} and \eqref{eq:OmegaBoundMuConvex} we have, for any $\lambda > 1$,
	\[ 1 \leq \frac{\mu(t)}{\omega(t)} \leq \frac{\omega(\lambda t) - \omega(t)}{\omega(t) \log \lambda} \to \frac{\lambda - 1}{\log \lambda} , \qquad t \to \infty . \]
Letting $\lambda \to 1^+$, we see that $\omega(t) \sim \mu(t)$.
Take any $N \in \Z_+$.
Let $x \in \R$ with $|x|$ sufficiently large be fixed and set $m = \lfloor \omega(|x|) \rfloor$.
As $|\sinc t| \leq 1$ for any $t \in \R$, we have
	\[ |F(x)|^{1/m} \leq \prod_{p=1}^{N m} \left| \sinc\left(\frac{x}{\mu_p}\right) \right|^{\frac{1}{m}} . \]
Taking the logarithm, we find
	\begin{align*} 
		\frac{1}{m} \log |F(x)| 
		&\leq \sum_{p=1}^{Nm} \frac{1}{m} \log\left|\sinc\left(\frac{x}{\mu_p}\right)\right| 
		= \sum_{p=1}^{Nm} \int_{(p-1) / m}^{p/m} \log\left|\sinc\left(\frac{x}{\mu_p}\right)\right| \, \mathrm{d}t \\
		&= \int_0^N \log\left|\sinc\left(\frac{x}{\mu_{\lceil t m \rceil}}\right)\right| \, \mathrm{d}t .
	\end{align*}
Since $\mu(t |x|) / \omega(|x|) \to t$, also $|x| / \mu_{\lceil t m \rceil} \to 1/t$ as $|x| \to \infty$.
Thus, by Fatou's lemma,
	\[ \liminf_{|x| \to \infty} \frac{- \log |F(x)|}{\omega(|x|)} \geq - \int_0^N \log\left|\sinc\left(\frac{1}{t}\right)\right| \, \mathrm{d}t = 
	- \int_{1/N}^\infty \frac{\log |\sinc(s)|}{s^2} \, \mathrm{d}s . \]
As $N$ was chosen arbitrarily, letting it tend to infinity, we find
	\[ \liminf_{|x| \to \infty} \frac{- \log |F(x)|}{\omega(|x|)} \geq - \int_0^\infty \frac{\log |\sinc(s)|}{s^2} \, \mathrm{d}s = \frac{\pi}{2} .  \]
We conclude
	\begin{equation} 
		\label{eq:FBoundRRegVar1}
		\log |F(x)| \leq - \left(\frac{\pi}{2} - \varepsilon\right) \omega(|x|) + O(1) , \qquad x \in \R . 
	\end{equation}

Next, we give a general estimate for $F$ on $\C$. 
Since $|\sinc(z)| \leq \sinh(|y|) / |y|$ for $z = x + iy \in \C$, we have
	\[ \log |F(z)| \leq \sum_{p = 1}^\infty \log \frac{\sinh(|y| / \mu_p)}{|y| / \mu_p} , \qquad z = x + iy \in \C . \]
When $t \geq 0$, then $\sinh(t) / t \leq e^t$.
Also, for $0 \leq t \leq 1$, we have
	\[ \log \frac{\sinh t}{t} = \log\left(1 + \sum_{k=1}^\infty \frac{t^{2k}}{(2k + 1)!}\right) \leq \log\left(1 + t^2 \sum_{k=1}^\infty \frac{1}{(2k + 1)!} \right) \leq (\sinh 1 - 1) t^2 . \]
Thus
	\[ \log |F(z)| \leq |y| \sum_{\mu_p \leq |y|} \frac{1}{\mu_p} + (\sinh 1 - 1) |y|^2 \sum_{\mu_p > |y|} \frac{1}{\mu_p^2} . \] 
We estimate each term individually.
Let $M$ be as in \eqref{eq:M(t)}.
By integration by parts, using \eqref{eq:muBoundomega}, we have, when $|y|$ is large enough,
	\begin{align*} 
		|y| \sum_{\mu_p \leq |y|} \frac{1}{\mu_p} 
		&= |y| \int_1^{|y|} \frac{\mathrm{d} M(u)}{u} 
		= M(|y|) + |y| \int_{1}^{|y|} \frac{M(u)}{u} \frac{\mathrm{d}u}{u} \\
		&\leq \mu(|y|) + \omega(|y|) + \varpi(|y|) 
		\leq \left(1 + \frac{\varepsilon}{2}\right) \varpi(|y|) . 
	\end{align*}
Also, by Karamata's theorem \cite[Theorem 1.5.11(ii), p.~28]{B-G-T-RegVar} and \eqref{eq:muBoundomega}, when $|y|$ is large enough,
	\begin{align*}
		|y|^2 \sum_{\mu_p > |y|} \frac{1}{\mu_p^2}
		&= |y|^2 \int_{|y|^+}^\infty \frac{\mathrm{d}M(u)}{u^2}
		= - M(|y|) + 2 |y|^2 \int_{|y|}^\infty \frac{M(u)}{u^3} \, \mathrm{d}u \\
		&\leq \left(1 + \frac{\varepsilon}{4}\right) \mu(|y|) 
		\leq \frac{\varepsilon}{2} \varpi(|y|) .
	\end{align*}
We conclude
	\begin{equation} 
		\label{eq:FBoundCRegVar1}
		\log |F(z)| \leq (1 + \varepsilon) \varpi(|y|) + O(1) , \qquad z = x+iy \in \C . 
	\end{equation}

Applying Lemma \ref{l:InterpolaEst} with \eqref{eq:FBoundRRegVar1} and \eqref{eq:FBoundCRegVar1}, with the weight functions $\left(\frac{\pi}{2} - \varepsilon\right)\omega$ and $(1+\varepsilon)\varpi$, we get that, for any $\theta \in (0, 1)$,
	\[ |F(z)| \lesssim e^{\theta (1 + \varepsilon) \varpi\left(\frac{|y|}{\theta}\right) - (1 - \theta) \left(\frac{\pi}{2} - \varepsilon\right) \omega\left(|x|\right)} , \qquad z = x + iy \in \C . \] 
As $\omega$ and $\varpi$ are regularly varying of index 1, we see that
	\[ |F(z)| \lesssim e^{(1 + \varepsilon)^2 \varpi\left(|y|\right) - (1 - \theta) (1 - \varepsilon) \omega\left(\frac{\pi}{2} |x|\right)} , \qquad z = x + iy \in \C . \] 
By letting $\varepsilon, \theta \to 0^+$, we conclude that $\Psi_{(\omega, \varpi)} \geq \pi/2$.

$(ii)$
Using Lemma \ref{l:RegVarEquivOverlineOmega}, this follows from \cite[Corollary 1]{S-ClassEntFuncRapDecRealAxis}.
\end{proof}

\section{An optimal Beurling-H\"ormander-type uncertainty principle}
\label{sec:UncertaintyPrinciples}

The primary aim of this section is to show Theorem \ref{t:MainUncertaintyPrincipleIntro}, that $(\omega, \varpi^*)$ is an optimal uncertainty pair.
Before continuing, let us first specify the multidimensional extension of Definition \ref{def:UncertaintyPair}.

\begin{definition}
\label{def:UncertaintyPairMultiDim}
Let $\omega, \eta$ be weight functions.
We call $(\omega, \eta)$ an \emph{uncertainty pair} if there exists a constant $\Phi_{(\omega, \eta)} \in [0, \infty)$ such that for any $a, b > 0$ the following hold:
	\begin{itemize}
		\item[(I)] If $a b > \Phi_{(\omega, \eta)}$, should $f \in L^2(\R^d)$ satisfy, for some $j \in \{1, \ldots, d\}$,
			\begin{equation}
				\label{eq:UncertaintyDecay}
		 		\int_{\R^d} |f(x)| e^{\omega(a |x_j|)} \, \mathrm{d}x < \infty \quad \text{and} \quad \int_{\R^d} |\widehat{f}(\xi)| e^{\eta(b |\xi_j|)} \, \mathrm{d}\xi < \infty 
			\end{equation}
		then $f = 0$.
		\item[(II)] If $a b < \Phi_{(\omega, \eta)}$, there exists a dense subspace of $L^2(\R^d)$ consisting of functions $f$ satisfying 
			\begin{equation}
				\label{eq:UncertaintyDecayDenseSubsp}
		 		\int_{\R^d} |f(x)| e^{\sum_{j=1}^d \omega(a |x_j|)} \, \mathrm{d}x < \infty \quad \text{and} \quad \int_{\R^d} |\widehat{f}(\xi)| e^{\sum_{j=1}^d \eta(b |\xi_j|)} \, \mathrm{d}\xi < \infty .
			\end{equation}
	\end{itemize}
We call $\Phi_{(\omega, \eta)}$ the \emph{uncertainty threshold} of $(\omega, \eta)$.
The uncertainty pair $(\omega, \eta)$ is called \emph{optimal} if $\Phi_{(\omega, \eta)} > 0$.
\end{definition}

We may now reformulate Theorem \ref{t:MainUncertaintyPrincipleIntro} as follows.

\begin{theorem}
\label{t:MainUncertaintyPrincipleGeneral}
Let $\omega$ be a moderate quasianalytic weight function.
Then $(\omega, \varpi^*)$ is an optimal uncertainty pair.
\end{theorem}

As a consequence of Theorem \ref{t:MainUncertaintyPrincipleGeneral} and Propositions \ref{p:SuffNeccCondEquivVarpi} and \ref{p:CharCondDomVarpi}, we get the following sufficient conditions for optimality and nonoptimality of the uncertainty pair $(\omega, \omega^*)$, that is, the Beurling-H\"ormander uncertainty principle \eqref{eq:BeurlingHormanderUncertainty}.
In particular, Theorem \ref{t:MainUncertaintyPrincipleRegVar} is a direct consequence of it in view of Lemma \ref{l:RegVarEquivOverlineOmega}.

\begin{theorem}
\label{t:SuffNeccCondBeurlingHormanderOptimal}
Let $\omega$ be a moderate convex weight function such that $t = o(\omega(t))$.
Then $(\omega, \omega^*)$ is an uncertainty pair with $\Phi_{(\omega, \omega^*)} \leq 1$.
	\begin{itemize}
		\item[1.] If $\beta(\omega) > 1$, then $(\omega, \omega^*)$ is optimal.
		\item[2.] If $\omega$ is geometrically persistent, then $(\omega, \omega^*)$ is not optimal.
	\end{itemize}
\end{theorem}

We now move on to prove Theorem \ref{t:MainUncertaintyPrincipleIntro}.
The idea is to translate it into Theorem \ref{t:CharNonTrivEntireWeightFunc}.
We need several intermediate steps.
Our first observation is the following explicit form of the uncertainty threshold.
In particular, it shows that it is independent of the dimension.

\begin{lemma}
\label{l:PhiExplicit}
Let $\omega, \eta$ be weight functions such that $(\omega, \eta)$ is an uncertainty pair.
Then
	\begin{equation}
		\label{eq:PhiExplicit} 
		\Phi_{(\omega, \eta)} = \sup \left[ \{ ab \mid \text{$\exists f \in L^2(\R) \setminus \{0\}$ satisfying \eqref{eq:UncertaintyDecayDenseSubsp} for $a$ and $b$} \} \cup \{0\} \right] .
	\end{equation}
Conversely, if the right-hand side of \eqref{eq:PhiExplicit} is finite, then $(\omega, \eta)$ is an uncertainty pair in any dimension for which the threshold satisfies \eqref{eq:PhiExplicit}.
\end{lemma}

\begin{proof}
If $f \in L^2(\R)$ satisfies \eqref{eq:UncertaintyDecayDenseSubsp} for certain $a_0, b_0 > 0$, and if $a_1, b_1 > 0$ are such that $a_0 b_0 = a_1 b_1$, then $g(x) = f\left(\frac{a_1}{a_0} x\right)$ satisfies \eqref{eq:UncertaintyDecayDenseSubsp} for $a_1$ and $b_1$.
Therefore, the supremum in \eqref{eq:PhiExplicit}, which we denote by $\Phi_{(\omega, \eta)}'$, is well-defined.

Suppose $f \in L^2(\R^d) \setminus \{0\}$ is such that \eqref{eq:UncertaintyDecay} holds for certain $ab > \Phi_{(\omega, \eta)}'$ and some $j \in \{1, \ldots, d\}$.
We revert to the one-dimensional case using a trick from \cite{B-R-UncertaintyPrinciplesLieGroups}.
Without loss of generality, we may assume that $j = 1$.
Then there exists some $x' \in \R^{d-1}$ such that $g(x) = f(x, x') \in L^2(\R) \setminus \{0\}$ and
	\[ \int_{-\infty}^\infty |g(x)| e^{\omega(a |x|)} \, \mathrm{d}x < \infty . \]
We have
	\[ \widehat{g}(\xi) = \frac{1}{(2\pi)^{d-1}} \int_{\R^{d-1}} \widehat{f}(\xi, \xi') e^{i \ev{x'}{\xi'}} \, \mathrm{d}\xi' , \]
so
	\[ \int_{-\infty}^\infty |\widehat{g}(\xi)| e^{\eta(b |\xi|)} \, \mathrm{d} \xi < \infty . \]
This contradicts the definition of $\Phi_{(\omega, \eta)}'$, so such $f$ cannot exist. 
Hence $\Phi'_{(\omega, \eta)} \geq \Phi_{(\omega, \eta)}$.

Suppose now that $0 < a < a_0$ and $0 < b < b_0$ are such that $a_0 b_0 < \Phi_{(\omega, \eta)}'$.
Then, there exists some $g_0 \in L^2(\R) \setminus \{0\}$ satisfying \eqref{eq:UncertaintyDecayDenseSubsp} for $a_0$ and $b_0$.
Putting $f_0 = g_0 \otimes \cdots \otimes g_0$, where we take $d$ tensor products, then $f_0 \in L^2(\R^d) \setminus \{0\}$ and $f_0$ also satisfies \eqref{eq:UncertaintyDecayDenseSubsp} in dimension $d$ for $a_0$ and $b_0$.
Note that necessarily $f_0 \in \S(\R^d)$ \cite{C-C-K-CharSchwartz}.
For any $\lambda = (t, \zeta) \in \R^{2d}$, we consider the time-frequency shift
	\[ [\pi(\lambda) f_0](x) := f_0(x - t) e^{i \ev{x}{\zeta}} , \qquad x \in \R^d . \]
Using the convexity of $t \mapsto \omega(e^t)$ and $t \mapsto \eta(e^t)$, for any fixed $x, \xi \in \R^d$, we may find a $C = C_{x, \xi} \geq 1$ such that $\omega(a |x_j|) \leq \omega(a_0 |x_j - t_j|) + \log C$ and $\eta(b |\xi_j|) \leq \eta(b_0 |\xi_j - \zeta_j|) + \log C$.
We calculate
	\[
		\int_{\R^d} |[\pi(\lambda) f_0](x)| e^{\sum_{j=1}^d \omega(a |x_j|)} \, \mathrm{d}x 
		\leq C^d \int_{\R^d} |f_0(x - t)| e^{\sum_{j=1}^d \omega(a_0 |x_j - t_j|)} \, \mathrm{d}x  
		< \infty
	\]
and	
	\[
		\int_{\R^d} |\mathcal{F}\{\pi(\lambda) f_0\}(\xi)| e^{\sum_{j=1}^d \eta(b |\xi_j|)} \, \mathrm{d}\xi 
		\leq C^d \int_{\R^d} |\widehat{f}_0(\xi - \zeta)| e^{\sum_{j=1}^d \eta(b_0 |\xi_j - \zeta_j|)} \, \mathrm{d}\xi  
		< \infty .
	\]
We conclude that the elements in the linear span of $\{ \pi(\lambda) f_0 \mid \lambda \in \R^{2d} \}$ all satisfy \eqref{eq:UncertaintyDecayDenseSubsp} for $a$ and $b$.
This is a dense subspace of $\S(\R^d)$, as otherwise the short-time Fourier transform with window $f_0$ would be non-injective on $\S^\prime(\R^d)$ contradicting \cite[Corollary 11.2.7, p.~229]{G-FoundTFAnal}. 
Therefore $\Phi'_{(\omega, \eta)} \leq \Phi_{(\omega, \eta)}$, so they must coincide.
\end{proof}

\begin{remark}
As follows from the proof of Lemma \ref{l:PhiExplicit}, in Definitions \ref{def:UncertaintyPair} and \ref{def:UncertaintyPairMultiDim} we may exchange ``dense subspace of $L^2(\R^d)$'' with ``dense subspace of $\S(\R^d)$'' to get an equivalent definition.
\end{remark}

When the frequency bounds dominate $t$, then our functions allow entire extensions satisfying Paley--Wiener-type bounds as in Section \ref{sec:EntireFunc}.

\begin{lemma}
\label{l:TFDecayAndEntireFunc}
Let $\omega, \eta$ be weight functions such that $\omega$ is moderate, $\eta$ is convex, $\omega \preceq \eta$, and $t = o(\eta(t))$.
\begin{itemize}
	\item[$(i)$] If $f \in L^2(\R)$ satisfies, for certain $a, b > 0$,
		\begin{equation} 
			\label{eq:TFDecayToEntireFunc}
			\int_{-\infty}^\infty |f(x)| e^{\omega(a |x|)} \, \mathrm{d}x < \infty \quad \text{and} \quad \int_{-\infty}^\infty |\widehat{f}(\xi)| e^{\eta^*(b |\xi|)} \, \mathrm{d}\xi < \infty , 
		\end{equation}
	then $f$ extends to an entire function satisfying, for any $\theta \in (0, 1)$ and for any $H_0 > \underline{\mathfrak{m}}_{\frac{2}{(1 + \kappa_\omega)(1 - \theta)}}(\omega)$ and $H_1 > \underline{\mathfrak{m}}_{\theta}(\eta)$,
		\[
			|f(z)| \lesssim \exp\left[ \eta\left(\frac{H_1}{\theta b} |y| \right) - \omega\left(\frac{a}{H_0} |x|\right) \right] , \qquad z = x + iy \in \C ,
		\]
	with $\kappa_\omega$ as in \eqref{eq:IsConvex}.
		
	\item[$(ii)$] If $f \in L^2(\R)$ extends to an entire function satisfying
		\begin{equation}
			\label{eq:EntireFuncToTFDecay}
			|f(z)| \lesssim \exp\left[ \eta\left(\frac{|y|}{b}\right) - \omega(a |x|) \right] , \qquad z = x + iy \in \C ,
		\end{equation}
	then $f$ satisfies \eqref{eq:TFDecayToEntireFunc} for any $0 < a' < a$ and $0 < b' < b$.
\end{itemize}
\end{lemma}

\begin{proof}
$(i)$ 
Suppose $f \in L^2(\R) \setminus \{0\}$ is such that \eqref{eq:TFDecayToEntireFunc} holds for $a, b > 0$.
Fix any $k \in \Z_+$.
Take a parametrix $\phi \in L^1(\R)$ and $\psi \in C^\infty(\R)$, having compact support in $[-1, 0]$, such that $\delta = \phi^{(k)} + \psi$.
For any $\varepsilon > 0$, put $\phi_\varepsilon = \varepsilon^{k-1} \phi(\cdot / \varepsilon)$ and $\psi_\varepsilon = \varepsilon^{-1} \psi(\cdot / \varepsilon)$.
By the Riemann-Lebesgue lemma and the assumption on $\widehat{f}$, we see that $f$ and all its derivatives vanish at infinity.
Then
	\[ f = \phi_\varepsilon * f^{(k)} + \psi_\varepsilon * f , \qquad \varepsilon > 0 . \]
Therefore, for $x \geq 0$,
	\[ |f(x)| \leq \varepsilon^k \|\phi\|_{L^1} \|f^{(k)}\|_{L^\infty} + \frac{\|\psi\|_{L^\infty}}{\varepsilon} \int_{x}^{x + \varepsilon} |f(t)| \, \mathrm{d}t \lesssim \varepsilon^k + \varepsilon^{-1} e^{- \omega(a |x|)} ,  \]
where the bound does not depend on $\varepsilon$.
Then, for any $x \geq 0$, choosing $\varepsilon = e^{-\omega(a |x|) / (k + 1)}$, it follows that
	\[ |f(x)| \lesssim e^{- \frac{k}{k + 1} \omega(a |x|)} . \]
The same bound for $x < 0$ can be found by reflecting the parametrix about the origin.

As $t = o(\eta^*(t))$, we may extend $f$ to an entire function as
	\[ f(z) = \frac{1}{2\pi} \int_{-\infty}^\infty \widehat{f}(\xi) e^{i \xi z} \, \mathrm{d}\xi , \qquad z \in \C . \]
Then, we find the bound, for $z = x + iy \in \C$,
	\begin{align*} 
		|f(z)| 
		&\lesssim \int_{-\infty}^\infty |\widehat{f}(\xi)| e^{- y \xi} \, \mathrm{d}\xi 
		\leq \left( \int_{-\infty}^\infty |\widehat{f}(\xi)| e^{\eta^*(b |\xi|)} \, \mathrm{d}\xi \right) \cdot \sup_{\xi \in \R} \exp[ |y \xi| - \eta^*(b |\xi|)]  
		\lesssim e^{\eta\left(\frac{|y|}{b}\right)} . 
	\end{align*}
Therefore, for $x \in \R$ and $z \in \C$, we have the bounds:
	\[ |f(x)| \lesssim e^{- \frac{k}{k+1} \omega\left(a |x|\right)} \qquad \text{and} \qquad |f(z)| \lesssim e^{\eta\left(\frac{|y|}{b}\right)} . \] 
Fix $\theta \in (0, 1)$.
Lemma \ref{l:InterpolaEst}, with the weight function $\frac{k}{k+1} \omega$ and $\eta$, yields
	\[ |f(z)| \lesssim \exp\left[ \theta \eta\left(\frac{|y|}{\theta b}\right) - \frac{(1 + \kappa_\omega)(1 - \theta)}{2} \frac{k}{k+1} \omega(a |x|) \right] , \qquad z = x + iy \in \C .   \]
Let $C_0, H_0 > 0$ be such that $\frac{2}{(1 + \kappa_\omega)(1 - \theta)} \frac{k+1}{k} \omega(t) \leq \omega(H_0 t) + C_0$.
Also, take any $C_1, H_1 > 0$ such that $\theta \eta(t) \leq \eta(H_1 t) + C_1$.
From here, it follows that
	\[ |f(z)| \lesssim \exp\left[ \eta\left(\frac{H_1}{\theta b} |y|\right) - \omega\left(\frac{a}{H_0} |x|\right) \right] . \]
Letting $k \to \infty$ and using Lemma \ref{l:Continuitymathfrakm}, the result follows.

$(ii)$ 
Fix any $0 < a' < a$.
Let $C > 0$ and $0 < H < 1$ be such that $\omega\left(\frac{a'}{a} t\right) \leq H \omega(t) + C$.
Then
	\[ \int_{-\infty}^\infty |f(x)| e^{\omega(a' |x|)} \, \mathrm{d}x \leq \int_{-\infty}^\infty |f(x)| e^{\omega(a |x|)} \cdot e^{-\left(1 - H\right) \omega(a |x|)} \mathrm{d}x < \infty . \]
On the other hand, by changing the contour
	\begin{multline*} 
		|\widehat{f}(\xi)| 
		= \left| \int_{-\infty}^\infty f(x) e^{- i x \xi} \, \mathrm{d}x \right| 
		= \inf_{y \in \R} \left| \int_{-\infty}^\infty f(x - iy) e^{- i x \xi - y \xi} \, \mathrm{d}x \right| \\
		\lesssim \left( \int_{-\infty}^\infty e^{- \omega(a |x|)} \, \mathrm{d}x \right) \cdot \inf_{y \in \R,  \, y\xi \geq 0} \exp\left[ \eta\left(\frac{|y|}{b}\right) - |y \xi| \right]
		\lesssim e^{- \eta^*(b |\xi|)} .
	\end{multline*}
Then, using the convexity of $\eta^*$, for any $0 < b' < b$,
	\[ \int_{-\infty}^\infty |\widehat{f}(\xi)| e^{\eta^*(b' |\xi|)} \, \mathrm{d}\xi \leq \int_{-\infty}^\infty |\widehat{f}(\xi)| e^{\eta^*(b |\xi|)} \cdot e^{-\left(1 - \frac{b'}{b}\right) \eta^*(b |\xi|)} \mathrm{d}\xi < \infty . \qedhere \]
\end{proof}

We are now sufficiently prepared to show Theorems \ref{t:MainUncertaintyPrincipleGeneral} and \ref{t:SuffNeccCondBeurlingHormanderOptimal}.

\begin{proof}[Proof of Theorem \ref{t:MainUncertaintyPrincipleGeneral}]
By combining Theorem \ref{t:CharNonTrivEntireWeightFunc} with Lemma \ref{l:TFDecayAndEntireFunc} it follows that $0 < \Phi_{(\omega, \varpi^*)} < \infty$.
\end{proof}

\begin{proof}[Proof of Theorem \ref{t:SuffNeccCondBeurlingHormanderOptimal}]
That $(\omega, \omega^*)$ is an uncertainty pair follows from the result of Beurling and H\"ormander \eqref{eq:BeurlingHormanderUncertainty}, see \cite{H-UniquenessThmBeurlingFourier}, and Lemma \ref{l:PhiExplicit}.

1. If $\beta(\omega) > 1$, then $\varpi \asymp \omega$ by Proposition \ref{p:SuffNeccCondEquivVarpi}. Combining Theorem \ref{t:MainUncertaintyPrincipleGeneral} with Lemma \ref{l:PhiExplicit} gives that $(\omega, \omega^*)$ is an optimal uncertainty pair.

2. If $\omega$ is geometrically persistent, then $\omega \triangleleft \varpi$ by Proposition \ref{p:CharCondDomVarpi}, so also $\varpi^* \triangleleft \omega^*$.
Therefore, if there would be some $f \in L^2(\R^d) \setminus \{0\}$ satisfying \eqref{eq:UncertaintyDecayDenseSubsp} for $\eta = \omega^*$ for certain $a_0, b_0$, then it would also satisfy \eqref{eq:UncertaintyDecayDenseSubsp} for $\eta = \varpi^*$ with $a = a_0$ and any $b > 0$.
But this would imply that $(\omega, \varpi^*)$ is not an optimal uncertainty pair, contradicting Theorem \ref{t:MainUncertaintyPrincipleGeneral}. 
Therefore, necessarily $\Phi_{(\omega, \omega^*)} = 0$, so that $(\omega, \omega^*)$ is not an optimal uncertainty pair.
\end{proof}

We end this section by considering upper and lower bounds for $\Phi_{(\omega, \varpi^*)}$.

\begin{theorem}
\label{t:UncertaintyThresholdBound}
Let $\omega$ be a moderate quasianalytic weight function.
Suppose $n \geq 1$ is such that $\alpha(\omega) < 2n$.
Then
	\[ \frac{1}{e \mathfrak{c}_+} \sup_{0<\theta<1} \frac{\theta}{\underline{\mathfrak{m}}_{\frac{1}{1-\theta}}(\omega) \underline{\mathfrak{m}}_{\theta K_\omega}(\omega)} \leq \Phi_{(\omega, \varpi^*)} \leq \underline{\mathfrak{m}}_{1/2}(\omega) \underline{\mathfrak{m}}_{4}(\omega) \pi , \]
where $\mathfrak{c}_+ = \max\{1, \mathfrak{c}(\omega, \varpi)\}$, and $K_\omega$ is defined depending on $n$ as follows.
	\begin{itemize}
		\item[$(i)$] When $n = 1$, then $K_\omega$ equals \eqref{eq:PsiLowerBoundAlpha<2General}. 
		Moreover, if $\omega$ is convex we have
			\begin{equation*}
				\frac{1}{e} \sup_{0<\theta<1} \frac{\theta}{\underline{\mathfrak{m}}_{\frac{1}{2(1 - \theta)}}(\omega) \underline{\mathfrak{m}}_{\theta K_{1, \omega}}(\omega)} \leq \Phi_{(\omega, \varpi^*)} \leq \underline{\mathfrak{m}}_{1/2}(\omega) \underline{\mathfrak{m}}_{2}(\omega) \pi 
			\end{equation*}
		where $K_{1, \omega}$ is as in \eqref{eq:K1omega}.
		
		\item[$(ii)$] When $n \geq 2$, then $K_\omega$ equals \eqref{eq:PsiLowerBoundAlpha>=2General}.
		Moreover, if $\omega$ is convex we have
			\begin{equation*}
				\frac{1}{\mathfrak{c}_+} \sup_{0<\theta<1} \frac{\theta}{\underline{\mathfrak{m}}_{\frac{1}{2(1 - \theta)}}(\omega)  \underline{\mathfrak{m}}_{\theta K_{2, \omega}}(\omega)} \leq \Phi_{(\omega, \varpi^*)} \leq \underline{\mathfrak{m}}_{1/2}(\omega) \underline{\mathfrak{m}}_{2}(\omega) \pi  , 
			\end{equation*}
		where $K_{2, \omega}$ is as in \eqref{eq:K2omega}.
	\end{itemize}
\end{theorem}

\begin{proof}
By Lemmas \ref{l:IndicesForVarpi} and \ref{l:TFDecayAndEntireFunc}, with $\theta = 1/2$, we see that
	\begin{equation}
		\label{eq:BoundsPsiPhi} 
		\Psi_{(\omega, \varpi)} \leq \Phi_{(\omega, \varpi^*)} \leq 2 \underline{\mathfrak{m}}_{1/2}(\omega) \underline{\mathfrak{m}}_{\frac{4}{1+\kappa_\omega}}(\omega)  \Psi_{(\omega, \varpi)} ,
	\end{equation}
with $\kappa_\omega$ as in \eqref{eq:IsConvex}.
The result now follows from Theorem \ref{t:PaleyWienerThresholdLowerBound}.
\end{proof}

When $\omega$ is regularly varying, we may find the following exact values for the uncertainty threshold.

\begin{theorem}
\label{t:LowBoundUncertaintyThresholdRegVar}
Let $\omega$ be a quasianalytic convex weight function that is regularly varying of index $p \geq 1$.
	\begin{itemize}
		\item[$(i)$] When $p = 1$, then
			\[ \Phi_{(\omega, \varpi^*)} = \frac{\pi}{2} . \]
			
		\item[$(ii)$] When $p > 1$, and we put $\rho = \min\{p, p^*\}$, then
			\[ \Phi_{(\omega, \varpi^*)} = \frac{\left[\sin\left(\frac{\pi}{2}(\rho - 1)\right)\right]^{1 /\rho}}{(p-1)^{1/p}} .  \]
	\end{itemize}
\end{theorem}

\begin{proof}
$(i)$ By \eqref{eq:BoundsPsiPhi} and Theorem \ref{t:LowBoundPaleyWienerThresholdRegVar} we have that $\Phi_{(\omega, \varpi^*)} \geq \Psi_{(\omega, \varpi)} = \pi/2$.
It thus suffices to show that $\Phi_{(\omega, \varpi^*)} \leq \pi/2$.
Take any $a, b > 0$ such that $ab > \pi/2$.
Suppose $f \in L^2(\R) \setminus \{0\}$ is such that it satisfies \eqref{eq:TFDecayToEntireFunc}, with $\eta = \varpi$, for $a$ and $b$.
By Lemma \ref{l:TFDecayAndEntireFunc}$(i)$, choosing $\theta$ arbitrarily close to $0$, we see that $f$ extends to an entire function satisfying, for any $\varepsilon > 0$,
	\[ |f(z)| \lesssim \exp\left[ \varpi\left(\frac{1+\varepsilon}{b} |y|\right) - \omega\left(\frac{a}{1 + \varepsilon} |x|\right) \right] , \qquad z = x + iy \in \C . \]
Proposition \ref{p:NonTrivEntireLowerBound} then shows that, for $t$ large enough,
	\[ \varpi(t) \leq \varpi\left(\frac{\pi}{2} \frac{(1 + \varepsilon)^3}{ab} t\right) . \]
Therefore, necessarily $ab \leq (\pi / 2) (1 + \varepsilon)^3$.
Letting $\varepsilon \to 0^+$, Lemma \ref{l:PhiExplicit} implies $\Phi_{(\omega, \varpi^*)} \leq \pi/2$.

$(ii)$ By Lemma \ref{l:RegVarEquivOverlineOmega} and \cite[Theorem 5]{S-ClassEntFuncRapDecRealAxis}.
\end{proof}

Finally, for the sake of completeness, we consider the limit case $\omega(t) = t$ explicitly.

\begin{theorem}
\label{t:MainUncertaintyPrincipleCaset}
$(t, e^{t}-1)$ is an optimal uncertainty pair with uncertainty threshold $\pi / 2$.
For any $ab = \pi / 2$, the space of $f \in L^2(\R^d)$ satisfying 
	\begin{equation}
		\label{eq:UncertaintyDecayDenseSubspCaset}
		\sup_{x \in \R^d} |f(x)| e^{\sum_{j=1}^d a |x_j|} < \infty \quad \text{and} \quad \sup_{\xi \in \R^d} |\widehat{f}(\xi)| e^{\sum_{j=1}^d e^{b |\xi_j|}}  < \infty 
	\end{equation}
is dense in $\S(\R^d)$.
\end{theorem}

\begin{proof}
As in the proof of Lemma \ref{l:PhiExplicit}, we may revert to the one-dimensional case with $a = \frac{\pi}{2}$.

Suppose $f \in L^2(\R)$ satisfies \eqref{eq:UncertaintyDecay} for $a = \frac{\pi}{2}$ and $b > 1$.
We can extend $\widehat{f}$ to a bounded analytic function on the strip $\{ z = x + iy \in \C \mid |y| < \frac{\pi}{2} \}$ with continuous extension to its boundary.
Then $g(z) = \widehat{f}(\log z)$ becomes a bounded and analytic function on the right half-plane $\C_+ = \{ z = x + iy \in \C \mid x > 0\}$ with continuous extension to its boundary.
Moreover, by a parametrix argument, we have $|g(x)| = |\widehat{f}(\log x)| \lesssim e^{- x^{b'}}$ for $x > 0$ and any $1 < b' < b$.
It follows from \cite[Theorem 1.4.3, p.~4]{B-EntFunc} that $g = 0$ and thus also $f = 0$.

Let now $f_0 \in L^2(\R)$ be such that $\widehat{f}_0(\xi) = e^{-2 \cosh \xi}$.
Then
	\[ |\widehat{f}_0(\xi)| = e^{- (e^{\xi} + e^{-\xi})} \leq e^{- e^{|\xi|}} , \qquad \xi \in \R . \]
As $\widehat{f}_0$ is even, we have by \cite[(7.94), p.~219]{Y-IndexTrans}, for $|x| \geq 1$,
	\[ 2\pi |f_0(x)| = \left| \int_{-\infty}^\infty \widehat{f}_0(\xi) e^{i x \xi} \, \mathrm{d}\xi \right| = \left| \int_{-\infty}^\infty e^{-2 \cosh \xi} \cos(x \xi) \, \mathrm{d}\xi \right| \lesssim \frac{1}{\sqrt{\sinh(\pi |x|)}} \lesssim e^{- \frac{\pi}{2} |x|} . \]
If we put $f_y(x) = f_0(x - y)$, $y \in \R$, then also
	\[ |f_y(x)| \lesssim e^{-\frac{\pi}{2} |x|} \qquad \text{and} \qquad |\widehat{f}_y(\xi)| \lesssim e^{- e^{|\xi|}} . \]
Therefore, the linear span $\mathcal{V}$ of the set $\{ f_y \mid y \in \R \}$ consists only of functions in $\S(\R)$ \cite{C-C-K-CharSchwartz} satisfying \eqref{eq:UncertaintyDecayDenseSubspCaset} for $a = \frac{\pi}{2}$ and $b = 1$. 
Suppose $T \in \S'(\R)$ is such that $\ev{T}{\widehat{g} \,} = 0$ for all $g \in \mathcal{V}$.
The Fourier transform of $\widehat{f}_0 T$ is given by
\[ \ev{\widehat{f}_0 T}{e^{- i \xi \cdot} \,} = \ev{T}{\widehat{f}_{\xi}} = 0 , \qquad \xi \in \R . \]
So $\widehat{f}_0 T = 0$.
As $\widehat{f}_0(\xi) > 0$ for every $\xi \in \R$, it follows that for any $\varphi \in \D(\R)$, the space of smooth functions with compact support, also $\varphi / \widehat{f}_0 \in \D(\R)$.
But then
	\[ \ev{T}{\varphi} = \left\langle T , \widehat{f}_0 \cdot \frac{\varphi}{\widehat{f}_0} \right\rangle = \left\langle \widehat{f}_0 T , \frac{\varphi}{\widehat{f}_0} \right\rangle = 0 , \qquad \forall \varphi \in \D(\R) . \]
We conclude that $T = 0$.
By the Hahn-Banach theorem and since the Fourier transform is an isomorphism on $\S(\R)$, we finally see that $\mathcal{V}$ is dense in $\S(\R)$.
Finally, tensoring $\mathcal{V}$ also gives density in $\S(\R^d)$.

To conclude our proof, observe that any element in $\mathcal{V}$ satisfies \eqref{eq:UncertaintyDecayDenseSubsp} for $(\omega, \eta) = (t, e^t-1)$ and any $0 < a < \pi / 2$ and $0 < b < 1$.
Therefore, the uncertainty threshold $\Phi_{(t, e^t-1)}$ is exactly $\pi / 2$.
\end{proof}

\section{Weight sequences}
\label{sec:WeightSeq}

In the remaining sections, we use weight sequences to define Gelfand-Shilov spaces. 
We therefore begin by reviewing the relevant properties of weight sequences and establishing counterparts of the results from Section \ref{sec:Weights}. 
We closely follow the classical notation of Komatsu \cite{K-UltradistI}.

\begin{definition}
\label{def:WeightSeq}
A \emph{weight sequence} is a sequence of positive real numbers $M = (M_p)_{p \in \N}$ so that $\lim_{p \to \infty} M_p^{1/p} = \infty$ and which is logarithmically convex:
	\begin{equation} 
		\label{eq:M1}
		M_{p}^2 \leq M_{p-1} M_{p+1} , \qquad p \in \Z_+ . 
	\end{equation}
We call $M$ \emph{normalized} when $M_0 = M_1 = 1$.

We write $m_p = M_p / M_{p-1}$ for $p \in \Z_+$.
Then $(m_p)_{p \in \Z_+}$ is a nondecreasing sequence by \eqref{eq:M1} and $m_1 = 1$ if $M$ is normalized.
\end{definition}

For two sequences $M = (M_p)_{p \in \N}$ and $N = (N_p)_{p \in \N}$ of positive real numbers we write $M \preceq N$ if, for some $h, C > 0$,
	\[ M_p \leq C h^p N_p , \qquad p \in \N . \]
If for any $h > 0$ there is a $C = C_h > 0$ such that the previous inequality holds, we write $M \triangleleft N$.
The two sequences are called \emph{equivalent}, denoted by $M \asymp N$, if $M \preceq N$ and $N \preceq M$.
Any weight sequence is equivalent to a normalized weight sequence.

We will often impose one or more of the following conditions on a weight sequence.

\begin{definition}
Let $M$ be a weight sequence.
We say $M$ satisfies
	\begin{itemize}
		\item[$(M.1)^*$] if the sequence $(p! / M_p)_{p \in \N}$ is a weight sequence;
		\item[$(M.2)$\hphantom{$^*$}] if there exist $C, H > 0$ such that $M_{p + q} \leq C H^{p+q} M_p M_q$ for all $p, q \in \N$;
		\item[$(M.2)^*$] if $2 m_{p} \leq m_{Np}$, $p \geq p_0$, for some $p_0, N \in \Z_+$;
		\item[$(\QA)$\hphantom{$^*$}] if $\sum_{p = 1}^\infty \frac{M_{p-1}}{M_p} = \infty$.
	\end{itemize}
The conditions $(M.2)$, $(M.2)^*$, and $(\QA)$ are invariant under equivalence.
\end{definition}

For any weight sequence $M$, we write $M^*$ for the conjugate sequence
	\[ M^*_p = \frac{p!}{M_p} , \qquad p \in \N . \]
Then $M^*$ is a weight sequence if $M$ satisfies $(M.1)^*$, and in this case also $M^*$ satisfies $(M.1)^*$.
If a weight sequence $M$ satisfies $(M.1)^*$, then necessarily $M \triangleleft p!$, so that $M$ satisfies $(QA)$ \cite[Lemma 4.1, p.~55]{K-UltradistI}.
Additionally, any weight sequence satisfying $(M.1)^*$ satisfies $(M.2)$ \cite[Corollary 4.7]{S-ConjWeight}.

Let $M$ be a weight sequence.
We define its \emph{associated weight function} as
	\begin{equation*}
		\label{eq:AssocFunc} 
		\omega_M(t) = \sup_{p \in \N} \log \frac{t^p M_0}{M_p} .
	\end{equation*}
Then $\omega_M$ is a weight function as in Definition \ref{def:WeightFunc}.
By \eqref{eq:M1} we have \cite[Proposition 3.2, p.~49]{K-UltradistI}:
	\begin{equation}
		\label{eq:SeqFromAssocFunc} 
		M_p = M_0 \sup_{t \geq 0} \frac{t^p}{e^{\omega_M(t)}} , \qquad p \in \N .
	\end{equation}
Therefore, two weight sequences $M$ and $N$ are equivalent if and only if $\omega_M \asymp \omega_N$.
The above conditions on a weight sequence $M$ can be detected through $\omega_M$.

\begin{lemma}
\label{l:WeighSeqCondViaWeighFunc}
Let $M$ be a weight sequence.
	\begin{itemize}
		\item[$(i)$] $M$ satisfies $(M.2)$ if and only if $\omega_M$ is lower moderate.
		
		\item[$(ii)$] $M$ satisfies $(M.2)^*$ if and only if $\omega_M$ is upper moderate.
		
		\item[$(iii)$] $M$ satisfies $(\QA)$ if and only if $\omega_M$ is quasianalytic.
	\end{itemize}
\end{lemma}

\begin{proof}
Statement $(i)$ is shown in \cite[Proposition 3.6, p.~51]{K-UltradistI}, $(ii)$ in \cite[Proposition 13 and its proof]{B-M-M-CompareUltradiff}, and $(iii)$ follows from \cite[Lemma 4.1, p.~55]{K-UltradistI}.
\end{proof}

For a weight sequence $M$, we may consider the weight function $\varpi_M := \varpi_{\omega_M}$.
We now construct an explicit weight sequence $\overline{M}$ such that $\varpi_M \asymp \omega_{\overline{M}}$.

Consider the counting function of $\{ m_p \mid p \in \Z_+ \}$:
	\[ m(t) = \sum_{m_p \leq t} 1 , \qquad t \geq 0 . \]
Then, $\omega_M(t) = 0$ for $t \in [0, m_1]$ and
	\begin{equation}
		\label{eq:AssocFuncAlternative} 
		\omega_M(t) = \int_{m_1}^t \frac{m(u)}{u} \, \mathrm{d}u = \sum_{m_p \leq t} \log \frac{t}{m_p} = \sum_{p = 1}^\infty \log_+ \frac{t}{m_p} , \qquad t \geq m_1 . 
	\end{equation}
This immediately gives, for any $\lambda > 1$,
	\begin{equation}
		\label{eq:UpperBoundm} 
		m(t) \leq \frac{1}{\log \lambda} \omega_M(\lambda t) , \qquad t \geq 0 . 
	\end{equation}
Moreover, when $M$ is normalized, we find
	\begin{equation}
		\label{eq:varpiMAsSum}
		\varpi_M(t) 
		= t \sum_{m_p \leq t} \int_{m_p}^t \frac{\log(s / m_p)}{s^2} \, \mathrm{d}s
		= \sum_{m_p \leq t} \left( \frac{t}{m_p} - 1 - \log(t / m_p) \right) , \quad t \geq 1 . 
	\end{equation}
As $\log x \leq x/2  - 1$ for $x \geq 6$, we get the estimate
	\begin{equation}
		\label{eq:varpiMUpperLowerBound}
		\frac12 \sum_{m_p \leq t/6} \frac{t}{m_p} \leq \varpi_M(t) \leq \sum_{m_p \leq t} \frac{t}{m_p} . 
	\end{equation}
These bounds inspire the following definition.
	
\begin{definition}
\label{def:overlineM}
For a weight sequence $M$ we let $\overline{m}_p$ be the nondecreasing rearrangement, with multiplicities, of
	\[ \{ k \cdot m_p \mid k, p \in \Z_+ \} . \]
We define the sequence $\overline{M} = \overline{M}_M = (\overline{M}_p)_{p \in \N}$ as
	\[ \overline{M}_0 = 1 , \qquad \overline{M}_p = \prod_{j = 1}^p \overline{m}_j , \quad p \in \Z_+ . \]
\end{definition}
	
When $M$ is normalized, then clearly $\overline{M}_p \leq M_p$ and $\overline{M}_p \leq p!$ for all $p \in \N$.
We denote by $\overline{m}(t)$ the counting function of the sequence $(\overline{m}_p)_{p \in \Z_+}$, that is,
	\begin{equation}
		\label{eq:overlinem} 
		\overline{m}(t) = \sum_{\overline{m}_p \leq t} 1 = \sum_{p = 1}^\infty \left\lfloor \frac{t}{m_p} \right\rfloor . 
	\end{equation}
Comparing \eqref{eq:varpiMUpperLowerBound} and \eqref{eq:overlinem}, we find the following result. 

\begin{lemma}
\label{l:overlineM}
Let $M$ be a weight sequence.
Then $\overline{M}$ is also a weight sequence satisfying $(M.2)$, $(\QA)$, and $\varpi_M \asymp \omega_{\overline{M}}$.
If $M$ satisfies $(M.2)^*$ then so does $\overline{M}$.
\end{lemma}

\begin{proof}
For any $h > 0$, eventually $ m_p \geq h$, so the same is true for $\overline{m}_p$.
This shows $\lim_{p \to \infty} \overline{M}_p^{1/p} = \infty$.
Also, as $\overline{m}_p$ is nondecreasing, one has \eqref{eq:M1}.
Thus $\overline{M}$ is a weight sequence.

Let us now show $\varpi_M \asymp \omega_{\overline{M}}$.
We may assume $M$ is normalized.
By \eqref{eq:UpperBoundm}, \eqref{eq:varpiMUpperLowerBound}, \eqref{eq:overlinem}, and the convexity of $\varpi_M$, we find the bound
	\[ \varpi_M\left(\frac{t}{2}\right) \leq \frac12 \sum_{m_p \leq t} \frac{t}{m_p}  \leq \sum_{m_p \leq t} \left\lfloor \frac{t}{m_p} \right\rfloor = \overline{m}(t) \leq \omega_{\overline{M}}(e t) . \]
On the other hand, combining Lemma \ref{l:PropertiesVarpi}$(iv)$, \eqref{eq:AssocFuncAlternative}, \eqref{eq:UpperBoundm}, \eqref{eq:overlinem}, and the convexity of $\varpi_M$, we find through integration by parts that
	\begin{align*} 
		\omega_{\overline{M}}(t) 
		&= \int_{1}^t \frac{\overline{m}(u)}{u} \, \mathrm{d}u 
		\leq \int_{1}^t \frac{1}{u} \left( \int_{1}^u \frac{u}{s} \, \mathrm{d}m(s) \right) \, \mathrm{d}u 
		\leq t \int_{1}^t \frac{\mathrm{d}m(s)}{s}  \\
		&\leq m(t) + t \int_{1}^t \frac{m(s)}{s} \frac{\mathrm{d}s}{s}
		\leq m(t) + \omega_M(t) + \varpi_M(t) \\
		&\leq \varpi_M((2e + 3)t) . 
	\end{align*}
	
The final statement now follows from Lemmas \ref{l:PropertiesVarpi} and \ref{l:WeighSeqCondViaWeighFunc}.
\end{proof}

The weight sequence $\overline{M}$ does not necessarily need to satisfy $(M.1)^*$, even if $M$ satisfies $(M.1)^*$.
However, its conjugate $\overline{M}^*$ satisfies almost all properties of a weight sequence.

\begin{lemma}
\label{l:ConjOverlineMAlmostWeightSeq}
Let $M$ be a weight sequence satisfying $(\QA)$.
Then $\lim_{p \to \infty} (\overline{M}^*_p)^{1/p} = \infty$, we have
	\begin{equation}
		\label{eq:overlineMSuperMult}
		\overline{M}^*_p \overline{M}^*_q \leq \overline{M}^*_{p+q} , \qquad p, q \in \N , 
	\end{equation}
and there exists some $H \geq 1$ such that
	\begin{equation}
		\label{eq:overlineM*InTermsOfWeightFunc}
		H^{-(p + 1)} \overline{M}^*_p \leq \sup_{t \geq 0} t^p e^{-\varpi_M^*(t)} \leq H^{p+1} \overline{M}^*_p , \qquad p \in \N . 
	\end{equation}
\end{lemma}

\begin{proof}
For any $N \in \N$ we have, by \eqref{eq:overlinem},
	\[ \frac{\overline{m}(t)}{t} \geq \sum_{p = 1}^N \frac{1}{m_p}  - \frac{N}{t} . \]
As $M$ satisfies $(\QA)$, we see that $\lim_{t \to \infty} \overline{m}(t) / t = \infty$.
Plugging in $t = \overline{m}_p / 2$, we see that $p / \overline{m}_p \to \infty$ as $p \to \infty$.
Thus
	\[ \left(\overline{M}_p^*\right)^{\frac{1}{p}} = \left( \prod_{q=1}^p \frac{q}{\overline{m}_q} \right)^{\frac{1}{p}} \to \infty , \qquad p \to \infty . \]

We now show \eqref{eq:overlineMSuperMult}.
By \eqref{eq:overlinem} we have, for any $t,s \geq 0$,
	\[ \overline{m}(t + s) = \sum_{p = 1}^\infty \left\lfloor \frac{t + s}{m_p} \right\rfloor \geq \sum_{p = 1}^\infty \left\lfloor \frac{t}{m_p} \right\rfloor + \sum_{p = 1}^\infty \left\lfloor \frac{s}{m_p} \right\rfloor = \overline{m}(t) + \overline{m}(s) .  \]
Plugging in $t = \overline{m}_p$ and $s = \overline{m}_q$, we have $\overline{m}(\overline{m}_p + \overline{m}_q) \geq p + q$, so 
	\begin{equation}
		\label{eq:overlinemSubAdd} 
		\overline{m}_{p+q} \leq \overline{m}_p + \overline{m}_q , \qquad p,q \geq 1 . 
	\end{equation}
We now prove \eqref{eq:overlineMSuperMult} by induction on $p + q$.
When $p = 0$ or $q = 0$, it is immediate.
Suppose now $p, q \geq 1$ and \eqref{eq:overlineMSuperMult} holds for $p$ and $q - 1$ and for $p - 1$ and $q$.
Then
	\begin{align*} 
		\overline{M}^*_{p+q} 
		&= \frac{p + q}{\overline{m}_{p + q}} \overline{M}^*_{p+q-1} 
		\geq  \frac{p + q}{\overline{m}_{p + q}} \max\left\{ \overline{M}^*_{p} \overline{M}^*_{q-1} , \overline{M}^*_{p-1} \overline{M}^*_q \right\}  \\
		&= \overline{M}^*_p \overline{M}^*_q \, \frac{p + q}{\overline{m}_{p + q}} \max\left\{ \frac{\overline{m}_q}{q} , \frac{\overline{m}_p}{p} \right\}
		\geq  \overline{M}^*_p \overline{M}^*_q \, \frac{\overline{m}_p + \overline{m}_q}{\overline{m}_{p+q}} \\
		&\geq \overline{M}^*_p \overline{M}^*_q . 
	\end{align*}
	
Finally, we verify \eqref{eq:overlineM*InTermsOfWeightFunc}.
By Lemma \ref{l:overlineM} there exist $C_j, H_j \geq 1$, $j = 0, 1$, such that $\varpi_M(t) \leq \omega_{\overline{M}}(H_0 t) + \log C_0$ and $\omega_{\overline{M}}(t) \leq \varpi_M(H_1 t) + \log C_1$.
We have by \eqref{eq:SeqFromAssocFunc}, for any $t \geq 0$,
	\[ e^{\varpi_M^*(t)} = \sup_{s \geq 0} e^{s t - \varpi_M(s)} \geq \frac{t^p}{p!} \sup_{s \geq 0} s^p e^{-\varpi_M(s)} \geq \frac{1}{C_0} \frac{t^p}{p!} \sup_{s \geq 0} s^p e^{- \omega_{\overline{M}}(H_0 s)} = \frac{(t H_0^{-1})^p}{C_0} \frac{\overline{M}_p}{p!} . \]
Hence
	\[ \sup_{t \geq 0} t^{p} e^{-\varpi_M^*(t)} \leq C_0 H_0^p \overline{M}^*_p . \]
On the other hand, for any $t \geq 0$ we find, by \eqref{eq:SeqFromAssocFunc}, 
	\[
		e^{\varpi_M^*(t)} = \sup_{s \geq 0} e^{st - \varpi_M(s)} \leq C_1 \sum_{k \in \N} \frac{t^k}{k!} \sup_{s \geq 0} s^k e^{-\omega_{\overline{M}}(s / H_1)} \leq C_1 \sum_{k \in \N} \frac{(H_1 t)^k}{\overline{M}^*_k} . 
	\]
Fix some $p \in \Z_+$, then we see that, for any $t \geq 0$,
	\begin{equation}
		\label{eq:overlineM*expvarpi_M*} 
		\overline{M}^*_p (H_1 t)^{-p} e^{\varpi_M^*(t)} \leq C_1 \sum_{k = 0}^{p - 1} (H_1 t)^{k-p} \frac{\overline{M}^*_p}{\overline{M}^*_k} + C_1 \sum_{k = p}^{\infty} (H_1 t)^{k-p} \frac{\overline{M}^*_p}{\overline{M}^*_k} .
	\end{equation}
Put
	\[ c = \inf_{j \geq p} \frac{j}{\overline{m}_j} . \]
For any $1 \leq k \leq p$ and $j \geq p$, write $j = q k + r$ for some $0 \leq r < k$.
Iterating \eqref{eq:overlinemSubAdd}, we find that $\overline{m}_j \leq (q + 1) \overline{m}_k$.
Therefore
	\[ \frac{k}{\overline{m}_k} \leq 2 \frac{j}{\overline{m}_j} .  \]
This shows that $k / \overline{m}_k \leq 2 c$ for any $1 \leq k \leq p$.
Plugging in $t = c / (2 H_1)$ in \eqref{eq:overlineM*expvarpi_M*} then gives
	\begin{align*} 
		\overline{M}^*_p \left(\frac{c}{2}\right)^{-p} e^{\varpi_M^*(c/2H_1)} 
		&\leq C_1 \sum_{k = 0}^{p-1} \left(\frac{c}{2}\right)^{k-p} \frac{\overline{M}^*_p}{\overline{M}^*_k} + C_1 \sum_{k = p}^{\infty} \left(\frac{c}{2}\right)^{k-p} \frac{\overline{M}^*_p}{\overline{M}^*_k} \\
		&\leq C_1 \sum_{k = 0}^{p-1} 4^{p - k} + C_1 \sum_{k = p}^\infty 2^{-(k -p)} \leq 4 C_1 4^{p} .
	\end{align*}
Therefore
	\[ 4^{-p} \overline{M}^*_p \leq 4 C_1 \left(\frac{c}{2}\right)^p e^{-\varpi_M^*(c/2H_1)} \leq 4 C_1 H_1^p \sup_{t \geq 0} t^p e^{-\varpi_M^*(t)} . \qedhere \]
\end{proof}

We now characterize when $M \asymp \overline{M}$ or $\overline{M} \triangleleft M$.
To do this, we introduce analogous notions as in Definitions \ref{def:Matuszewska} and \ref{def:GeomPers} for weight sequences.
Here $GM(m_1, \ldots, m_p) = (m_1 \cdots m_p)^{1/p}$ denotes the geometric mean and $HM(m_1, \ldots, m_p) = p / \sum_{j=1}^p m_j^{-1}$ denotes the harmonic mean of $\{m_1, \ldots, m_p\}$.

\begin{definition}
Let $M$ be a weight sequence.
	\begin{itemize}
		\item[$(i)$] The \emph{upper Matuszewska index $\alpha(M)$} of $M$ is defined as
			\[
				\alpha(M) = \inf \left\{ \alpha \in \R \mid \lim_{k \to \infty} \limsup_{p \to \infty} \frac{m_{kp}}{k^\alpha m_p} = 0 \right\} .
			\]
		
		\item[$(ii)$] The weight sequence $M$ is called \emph{geometrically persistent} if
			\begin{equation*}
				\label{eq:GeomPersSeq}
				\lim_{p \to \infty} \frac{GM(m_1, \ldots, m_p)}{HM(m_1, \ldots, m_p)} = \infty . 
			\end{equation*}
	\end{itemize}
\end{definition}

We have $\alpha(M) < \infty$ if and only if $M$ satisfies $(M.2)$ \cite[Corollary 3.17]{JG-S-S-ORegVar}.
If $M$ satisfies $(M.1)^*$, then always $\alpha(M) \leq 1$.
The sequence $M = (p!)_{p \in \N}$ has $\alpha(M) = 1$ and is geometrically persistent.

\begin{lemma}
\label{l:CondFromSeqToFunc}
Let $M$ be a weight sequence.
	\begin{itemize}
		\item[$(i)$] $\alpha(M) = 1 / \beta(\omega_M)$.
		\item[$(ii)$] $M$ is geometrically persistent if and only if $\omega_M$ is geometrically persistent.
	\end{itemize}
\end{lemma}

\begin{proof}
$(i)$ This is shown in \cite[Theorem 3.16 and Corollary 4.6]{JG-S-S-ORegVar}.

$(ii)$ We may assume $M$ is normalized.
By Proposition \ref{p:CharCondDomVarpi}, it suffices to prove that $M$ is geometrically persistent if and only if $\omega_M(t) = o(\varpi_M(t))$.
Fix some $t > m_1$.
Then $m_p \leq t < m_{p+1}$ for some $p \in \Z_+$.
By \eqref{eq:varpiMAsSum}, we have
	\begin{equation}
		\label{eq:ExpansionVarpiM}
		\varpi_M(t)
		= \sum_{q = 1}^p \left( \frac{t}{m_q} - 1 - \log \frac{t}{m_q} \right)
		= \frac{p t}{HM(m_1, \ldots, m_p)} - p - \omega_M(t) .
	\end{equation}
Now
	\[ \omega_M(t) = \log \frac{t^p}{M_p} = p \log \frac{t}{GM(m_1, \ldots, m_p)} = p \cdot \tau , \]
where we have put $\tau = \log (t / GM(m_1, \ldots, m_p))$.
Then
	\[ \frac{\varpi_M(t)}{\omega_M(t)} = \frac{GM(m_1, \ldots, m_p)}{HM(m_1, \ldots, m_p)} \frac{e^{\tau}}{\tau} - \frac{1}{\tau} - 1 . \]
Since $e^\tau \geq \tau + 1$, we get
	\[ \frac{\varpi_M(t)}{\omega_M(t)} \geq \left( \frac{GM(m_1, \ldots, m_p)}{HM(m_1, \ldots, m_p)} - 1 \right) \left( 1 + \frac{1}{\tau} \right) \geq \frac{GM(m_1, \ldots, m_p)}{HM(m_1, \ldots, m_p)} - 1  . \]
Thus, if $M$ is geometrically persistent, one has $\omega_M(t) = o(\varpi_M(t))$.

Suppose now that $\omega_M(t) = o(\varpi_M(t))$.
Choose the sequence
	\[ t_p = e GM(m_1, \ldots, m_p) = e M_p^{1/p} , \qquad p \in \Z_+ . \]
Then
	\[ \omega_M(t_p) = \sup_{q \in \N} \log \frac{t^q_p}{M_q} \geq \log \frac{e^p M_p}{M_p} = p , \qquad p \in \Z_+ .  \]
Fix some $p \in \Z_+$ and let $q \in \Z_+$ be such that $m_q \leq t_p < m_{q+1}$.
By \eqref{eq:ExpansionVarpiM}, we have
	\begin{align*} 
		\varpi_M(t_p) 
		&\leq \sum_{j = 1}^{\min\{p, q\}} \frac{t_p}{m_j} + \sum_{p < j \leq q} \left[ \frac{t_p}{m_j} - 1 - \log \frac{t_p}{m_j} \right] \\
		&\leq e \frac{GM(m_1, \ldots, m_p)}{HM(m_1, \ldots, m_p)} p + \sum_{p < j \leq q} \left[ \frac{t_p}{m_j} - 1 - \log \frac{t_p}{m_j} \right] 
		. 
	\end{align*}
We estimate the second term.
We assume $p < q$, as otherwise it is empty. 
For $j > p$, we have $m_j \geq m_p \geq GM(m_1, \ldots, m_p)$.
Therefore
	\[ 1 \leq \frac{t_p}{m_j} \leq e , \qquad p < j \leq q . \]
Since $x - 1 - \log x \leq e \log x$ for $1 \leq x \leq e$, it follows that
	\[ \sum_{p < j \leq q} \left[ \frac{t_p}{m_j} - 1 - \log \frac{t_p}{m_j} \right] \leq e \sum_{p < j \leq q} \log \frac{t_p}{m_j} \leq e \omega_M(t_p) .  \]
We conclude
	\[ \frac{\varpi_M(t_p)}{\omega_M(t_p)} \leq e \left( \frac{GM(m_1, \ldots, m_p)}{HM(m_1, \ldots, m_p)}  + 1 \right) , \qquad p \in \Z_+ . \]
Since $\varpi_M(t_p) / \omega_M(t_p) \to \infty$ as $p \to \infty$, we see that $M$ is geometrically persistent.
\end{proof}

We may now characterize, in terms of $M$, when $M$ and $\overline{M}$ are equivalent, or when $\overline{M} \triangleleft M$.

\begin{proposition}
\label{p:CondEquivSeq}
Let $M$ be a weight sequence.
	\begin{itemize}
		\item[$(i)$] $M \asymp \overline{M}$ if and only if $\alpha(M) < 1$.
		
		\item[$(ii)$] $\overline{M} \triangleleft M$ if and only if $M$ is geometrically persistent.
	\end{itemize}
\end{proposition}

\begin{proof}
Follows from a combination of Propositions \ref{p:SuffNeccCondEquivVarpi} \and \ref{p:CharCondDomVarpi}, \eqref{eq:SeqFromAssocFunc}, and Lemmas \ref{l:overlineM} and \ref{l:CondFromSeqToFunc}. 
\end{proof}

\section{An optimal result for the decay and regularity of smooth functions}
\label{sec:OptimalDecayRegularity}

The goal of this section is to show an optimality result for the decay and regularity as in \eqref{eq:GelfandShilov}.
Concretely, we will prove the following result.

\begin{theorem}
\label{t:UncertaintyGS}
Let $M$ be a weight sequence satisfying $(M.2)$, $(M.2)^*$, and $(\QA)$.
There exists a constant $\mathscr{R}_M \in (0, \infty)$ such that for any $h, k > 0$:
	\begin{itemize}
		\item[(I)] If $h k < \mathscr{R}_M$, should $f \in C^\infty(\R^d)$ satisfy
			\begin{equation} 
				\label{eq:GSCond}
				\sup_{x \in \R^d} | x^\alpha f^{(\beta)}(x) | \lesssim k^\alpha \overline{M}^*_{|\alpha|} \cdot h^\beta M_{|\beta|} , \qquad \alpha, \beta \in \N^d ,
			\end{equation}
		then $f = 0$.
		
		\item[(II)] If $h k > \mathscr{R}_M$, then the set of $f \in C^\infty(\R^d)$ satisfying \eqref{eq:GSCond} is dense in $\S(\R^d)$.
	\end{itemize}
\end{theorem}

As a consequence of Proposition \ref{p:CondEquivSeq}, we may then also establish sufficient conditions under which $\overline{M}$ may or may not be replaced by $M$ in \eqref{eq:GSCond}.

\begin{theorem}
\label{t:UncertaintyGSOriginalM}
Let $M$ be a weight sequence satisfying $(M.2)$, $(M.2)^*$, and $(\QA)$.
	\begin{itemize}
		\item[1.] If $\alpha(M) < 1$, when $k,h > 0$ are such that $kh$ is small enough, then any $f \in C^\infty(\R^d)$ satisfying
			\begin{equation}
				\label{eq:GSCondOriginalM} 
				\sup_{x \in \R^d} | x^\alpha f^{(\beta)}(x) | \lesssim k^\alpha M^*_{|\alpha|} \cdot h^\beta M_{|\beta|} , \qquad \alpha, \beta \in \N^d , 
			\end{equation}
		must necessarily be trivial.
		On the other hand, if $hk$ is large enough, there exists a dense subspace of $\S(\R^d)$ of functions satisfying \eqref{eq:GSCondOriginalM}.

		\item[2.] If $M$ is geometrically persistent, there exists no nontrivial $f \in C^\infty(\R^d)$ satisfying \eqref{eq:GSCondOriginalM} for some $h, k > 0$. 
	\end{itemize}
\end{theorem}

The idea for proving Theorem \ref{t:UncertaintyGS} is to translate it into the terms of Theorem \ref{t:CharNonTrivEntireWeightFunc}.
We first show two intermediate results that allow this.
We start by observing that the proof boils down to showing the existence and nonexistence of a nontrivial smooth function satisfying \eqref{eq:GSCond}.

\begin{lemma}
\label{l:RExplicit}
Let $M$ be a weight sequence satisfying $(M.2)$, $(M.2)^*$, and $(\QA)$.
Then statements (I) and (II) of Theorem \ref{t:UncertaintyGS} hold for
	\begin{equation}
		\label{eq:RExplicit} 
		\mathscr{R}_M = \inf\{ h k \mid \text{there exists $f \in C^\infty(\R) \setminus \{0\}$ satisfying \eqref{eq:GSCond} for $h$ and $k$} \} . 
	\end{equation}
\end{lemma}

\begin{proof}
If $f \in \S(\R)$ satisfies \eqref{eq:GSCond} for certain $h_0, k_0 > 0$, and if $h_1, k_1 > 0$ are such that $h_0 k_0 = h_1 k_1$, then $g(x) = f\left(\frac{h_1}{h_0} x\right)$ satisfies \eqref{eq:GSCond} for $h_1$ and $k_1$.
Therefore, the infimum in \eqref{eq:RExplicit} is well-defined.

Suppose $f \in C^\infty(\R^d) \setminus \{0\}$ satisfies \eqref{eq:GSCond} for certain $h, k > 0$.
Then, there is some $x' \in \R^{d-1}$ such that the smooth function $f_{x'}(x) = f(x, x')$ is nontrivial.
Clearly $f_{x'}$ also satisfies \eqref{eq:GSCond} for $h$ and $k$.
Therefore, it cannot be that $hk$ is smaller than the value on the right-hand side of \eqref{eq:RExplicit}.
Thus, this value can be at most $\mathscr{R}_M$.

Suppose now that $0 < h_0 < h$ and $0 < k_0 < k$ are such that there exists some $g_0 \in C^\infty(\R) \setminus \{0\}$ satisfying \eqref{eq:GSCond} for $h_0$ and $k_0$.
Putting $f_0 = g_0 \otimes \cdots \otimes g_0$, where we take $d$ tensor products, then $f_0 \in \S(\R^d) \setminus \{0\}$ and $f_0$ also satisfies \eqref{eq:GSCond} in dimension $d$ for $h_0$ and $k_0$.
For any $\lambda = (t, \xi) \in \R^{2d}$, we consider the time-frequency shift
	\[ [\pi(\lambda) f_0](x) := f_0(x - t) e^{i \ev{x}{\xi}} , \qquad x \in \R^d . \]
For any $\alpha, \beta \in \N^d$, by Lemma \ref{l:ConjOverlineMAlmostWeightSeq},
	\begin{align*}
		| x^\alpha [\pi(\lambda) f_0]^{(\beta)}(x)|
		&\leq \sum_{\alpha' \leq \alpha} \sum_{\beta' \leq \beta} {\alpha \choose \alpha'} {\beta \choose \beta'} |t|^{\alpha'} |\xi|^{\beta'} \cdot |(x - t)^{\alpha - \alpha'} f_0^{(\beta - \beta')}(x - t)| \\
		&\lesssim \sum_{\alpha' \leq \alpha} \sum_{\beta' \leq \beta} {\alpha \choose \alpha'} {\beta \choose \beta'} \left(k - k_0\right)^{|\alpha'|} \overline{M}^*_{|\alpha'|} \cdot \left(h - h_0\right)^{|\beta'|} M_{|\beta'|} \\
		&\qquad \qquad \qquad \qquad \cdot k_0^{|\alpha - \alpha'|} \overline{M}^*_{|\alpha - \alpha'|} \cdot h_0^{|\beta - \beta'|}  M_{|\beta - \beta'|} \\
		&\lesssim k^{|\alpha|} \overline{M}^*_{|\alpha|} \cdot h^{|\beta|} M_{|\beta|} .
	\end{align*}
We conclude that the elements in the linear span of $\{ \pi(\lambda) f_0 \mid \lambda \in \R^{2d} \}$ all satisfy \eqref{eq:GSCond} for $h$ and $k$.
This is a dense subspace of $\S(\R^d)$, as otherwise the short-time Fourier transform with window $f_0$ 
	\[ V_{f_0} f(\lambda) = \ev{f}{\pi(\lambda) f_0} , \qquad f \in \S^\prime(\R^d) , \]
would not be injective, contradicting \cite[Corollary 11.2.7, p.~229]{G-FoundTFAnal}. 
We conclude that the right-hand side in \eqref{eq:RExplicit} is at least $\mathscr{R}_M$.
\end{proof}

Next, we observe that when the decay and regularity of a smooth function grow sufficiently small, its Fourier transform extends to an entire function satisfying Paley--Wiener-type bounds as in Section \ref{sec:EntireFunc}, and vice versa.
We will need this result in the following form.

\begin{lemma}
\label{l:GSToEntire}
Let $M$ be a weight sequence satisfying $(M.2)$ and let $\eta$ be a convex weight function such that $t = o(\eta(t))$.
\begin{itemize}
	\item[$(i)$] If $f \in C^\infty(\R)$ satisfies, for certain $h, k > 0$,
		\[ \sup_{x \in \R} |e^{\eta^*\left(\frac{|x|}{k}\right)} f^{(n)}(x)| \lesssim h^n M_n , \qquad n \in \N , \]
	then $\widehat{f}$ extends to an entire function satisfying, for any $k' > k$,
		\[ |\widehat{f}(z)| \lesssim \exp\left[ \eta\left(k' |y|\right) - \omega_M\left(\frac{|z|}{h}\right) \right] , \qquad z = x + iy \in \C . \]
		
	\item[$(ii)$] Suppose $\omega_M \preceq \eta$.
	Then there exists some $H > 1$ such that if $f \in L^2(\R)$ extends to an entire function satisfying, for certain $a, b > 0$,
		\[ |f(z)| \lesssim \exp\left[ \eta\left(\frac{|y|}{b}\right) - \omega_M(a |x|) \right] , \qquad z = x + iy \in \C ,  \]
	then $\widehat{f} \in C^\infty(\R)$ and
		\[ \sup_{x \in \R} |e^{\eta^*\left(\left[Ha + \frac{1}{b}\right]^{-1} |x|\right)} \widehat{f}^{(n)}(x)| \lesssim \left(\frac{H}{a}\right)^{n} M_n , \qquad n \in \N . \]
\end{itemize}
\end{lemma}

\begin{proof}
We fix $0 < \varepsilon < 1$ throughout the proof.

$(i)$ 
By the convexity of $\eta^*$ we have
	\[ \int_{-\infty}^\infty |f^{(n)}(x)| e^{\eta^*\left(\frac{1 - \varepsilon}{k} |x|\right)} \, \mathrm{d}x \leq \left[ \sup_{x \in \R} |f^{(n)}(x)| e^{\eta^*\left(\frac{|x|}{k}\right)} \right] \int_{-\infty}^\infty e^{-\varepsilon \eta^*\left(\frac{|x|}{k}\right)} \, \mathrm{d}x \lesssim h^n M_n . \]
As $t = o(\eta^*(t))$, we may extend $\widehat{f}$ to an entire function as
	\[ \widehat{f}(z) = \int_{-\infty}^\infty f(t) e^{-i z t} \, \mathrm{d}t . \]
Then, for $z = x + iy \in \C \setminus \{0\}$,
		\begin{align*}
		|\widehat{f}(z)|
		&= \inf_{n \in \N} |z|^{-n} \left| \int_{-\infty}^\infty f(t) (-iz)^n e^{-izt} \, \mathrm{d}t \right|
		\leq \inf_{n \in \N} |z|^{-n} \int_{-\infty}^\infty |f^{(n)}(t)| e^{|y t|} \, \mathrm{d}t \\
		&\lesssim \left( \inf_{n \in \N} \frac{h^n M_n}{|z|^n} \right) \sup_{t \in \R} \exp\left[ |y| |t| - \eta^*\left(\frac{1-\varepsilon}{k} |t|\right) \right] \\
		&\lesssim \exp\left[ \eta\left( \frac{k}{1 - \varepsilon} |y| \right) - \omega_M\left(\frac{|z|}{h}\right) \right] . 
	\end{align*}

$(ii)$
Let $C_0, H_0 > 0$ be such that $\omega_M(t) \leq \eta(H_0 t) + C_0$.
Also, as $M$ satisfies $(M.2)$, let $C_1, H_1 \geq 1$ be such that $2 \omega_M(t) \leq \omega_M(H_1 t) + C_1$.
The function $e^{-\omega_M(|\cdot| / H_1)}$ is integrable.
For any $x, y \in \R$ such that $xy \leq 0$, by changing the contour and using \eqref{eq:SeqFromAssocFunc} and the convexity of $\eta$, we find, for any $n \in \N$,
	\begin{align*} 
		&|\widehat{f}^{(n)}(x)| 
		=  \left| \int_{-\infty}^\infty (- i \xi)^n f(\xi) e^{- i \xi x} \, \mathrm{d}\xi \right| 
		=  \left| \int_{-\infty}^\infty (y - i \xi)^n f(\xi + iy) e^{- i \xi x + y x} \, \mathrm{d}\xi \right|  \\
		&\quad \leq \sum_{m \leq n} {n \choose m} |y|^{n-m} \int_{-\infty}^\infty |\xi|^m |f(\xi + iy)| e^{y x} \, \mathrm{d}\xi \\
		&\quad \lesssim \exp\left[ \eta\left(\frac{1}{b} |y|\right) - |y x|\right] \sum_{m \leq n} {n \choose m} \left( \frac{(a |y|)^{n-m}}{M_{n-m}} a^{-(n-m)} M_{n-m} \right) \cdot \left( \sup_{\xi \in \R} \frac{|\xi|^{m}}{e^{\omega_M(a |\xi| / H_1)}} \right) \\
		&\quad \lesssim \exp\left[ \eta\left(\frac{1}{b} |y|\right) + \omega_M(a |y|) - |y x|\right] \cdot \left(\frac{1 + H_1}{a}\right)^{n} M_n \\
		&\quad \lesssim \exp\left[ \eta\left(\left[H_0 a + \frac{1}{b}\right]  |y|\right) - |yx| \right] \cdot \left(\frac{1 + H_1}{a}\right)^{n} M_n .
	\end{align*}
Taking the infimum over $y \in \R$, we find
	\[ |\widehat{f}^{(n)}(x)| \lesssim \exp\left[- \eta^* \left(\left[H_0 a + \frac{1}{b}\right]^{-1} |x|\right) \right] \cdot \left(\frac{1 + H_1}{a}\right)^{n} M_n , \qquad n \in \N . \qedhere \]
\end{proof}

We are now sufficiently prepared to show Theorems \ref{t:UncertaintyGS} and \ref{t:UncertaintyGSOriginalM}.

\begin{proof}[Proof of Theorem \ref{t:UncertaintyGS}]
We have that $\omega_M$ is moderate and quasianalytic by Lemma \ref{l:WeighSeqCondViaWeighFunc}.
We may therefore let $\Psi_{(\omega_M, \varpi_M)} \in (0, \infty)$ be as in Theorem \ref{t:CharNonTrivEntireWeightFunc}.

First, take any $a, b > 0$ with $ab < \Psi_{(\omega_M, \varpi_M)}$.
There exists some nontrivial entire function $F$ satisfying
	\[ |F(z)| \lesssim \exp\left[ \varpi_M\left(\frac{|y|}{b}\right) - \omega_M(a |x|) \right] , \qquad z \in \C . \]
Putting $f = \widehat{F}$, it follows from Lemma \ref{l:GSToEntire}$(ii)$ with $\eta = \varpi_M$ that for some $H_0 \geq 1$,
	\[ \left| e^{\varpi_M^*\left(\left[H_0 a + \frac{1}{b}\right]^{-1} |x|\right)} f^{(n)}(x) \right| \lesssim (a / H_0)^{-n} M_n , \qquad n \in \N . \]
Suppose $H_1 \geq 1$ is such that \eqref{eq:overlineM*InTermsOfWeightFunc} holds.
Then, for any $m, n \in \N$,
	\begin{align*} 
		|x^m f^{(n)}(x)| 
		&\leq |e^{\varpi_M^*\left(\left[H_0 a + \frac{1}{b}\right]^{-1} |x|\right)} f^{(n)}(x)| \left[\sup_{t \geq 0} t^m e^{-\varpi_M^*\left(\left[H_0 a + \frac{1}{b}\right]^{-1} t\right)}\right] \\
		&\lesssim \left[H_0 H_1 a + \frac{H_1}{b}\right]^{m} \overline{M}^*_m \cdot (a / H_0)^{-n} M_n . 
	\end{align*}
By Lemma \ref{l:RExplicit}, it follows that $\mathscr{R}_M < \infty$.

Now, suppose $f \in C^\infty(\R)$ is nontrivial and such that, for certain $h, k > 0$,
	\[ |x^m f^{(n)}(x)| \lesssim k^m  \overline{M}^*_m \cdot h^n M_n , \qquad m, n \in \N . \]
We have, by \eqref{eq:SeqFromAssocFunc}, for $t \geq 0$,
	\[ e^{\omega_{\overline{M}}^*(t)} = \sup_{s \geq 0} e^{s t - \omega_{\overline{M}}(s)} \leq \sum_{n \in \N} \frac{t^n}{n!} \sup_{s \geq 0}  s^n e^{-\omega_{\overline{M}}(s)} = \sum_{n \in \N} \frac{t^n}{\overline{M}^*_n} \leq 2 e^{\omega_{\overline{M}^*}(2 t)} . \]
Therefore, for $n \in \N$ and $x \neq 0$,
	\[ |f^{(n)}(x)| \lesssim \left( \inf_{m \in \N} \frac{k^m}{|x|^m} \overline{M}^*_m \right) \cdot h^n M_n = e^{-\omega_{\overline{M}^*}\left(\frac{|x|}{k}\right)} \cdot h^n M_n \leq 2 e^{-\omega_{\overline{M}}^*\left(\frac{|x|}{2k}\right)} \cdot h^n M_n . \]
Since $\varpi_M \asymp \omega_{\overline{M}}$ by Lemma \ref{l:overlineM}, we then find some $K > 0$ such that
	\[ |f^{(n)}(x)| \lesssim e^{-\varpi^*_M\left( \frac{|x|}{2 k K} \right)} \cdot h^n M_n , \qquad n \in \N , ~ x \in \R . \]
Lemma \ref{l:GSToEntire}$(i)$ with $\eta = \varpi_M$ shows that $\widehat{f}$ extends to an entire function satisfying, for every $k' > k$,
	\[ |\widehat{f}(z)| \lesssim \exp\left[ \varpi_M\left(2 k' K |y|\right) - \omega_M\left(\frac{|x|}{h}\right) \right] , \qquad z \in \C . \]
Then Lemma \ref{l:PsiExplicit} tells us that necessarily $1 / (2 K k h) \leq \Psi_{(\omega_M, \varpi_M)}$.
Therefore, by Lemma \ref{l:RExplicit}, 
	\[ \mathscr{R}_M \geq \frac{1}{2 K \Psi_{(\omega_M, \varpi_M)}} > 0 . \qedhere \]
\end{proof}

\begin{proof}[Proof of Theorem \ref{t:UncertaintyGSOriginalM}]
1. By Theorem \ref{t:UncertaintyGS}, using that $\overline{M} \asymp M$ by Proposition \ref{p:CondEquivSeq}.

2. Proposition \ref{p:CondEquivSeq} shows that $\overline{M} \triangleleft M$ and thus $M^* \triangleleft \overline{M}^*$.
Therefore, any nontrivial $f \in C^\infty(\R^d)$ satisfying \eqref{eq:GSCondOriginalM} for certain $k, h > 0$ would also have, for any $\varepsilon > 0$,
	\[ \sup_{x \in \R^d} |x^\alpha f^{(\beta)}(x)| \lesssim \varepsilon^{\alpha} \overline{M}^*_{|\alpha|} \cdot h^\beta M_{|\beta|} , \qquad \alpha, \beta \in \N^d . \]
Taking $\varepsilon < \mathscr{R}_M / h$, this contradicts Theorem \ref{t:UncertaintyGS}.
Hence, no such $f$ can exist.
\end{proof}

\section{The nontriviality of the Gelfand-Shilov spaces}
\label{sec:NonTrivGS}

We consider the nontriviality of the Gelfand-Shilov spaces $\S^M_A(\R^d)$ and $\Sigma^M_A(\R^d)$ as defined in \eqref{eq:GSRoumieuDef} and \eqref{eq:GSBeurlingDef}.
In particular, we will show Theorems \ref{t:CharNonTrivGSIntro} and \ref{t:CharNonTrivGSIntro2} from the Introduction \ref{sec:Intro}.
When imposing the weakest assumptions on the weight sequences, we may give sufficient conditions for triviality and nontriviality.
Theorem \ref{t:CharNonTrivGSIntro} will then be a refinement of it.

\begin{theorem}
\label{t:GeneralCharNonTrivGS}
Let $M$ be a weight sequence satisfying $(M.2)$ and $(M.2)^*$, and let $A$ be a weight sequence satisfying $(M.2)$.
\begin{itemize}
	\item[1.] If $p! \preceq \overline{M} \cdot A$, then
		\[ \S^M_A(\R^d) \neq \{0\} \neq \S^A_M(\R^d) ; \]
	If $\overline{M} \cdot A \triangleleft p!$, then
		\[ \S^M_A(\R^d) = \{0\} = \S^A_M(\R^d) . \]
		
	\item[2.] If $p! \triangleleft \overline{M} \cdot A$, then
		\[ \Sigma^M_A(\R^d) \neq \{0\} \neq \Sigma^A_M(\R^d) ; \]
	If $\overline{M} \cdot A \preceq p!$, then
		\[ \Sigma^M_A(\R^d) = \{0\} = \Sigma^A_M(\R^d) . \]
\end{itemize}
\end{theorem}

The following result shows when condition \eqref{eq:GSNonTrivNecc} becomes a characterization.
In particular, it covers Theorem \ref{t:CharNonTrivGSIntro2}.
Note that when $M$ is geometrically persistent and $A$ does not satisfy $(M.2)^*$, sufficiency in the following theorem can fail; see Theorem \ref{t:UncertaintyGSOriginalM}.
It is, however, always necessary; see Proposition \ref{p:NonTrivGSNeccCond} below.

\begin{theorem}
\label{t:GSNonTrivCharClassicalCond}
Let $M$ be a weight sequence satisfying $(M.2)$ and $(M.2)^*$, and let $A$ be a weight sequence satisfying $(M.2)$.
Suppose that either $\alpha(M) < 1$ or $A$ satisfies $(M.2)^*$.
Then
	\[ \S^M_A(\R^d) \neq \{0\} \neq \S^A_M(\R^d) \quad \Longleftrightarrow \quad p! \preceq M \cdot A  , \]
and
	\[ \Sigma^M_A(\R^d) \neq \{0\} \neq \Sigma^A_M(\R^d) \quad \Longleftrightarrow \quad p! \triangleleft M \cdot A . \]
\end{theorem}

To prove Theorems \ref{t:CharNonTrivGSIntro}, \ref{t:GeneralCharNonTrivGS}, and \ref{t:GSNonTrivCharClassicalCond}, we will make use of Theorem \ref{t:UncertaintyGS}.
Before doing so, we need several results in preparation.
First, we observe that the conditions on the weight sequences as in Theorem \ref{t:GSNonTrivCharClassicalCond} are always necessary, even under considerably less restrictive assumptions.

\begin{proposition}
\label{p:NonTrivGSNeccCond}
Let $M$ and $A$ be weight sequences satisfying $(M.2)$.
Then
	\[ \S^M_A(\R^d) \neq \{0\} \quad \Longrightarrow \quad p! \preceq M \cdot A \]
and
	\[ \Sigma^M_A(\R^d) \neq \{0\} \quad \Longrightarrow \quad p! \triangleleft M \cdot A . \]
\end{proposition}

\begin{proof}
As before, we may assume $d = 1$.
Suppose $f \in C^1(\R)$ is such that $x^p f(x) \in L^2(\R)$ and $x^{p+1} f'(x) \in L^2(\R)$ for some $p \in \N$.
Then
	\[ \int_{-\infty}^\infty (x^{2p + 1} |f(x)|^2)^{\prime} \, \mathrm{d}x = 0 , \]
and so
	\[ (2p + 1) \| x^p f \|_{L^2}^2 = - 2 \re \int_{-\infty}^\infty x^{2p + 1} f'(x) \overline{f}(x) \, \mathrm{d}x . \]
The Cauchy-Schwarz inequality yields
	\begin{equation}
		\label{eq:L2Est} 
		\frac{2p + 1}{2} \| x^p f \|_{L^2} \leq \| x^{p+1} f' \|_{L^2} .   
	\end{equation}
Suppose now $f \in \S^M_A(\R)$, respectively $f \in \Sigma^M_A(\R)$, is nontrivial.
By \cite[Theorem 1.2]{D-N-V-NuclGSSpaces}, we have that for some, respectively all, $h > 0$, 
	\[ \| x^p f^{(p)} \|_{L^2} \lesssim h^p M_p A_p , \qquad p \in \N . \]
Iterating \eqref{eq:L2Est}, we see
	\[ \frac{(2p - 1) (2p - 3) \cdots 1}{2^p} \| f \|_{L^2} \leq \| x^p f^{(p)} \|_{L^2} \lesssim h^p M_p A_p , \qquad p \in \Z_+ . \]
As
	\[ \frac{(2p - 1) (2p - 3) \cdots 1}{2^p} = \frac{(2p)!}{4^p p!} = \frac{p!}{4^p} {2p \choose p} \geq \frac{p!}{4^p} , \]
we find
	\[ p! \lesssim (4h)^p M_p A_p , \qquad p \in \N . \qedhere \]
\end{proof}

Next, we consider the Fourier transform on the Gelfand-Shilov spaces.
The following result is well-known; we need it in the following specific form.

\begin{lemma}
\label{l:GSFourierInvariant}
Let $M$ and $A$ be weight sequences with $A$ satisfying $(M.2)$.
Suppose $p! \lesssim \kappa^{2p} M_p A_p$ for some $\kappa > 0$.
There exists an $L > 0$ such that if $f \in C^\infty(\R^d)$ satisfies, for some $k, h > 0$,
	\begin{equation} 
		\label{eq:GSBeurlingStep}
		\sup_{x \in \R^d} |x^\alpha f^{(\beta)}(x)| \lesssim k^{\alpha} A_{|\alpha|} \cdot h^{\beta} M_{|\beta|} , \qquad \alpha, \beta \in \N^d ,
	\end{equation}
then
	\begin{equation*} 
		\label{eq:GSBeurlingStepFourier}
		\sup_{x \in \R^d} |x^\alpha \widehat{f}^{(\beta)}(x)| \lesssim (L \max\{h, \kappa\})^{\alpha} M_{|\alpha|} \cdot (L \max\{k, \kappa\})^{\beta} A_{|\beta|} , \quad \alpha, \beta \in \N^d . 
	\end{equation*}
In particular, if both sequences satisfy $(M.2)$, the Fourier transform is a bijection from $\S^M_A(\R^d)$ to $\S^A_M(\R^d)$.
If additionally $p! \triangleleft \, M \cdot A$, then the Fourier transform is a bijection from $\Sigma^M_A(\R^d)$ to $\Sigma^A_M(\R^d)$.
\end{lemma}

\begin{proof}
Since $A$ satisfies $(M.2)$, there exists $H \geq 1$ such that $A_{p+q} \leq H^{p+q} A_p A_q$ for any $p, q \in \N$.
Suppose $f \in C^\infty(\R^d)$ satisfies \eqref{eq:GSBeurlingStep} for some $k, h > 0$.
Then, for any $\alpha, \beta \in \N^d$,
	\[
		\int_{\R^d} |x^\alpha f^{(\beta)}(x)| \, \mathrm{d}x
		= \int_{\R^d} \prod_{j=1}^d (1 + x_j^2) |x^\alpha f^{(\beta)}(x)| \, \frac{\mathrm{d} x}{\prod_{j=1}^d (1 + x_j^2)} 
		\lesssim (H k)^\alpha A_{|\alpha|} \cdot h^\beta M_{|\beta|} .
	\]
We find, for any $\alpha, \beta \in \N^d$,
	\begin{align*}
		|\xi^\alpha \widehat{f}^{(\beta)}(\xi)|
		&= \left| \xi^\alpha \int_{\R^d} (- i x)^\beta f(x) e^{- i \ev{x}{\xi}} \, \mathrm{d}x \right| \\
		&\leq \sum_{\gamma \leq \alpha, \beta} {\alpha \choose \gamma} \frac{\beta!}{(\beta - \gamma)!} \int_{\R^d} |x^{\beta - \gamma} f^{(\alpha - \gamma)}(x)| \, \mathrm{d}x \\
		&\lesssim \sum_{\gamma \leq \alpha, \beta} {\alpha \choose \gamma} {\beta \choose \gamma} \gamma! \cdot h^{\alpha - \gamma} M_{|\alpha - \gamma|} \cdot (H k)^{\beta - \gamma} A_{|\beta - \gamma|} \\
		&\lesssim \sum_{\gamma \leq \alpha, \beta} {\alpha \choose \gamma} {\beta \choose \gamma}  \kappa^\gamma h^{\alpha - \gamma} M_{|\gamma|} M_{|\alpha - \gamma|} \cdot \kappa^\gamma (H k)^{\beta - \gamma} A_{|\gamma|} A_{|\beta - \gamma|} \\
		&\lesssim (2 \max\{h, \kappa\})^{\alpha} M_{|\alpha|} \cdot (2 H \max\{k, \kappa\})^{\beta} A_{|\beta|} .
	\end{align*}
The last statement is now clear.
\end{proof}

Our next observation is that $\Sigma^M_A(\R^d) \neq \{0\}$ already if for each $k, h > 0$ there exists a nontrivial $f \in C^\infty(\R^d)$ satisfying \eqref{eq:GSBeurlingStep}.

\begin{lemma}
\label{l:NonTrivBySteps}
Let $M$ and $A$ be weight sequences satisfying $(M.2)$.
The space $\Sigma^M_A(\R^d)$ is nontrivial if and only if for every $k, h > 0$ there exists $f \in C^\infty(\R^d) \setminus \{0\}$ satisfying \eqref{eq:GSBeurlingStep}.
\end{lemma}

\begin{proof}
Necessity is immediate.
Let us now show sufficiency, and therefore suppose we can find for any $k, h > 0$ a function $f_{k, h} \in C^\infty(\R^d) \setminus \{0\}$ satisfying \eqref{eq:GSBeurlingStep}.
Let $X_h$, $h > 0$, denote the Banach space of all $f \in C^\infty(\R^d)$ satisfying \eqref{eq:GSBeurlingStep} for $h$ in both parameters, endowed with the norm 
	\[ \|f\|_{X_h} = \sup_{\alpha, \beta \in \N^d} \frac{\|x^\alpha f^{(\beta)}(x)\|_{L^\infty}}{h^{\alpha + \beta} A_{|\alpha|} M_{|\beta|}} . \]
Since $f_{h, h} \in X_h$, each space is nontrivial.
As $M$ and $A$ satisfy $(M.2)$, there is $C, H \geq 1$ such that $M_{p + q} \leq C H^{p+q} M_p M_q$ and $A_{p + q} \leq C H^{p+q} A_p A_q$ for any $p, q \in \N$.
We claim that
	\begin{equation}
		\label{eq:DenseEmbedding} 
		X_{h} \subseteq \overline{X_{h/2}}^{X_{4 \sqrt{d} H h}} , \qquad h > 0 . 
	\end{equation}
Before proving this is true, let us first show how it yields the result.
We define the Banach spaces
	\[ Y_n = \overline{X_{2^{-(n+1)}}}^{X_{\sqrt{d}H 2^{-(n - 2)}}} , \qquad n \geq 1 . \]
By \eqref{eq:DenseEmbedding} we have that $X_{2^{-n}} \subseteq Y_n$ for any $n \geq 1$, so that each embedding $Y_{n + 1} \subseteq Y_n$ is dense.
By the Mittag-Leffler procedure (see, e.g., \cite{D-ProjIndLimBanach}), the intersection of the sequence $(Y_n)_{n \geq 1}$ is then nontrivial.
As this intersection coincides with $\Sigma^M_A(\R^d)$, the proof would be complete.

We now verify \eqref{eq:DenseEmbedding}.
Consider, for $h > 0$, the Banach spaces
	\[ \mathscr{X}_h = \{ \Phi \in L^1(\R^{2d}) \mid \|\Phi\|_{\mathscr{X}_h} = \iint_{\R^{2d}} |\Phi(t, \xi)| e^{\omega_A\left(\frac{|t|}{h}\right) + \omega_M\left(\frac{|\xi|}{h}\right)} \, \mathrm{d}t \mathrm{d}\xi < \infty \} . \]
If $\psi \in C^\infty(\R^d)$ is such that, for some $h, C_\psi > 0$ and any $\alpha, \beta \in \N^d$,
	\begin{equation}
		\label{eq:PsiLinftyL1} 
		|x^\alpha \psi^{(\beta)}(x)| \leq C_\psi h^{\alpha + \beta} A_{|\alpha|} M_{|\beta|} \quad \text{and} \quad \int_{\R^d} |x^\alpha \psi^{(\beta)}(x)| \, \mathrm{d}x \leq C_\psi h^{\alpha + \beta} A_{|\alpha|} M_{|\beta|} , 
	\end{equation}
then we have continuity of the short-time Fourier transform
	\begin{equation}
		\label{eq:STFT} 
		V_\psi : X_h \to \mathscr{X}_{2 \sqrt{d} H h} , \quad f \mapsto \left( (t, \xi) \mapsto V_\psi f (t, \xi) := \int_{\R^d} f(x) \overline{\psi(x - t)} e^{- i \ev{x}{\xi}} \, \mathrm{d}x \right) ,
	\end{equation}
and the adjoint short-time Fourier transform
	\begin{equation}
		\label{eq:adjointSTFT}
		V^*_\psi : \mathscr{X}_{h} \to X_{2 h} , ~ \Phi \mapsto \left( x \mapsto V^*_\psi \Phi(x) :=\iint_{\R^{2d}} \Phi(t, \xi) \psi(x - t) e^{i \ev{x}{\xi}} \, \mathrm{d}t \mathrm{d}\xi \right) . 
	\end{equation}
Indeed, for any $f \in X_h$ we have, for $\alpha, \beta \in \N^d$,	
	\begin{align*}
		| t^{\alpha} \xi^{\beta} V_\psi f(t, \xi)| 
		&\leq \sum_{\alpha' \leq \alpha} \sum_{\beta' \leq \beta} {\alpha \choose \alpha'} {\beta \choose \beta'} \int_{\R^{d}} |x^{\alpha'} f^{(\beta')}(x)| \cdot |(x - t)^{\alpha - \alpha'} \psi^{(\beta - \beta')}(x - t)| \, \mathrm{d} x \\
		&\lesssim (2 h)^{\alpha + \beta} A_{|\alpha|} M_{|\beta|} \|f\|_{X_h} ,
	\end{align*}
so that
	\[ |V_\psi f(t, \xi)| \lesssim \|f\|_{X_h} \exp\left[-\omega_A\left(\frac{|t|}{2\sqrt{d}h}\right) - \omega_M\left(\frac{|\xi|}{2\sqrt{d}h}\right) \right] . \]
As $M$ and $A$ satisfy $(M.2)$ we may use lower moderateness \cite[Proposition 3.6, p.~51]{K-UltradistI} to find
	\[
		\| V_\psi f \|_{\mathscr{X}_{2 \sqrt{d} H h}} 
		= \iint_{\R^{2d}} |V_\psi f(t, \xi)| e^{\omega_A\left(\frac{|t|}{2 \sqrt{d} H h}\right) + \omega_M\left(\frac{|\xi|}{2 \sqrt{d} H h}\right)} \, \mathrm{d}t \mathrm{d}\xi 
		\lesssim \|f\|_{X_h} ,
	\]
which shows the continuity of \eqref{eq:STFT}.
On the other hand, for any $\Phi \in \mathscr{X}_h$ and $\alpha, \beta \in \N^d$,
	\begin{align*}
		&\left|x^\alpha \frac{\partial^\beta}{\partial x^\beta} V^*_\psi \Phi(x)\right| \\
		&\qquad \leq \sum_{\alpha' \leq \alpha} \sum_{\beta' \leq \beta} {\alpha \choose \alpha'} {\beta \choose \beta'} \iint_{\R^{2d}} |t^{\alpha'} \xi^{\beta'} \Phi(t, \xi)| |(x - t)^{\alpha - \alpha'} \psi^{(\beta - \beta')}(x - t)| \, \mathrm{d}t \mathrm{d}\xi \\
		&\qquad \lesssim (2 h)^{\alpha + \beta} A_{|\alpha|} M_{|\beta|} \|\Phi\|_{\mathscr{X}_h} ,
	\end{align*}
proving the continuity of \eqref{eq:adjointSTFT}.
We can now show \eqref{eq:DenseEmbedding}.
Fix $h > 0$ and take any $\psi \in C^\infty(\R^d) \setminus \{0\}$ satisfying \eqref{eq:PsiLinftyL1} for $h / 4$, which is possible by our assumption.
Since any $\Phi \in C_c(\R^{2d})$ belongs to $\mathscr{X}_{h / 4}$, it follows from \eqref{eq:adjointSTFT} that $\{ V^*_\psi \Phi \mid \Phi \in C_c(\R^{2d}) \} \subseteq X_{h/2}$.
Take any $f \in X_h$, then $V_\psi f \in \mathscr{X}_{2\sqrt{d}Hh}$ by \eqref{eq:STFT}.
Since $C_c(\R^{2d})$ is dense in $\mathscr{X}_{2\sqrt{d}Hh}$, there is a sequence $\Phi_n \in C_c(\R^{2d})$ such that $\Phi_n \to V_\psi f$ in $\mathscr{X}_{2\sqrt{d}Hh}$.
Then $V^*_\psi \Phi_n \to V^*_\psi \circ V_\psi f$ in $X_{4\sqrt{d}Hh}$ by \eqref{eq:adjointSTFT}.
By \cite[Corollary 3.2.3, p.~44]{G-FoundTFAnal}
	\[ V^*_\psi \circ V_\psi f = (2\pi)^d (\psi, \psi)_{L^2} f , \]
which proves \eqref{eq:DenseEmbedding}.
\end{proof}

We may now show Theorem \ref{t:GeneralCharNonTrivGS}.

\begin{proof}[Proof of Theorem \ref{t:GeneralCharNonTrivGS}]
If $M$ does not satisfy $(\QA)$, then $\overline{M} \asymp p!$.
Therefore, automatically $p! \triangleleft \overline{M} \cdot A$, while all spaces are nontrivial \cite[Theorem 4.2, p.~56]{K-UltradistI}.
Hence, in this case, the result holds.
We may now assume $M$ satisfies $(\QA)$.

Suppose there exists some $h_0 > 0$ such that $p! \lesssim h_0^p \overline{M}_p A_p$ for $p \in \N$.
Note that then also $p! \lesssim h_0^p M_p A_p$ as $\overline{M}_p \leq M_p$ for all $p \in \N$.
We have $\overline{M}^*_p \lesssim h_0^p A_p$ for $p \in \N$.
Take any $h > 0$ and any $k > \mathscr{R}_M / h$, with $\mathscr{R}_M$ as in Theorem \ref{t:UncertaintyGS}.
Then, there exists some $f \in C^\infty(\R^d) \setminus \{0\}$ such that
	\[ \sup_{x \in \R^d} |x^\alpha f^{(\beta)}(x)| \lesssim k^{\alpha} \overline{M}^*_{|\alpha|} \cdot h^{\beta} M_{|\beta|} \lesssim (h_0 k)^{\alpha} A_{|\alpha|} \cdot h^{\beta} M_{|\beta|} . \]
Therefore $f \in \S^M_A(\R^d)$, so $\S^M_A(\R^d) \neq \{0\}$, and also $\S^A_M(\R^d) \neq \{0\}$ by Lemma \ref{l:GSFourierInvariant}.
Moreover, if $h_0$ can be taken arbitrarily small, so can $h_0 k$, and in this case Lemma \ref{l:NonTrivBySteps} shows that $\Sigma^M_A(\R^d) \neq \{0\}$. 
Then also $\Sigma^A_M(\R^d) \neq \{0\}$ by Lemma \ref{l:GSFourierInvariant}. 
This shows the first statements in (1) and (2).

Suppose now that for some $h_0 > 0$ we have $\overline{M}_p A_p \lesssim h_0^p p!$ for $p \in \N$.
Then $A_p \lesssim h_0^p \overline{M}^*_p$ for $p \in \N$.
Suppose that $f \in C^\infty(\R^d) \setminus \{0\}$ is such that
	\[ \sup_{x \in \R^d}  |x^\alpha f^{(\beta)}(x)| \lesssim k^\alpha A_{|\alpha|} \cdot h^{\beta} M_{|\beta|} , \qquad \alpha, \beta \in \N^d . \]
Then also
	\[ \sup_{x \in \R^d} |x^\alpha f^{(\beta)}(x)| \lesssim (h_0 k)^{\alpha} \overline{M}^*_{|\alpha|} \cdot h^{\beta} M_{|\beta|} , \qquad \alpha, \beta \in \N^d . \]
By Theorem \ref{t:UncertaintyGS}, we must necessarily have that $h_0 h k \geq \mathscr{R}_{M}$.
This leads to a contradiction if we can choose $h_0$, $h$, or $k$ to be arbitrarily small.
This shows the second statements in (1) and (2).
\end{proof}

Continuing to the proof of Theorem \ref{t:CharNonTrivGSIntro}, we first show the following crucial observation on the conditions for nontriviality we have encountered so far.

\begin{lemma}
\label{l:NonTrivCondSymmetric}
Let $M$ and $A$ be weight sequences such that $M$ satisfies $(M.2)^*$ or $\alpha(A) < 1$.
Then
	\[ 
		p! \preceq M \cdot A ~ \Longrightarrow ~ p! \preceq M \cdot \overline{A} 
		\qquad \text{and} \qquad
		p! \triangleleft M \cdot A ~ \Longrightarrow ~ p! \triangleleft M \cdot \overline{A}
		.
	\]	
\end{lemma} 

\begin{proof}
When $\alpha(A) < 1$, then $\overline{A} \asymp A$ by Proposition \ref{p:CondEquivSeq}, and the conclusion is direct.
Let us therefore now assume $M$ satisfies $(M.2)^*$.
We may assume $M$ and $A$ are normalized.
By potentially going to a higher coefficient, suppose $p! \leq h^{p} M_p A_p$ for certain $h > 0$ and $p \geq p_h$.
Then, by Stirling's approximation, for $p \geq p_h$,
	\[ m_p a_p \geq \left( M_p A_p \right)^{\frac{1}{p}} \geq \frac{1}{h} (p!)^{\frac{1}{p}} \geq \frac{p}{e h}  ,  \]
that is
	\begin{equation}
		\label{eq:UpperBound1/a_p} 
		\frac{1}{a_p} \leq (e h) \frac{m_p}{p}  . 
	\end{equation}
If $\overline{a}(t)$ is the counting function of $(\overline{a}_p)_{p \in \Z_+}$, then, for $t \geq a_{p_h}$,
	\[ 
		\overline{a}(t) 
		= \sum_{a_p \leq t} \left\lfloor \frac{t}{a_p} \right\rfloor \leq t \sum_{a_p \leq t} \frac{1}{a_p} 
		\leq t \left( \sum_{p = 1}^{p_h - 1} \frac{1}{a_p} + \sum_{\substack{a_p \leq t \\ p \geq p_h}} \frac{1}{a_p} \right)  
		\leq t \left( \sum_{p = 1}^{p_h - 1} \frac{1}{a_p} \right) +  (e h) t \sum_{a_p \leq t} \frac{m_p}{p}  . 
	\]
Now suppose $Q \geq 2$ is such that $2 m_{p} \leq m_{Q p}$ for $p \geq p_0 \geq 1$.
Then, for $t < a_{q}$ with $q \geq p_0$,
	\begin{align*}
		\sum_{a_p \leq t} \frac{m_p}{p} 
		&\leq \sum_{p = 1}^{p_0 - 1} \frac{m_p}{p} + \sum_{k = 0}^{\lceil \log_Q q \rceil} \sum_{\substack{\frac{q}{Q^{k+1}} < p \leq \frac{q}{Q^k} \\ p_0 \leq p}} \frac{m_p}{p} \\
		&\leq \left(\sum_{p=1}^{p_0 - 1} \frac{1}{p} \right) m_{p_0}  + m_q \sum_{k = 0}^{\lceil \log_Q q \rceil} \frac{Q + 1}{2^{k}} 
		\leq C m_q , 
	\end{align*}
for some $C \geq 1$ independent of $q$.
There is some $p'_h \geq \max\{p_h, p_0\}$, such that for any $a_{p'_h} \leq t < a_q$, with $q > p'_h$,
	\[ \overline{a}(t) \leq (Ce + 1)h m_q t . \]
Put
	\[ t_p = \frac{1}{2(Ce + 1)h} \frac{p}{m_p} , \qquad p > p'_h , \]
then $t_p < a_p$.
If also $t_p \geq a_{p'_h}$, then the preceding estimate with $q=p$ gives
	\[ \overline{a}(t_p) \leq \frac{p}{2} . \]
Consequently, for sufficiently large $p$, we have $\overline{a}(t_p) \leq p/2$, and hence
	\[ \overline{a}_p \geq t_p = \frac{1}{2(Ce + 1)h} \frac{p}{m_p}  . \]
Multiplying, we find for $k = 2 (Ce + 1) h$  some $C_k > 0$ for which
	\[ p! \leq C_k k^p M_p \overline{A}_p , \qquad p \in \N . \]
As $k$ becomes arbitrarily small when $h$ does, this finishes the proof.
\end{proof}

\begin{proof}[Proof of Theorem \ref{t:CharNonTrivGSIntro}]
By Proposition \ref{p:NonTrivGSNeccCond} and Lemma \ref{l:GSFourierInvariant}, we may assume without loss of generality that $M$ satisfies $(M.2)$ and $(M.2)^*$, while $A$ satisfies $(M.1)^*$.
Also, as in the proof of Theorem \ref{t:GeneralCharNonTrivGS}, we only have to consider the case where $M$ satisfies $(\QA)$.
The sufficiency of the conditions in \eqref{eq:NiceCharNonTrivRoumieu} and \eqref{eq:NiceCharNonTrivBeurling} then follows from Theorem \ref{t:GeneralCharNonTrivGS}, so it remains to show necessity.
As before, we may assume without loss of generality that $d = 1$.

Let $f \in C^\infty(\R) \setminus \{0\}$ satisfy \eqref{eq:GSBeurlingStep} for certain $k, h > 0$.
Then, for any $x \neq 0$ and $n \in \N$,
	\[ |f^{(n)}(x)| \lesssim \left( \inf_{m \in \N} \frac{k^m A_m}{|x|^m} \right) \cdot h^n M_n = e^{-\omega_A\left(\frac{|x|}{k}\right)} \cdot h^n M_n . \]
As $A$ satisfies $(M.1)^*$, there is a convex weight function $\eta$ such that $\eta^* \asymp \omega_A$ \cite[Theorem 4.15]{S-ConjWeight}.
Therefore, for some $K > 0$, we have
	\[ |e^{\eta^*\left(\frac{|x|}{K k}\right)} f^{(n)}(x)| \lesssim h^n M_n , \qquad n \in \N . \]
By Lemma \ref{l:GSToEntire}$(i)$, we now have that $\widehat{f}$ may be extended to an entire function satisfying, for every $k' > k$,
	\[ |\widehat{f}(z)| \lesssim \exp\left[ \eta(K k' |y|) - \omega_M\left(\frac{|x|}{h}\right) \right] , \qquad z = x+ iy \in \C . \]
Applying Proposition \ref{p:NonTrivEntireLowerBound} gives, for $t$ large enough,
	\[ \varpi_M(t) \leq \eta(\pi K k' h t) . \]
So, for $t$ large enough,
	\[ \eta^*(t) \leq \varpi_M^*\left(\pi K k' h t\right) . \]
Combining this with \eqref{eq:SeqFromAssocFunc} and Lemma \ref{l:ConjOverlineMAlmostWeightSeq}, we find $H > 0$ (independent of $k$ and $h$) and $C > 0$ such that
	\[ \overline{M}^*_p \leq C (H kh)^p A_p , \qquad p \in \N , \]
or
	\[ p! \leq C (H kh)^p \overline{M}_p A_p , \qquad p \in \N . \]
This means $p! \preceq \overline{M} \cdot A$.
When $h$ or $k$ can be chosen arbitrarily small, we see that $p! \triangleleft \overline{M} \cdot A$.
Lemma \ref{l:NonTrivCondSymmetric} now completes the proof.
\end{proof}

We end this section with the proof of Theorem \ref{t:GSNonTrivCharClassicalCond}, and thus also Theorem \ref{t:CharNonTrivGSIntro2}.

\begin{proof}[Proof of Theorem \ref{t:GSNonTrivCharClassicalCond}]
Necessity of the conditions follows from Proposition \ref{p:NonTrivGSNeccCond}.
Suppose now $p! \preceq M \cdot A$; the case of $p! \triangleleft \, M \cdot A$ can be done analogously.
As either $\alpha(M) < 1$ or $A$ satisfies $(M.2)^*$, Lemma \ref{l:NonTrivCondSymmetric} shows that $p! \preceq \overline{M} \cdot A$.
Then Theorem \ref{t:GeneralCharNonTrivGS} yields the nontriviality of $\S^M_A(\R^d)$ and $\S^A_M(\R^d)$.
\end{proof}

\end{document}